\def\VERSIONNUMBER{1.4.2t}
\def\VERSION{Version \VERSIONNUMBER}
\def\DATE{February 1, 1997}
\def\EPRINT{{\tt q-alg/9612033}}
\documentstyle[amssymb,amscd]{amsart}

\newcommand\ad{\operatorname{ad}}

\newcommand\bbbP{{\Bbb P}}

\newcommand\bcdot{{\raise0.2ex \hbox{\bf.}}}
\newcommand\bdot{{\raise0.4ex \hbox{\bf.}}}
\newcommand\C{{\Bbb C}}
\newcommand\caff{{{\mathrm c}_{\mathrm a}}}
\newcommand\cvir{{{\mathrm c}_{\mathrm V}}}
\newcommand\CB{\operatorname{CB}}

\newcommand\chat{{\hat c}} 
\newcommand\Comp{{\Bbb C}}
\newcommand\CC{\operatorname{CC}}
\newcommand\D{{\cal D}}
\newcommand\END{\mathop{{\cal E}\text{\it nd}}\nolimits}

\newcommand\eps{\varepsilon}
\newcommand\F{{\cal F}}

\newcommand\g{{\frak g}}
\newcommand\Gammatilde{{\widetilde\Gamma}}
\newcommand\ghat{\hat\g}
\newcommand\gtw{\g^{{\mathrm t}{\mathrm w}}}
\newcommand\gDpr{\g^D_{\prin}}    
\newcommand\gDPr{\g^D_{\Prin}} 

\newcommand\Gtw{G^{{\mathrm t}{\mathrm w}}}

\newcommand\Heis{{\cal H}}
\newcommand\HOM{\mathop{{\cal H}{\mathit o}{\mathit m}}\nolimits}
\newcommand\Hom{\operatorname{Hom}}
\newcommand\hvee{h^\vee}
\newcommand\id{{{\mathrm i}{\mathrm d}}}
\renewcommand\Im{\operatorname{Im}}
\newcommand\injto{\hookrightarrow}

\newcommand\isoto{\overset\sim\rightarrow}
\newcommand\isoot{\overset\sim\leftarrow}
\newcommand\Ind{\operatorname{Ind}}
\newcommand\khat{{\hat k}} 
\newcommand\K{{\cal K}}

\newcommand\m{{\frak m}}
\newcommand\M{{\cal M}}
\ifx\mathit\undefined
\newcommand\mathit{} 
\fi

\newcommand\twist{{{\mathrm t}{\mathrm w}}}
\newcommand\nontw{{{\mathrm n}{\mathrm t}}}
\newcommand\NP{\hbox{$\circ \atop \circ$}} 
\renewcommand\O{{\cal O}}
\newcommand\od[2]{\frac{d #1}{d #2}}
\newcommand\onto{\twoheadrightarrow}

\ifx\pd\undefined
\newcommand\pd[2]{\frac{\partial #1}{\partial #2}}
\else
\renewcommand\pd[2]{\frac{\partial #1}{\partial #2}} 
\fi
\newcommand\PD[2]{{\partial #1}/{\partial #2}}
\newcommand\pr{\pi}
\newcommand\prin{{\dot X}} 
\newcommand\Prin{{\dot \X}} 
\newcommand\Res{\mathop{\operatorname{Res}}\nolimits}
\newcommand\sCB{{\cal C}{\cal B}}
\newcommand\sCC{{\cal C}{\cal C}}
\renewcommand\setminus{\smallsetminus}
\newcommand\simeqq{\cong}
\newcommand\sWeyl{{\cal M}}
\newcommand\T{{\cal T}}
\newcommand\TDPr{{\cal T}^D_{\Prin}} 
\newcommand\tensor{\otimes}
\newcommand\tr{\operatorname{tr}}

\newcommand\UHP{{\frak H}} 

\newcommand\Vir{{\cal V}{\mathit i}{\mathit r}}
\newcommand\VirDPr{\Vir^D_{\Prin}} 
\newcommand\Weyl{M}
\newcommand\X{{\frak X}}
\newcommand\Xtilde{{\widetilde{\frak X}}}
\newcommand\Z{{\Bbb Z}}
\newcommand\Integer{{\Bbb Z}}

\newcommand\sn{\operatorname{sn}}
\newcommand\cn{\operatorname{cn}}
\newcommand\dn{\operatorname{dn}}

\newtheorem{thm}{Theorem}[section]
\newtheorem{lem}[thm]{Lemma}
\newtheorem{prop}[thm]{Proposition}
\newtheorem{cor}[thm]{Corollary}

\theoremstyle{definition}
\newtheorem{defn}[thm]{Definition}

\theoremstyle{remark}
\newtheorem{example}[thm]{Example}
\newtheorem{rem}[thm]{Remark}

\numberwithin{equation}{section}

\newcommand\thmref[1]{Theorem~\ref{#1}}
\newcommand\propref[1]{Proposition~\ref{#1}}
\newcommand\secref[1]{\S\ref{#1}}
\newcommand\corref[1]{Corollary~\ref{#1}}
\newcommand\lemref[1]{Lemma~\ref{#1}}
\newcommand\remref[1]{Remark~\ref{#1}}
\newcommand\exref[1]{Example~\ref{#1}}
\newcommand\tabref[1]{Table~\ref{#1}}
\newcommand\appref[1]{Appendix~\ref{#1}}

\begin{document}

\title[Twisted WZW Models on elliptic curves]
{Twisted Wess-Zumino-Witten models \\ on elliptic curves}


\author{Gen KUROKI} %
\address{Mathematical Institute, Tohoku University, Sendai 980, JAPAN}
\author{Takashi TAKEBE}%
\address{Department of Mathematical Sciences,
the University of Tokyo, Komaba, Tokyo 153, JAPAN}
\curraddr{Department of Mathematics, the University of California,
Berkeley, CA 94720, U.S.A.}

\date{\DATE{}\quad (\VERSION, \EPRINT)}

\ifx\ClassWarning\undefined
\maketitle 
\fi

\begin{abstract}
  Investigated is a variant of the Wess-Zumino-Witten model called a twisted
  WZW model, which is associated to a certain Lie group bundle on a family of
  elliptic curves.  The Lie group bundle is a non-trivial bundle with flat
  connection and related to the classical elliptic $r$-matrix.  (The usual
  (non-twisted) WZW model is associated to a trivial group bundle with trivial
  connection on a family of compact Riemann surfaces and a family of its
  principal bundles.)  The twisted WZW model on a fixed elliptic curve at the
  critical level describes the XYZ Gaudin model.  The elliptic
  Knizhnik-Zamolodchikov equations associated to the classical elliptic
  $r$-matrix appear as flat connections on the sheaves of conformal blocks in
  the twisted WZW model.
\end{abstract}

\ifx\ClassWarning\undefined\else
\maketitle
\fi



\tableofcontents


\setcounter{section}{-1}

\section{Introduction}
\label{intro}

In this paper, we deal with a variant of the (chiral) Wess-Zumino-Witten model
(WZW model, for short) on elliptic curves, which shall be called a twisted WZW
model.  The usual (non-twisted) WZW model on a compact Riemann surface $X$
gives rise to the sheaves of vector spaces (of conformal blocks and of
conformal coinvariants) on any family (or the moduli stack) of principal
$G$-bundles, where $G$ is a semisimple complex algebraic group. Note that the
notion of principal $G$-bundles is equivalent to that of $G^\nontw$-torsors,
where $G^\nontw$ denotes the trivial group bundle $G\times X$ on $X$.  (The
symbol $(\cdot)^\nontw$ stands for ``non-twisted''.)  This suggests that there
exists a model associated to a non-trivial group bundle $\Gtw$ with a flat
connection on $X$, which gives sheaves of conformal blocks and conformal
coinvariants on a family of $\Gtw$-torsors. (The symbol $(\cdot)^\twist$
stands for ``twisted''.) We call such a model a {\em twisted WZW model}
associated to $\Gtw$.

The aim of this work is not to establish a general theory of the twisted
WZW models but to describe certain interesting examples of the twisted WZW
models related to the elliptic classical $r$-matrices (\cite{bel-dr:82}).
In this introduction, we explain our motivations and clarify
the relationship between the twisted WZW models and various problems in
mathematics and physics.

One of the motivations is the viewpoint of representation theory where the
WZW model is formulated as an analogue of a theory of automorphic forms
due to Langlands. We list corresponding ingredients of both theories in
\tabref{table:automorphic-WZW}. Notations shall be explained in the main
part of this paper.

\begin{table}[ht]
\caption{Analogy between automorphic forms and conformal blocks} 
\begin{center}
\begin{tabular}{|p{0.45\textwidth}|p{0.45\textwidth}|}
\hline

Theory of automorphic forms & 
Theory of the WZW models \\ 

\noalign{\hrule height 1.5pt}

a global field, i.e., a number field or the function field of an algebraic
curve over a finite field &
the function field of a compact Riemann surface $X$ \\
\hline

a local field &
$\C((\xi))$, a field of formal Laurent series \\
\hline

a reductive group over the global field &
a semisimple group bundle on $X$ with flat connection 
or the associated Lie algebra bundle $\gtw$ \\
\hline

a non-split reductive group over the global field &
a semisimple group bundle with flat connection on $X$ which is not locally
trivial under the Zariski topology \\ 
\hline

the ad\`ele group associated to the reductive group & 
the affine Lie algebra \hfill\break
$(\g^{\oplus L})^\wedge=\bigoplus_{i=1}^L\g\tensor\C((\xi_i))\oplus\C\khat$ \\
\hline

the principal ad\`ele subgroup of the ad\`ele group &
the subalgebra $\gDpr=H^0(X,\gtw(\ast D))$ of $(\g^{\oplus L})^\wedge$ \\
\hline

a unitary representation of the ad\`ele group  &
a representation $M$ of $(\g^{\oplus L})^\wedge$ 
or its algebraic dual $M^\ast$ \\ 
\hline

the space of automorphic forms in the representation space,
i.e., the invariant subspace of the representation space with respect to the
principal ad\`ele subgroup &
the space $\CB(M)$ of conformal blocks,
i.e., the invariant subspace of $M^\ast$ with respect to $\gDpr$ \\
\hline

\end{tabular}
\end{center}
\label{table:automorphic-WZW}
\end{table}
There is a theory of automorphic forms for arbitrary (possibly non-split)
reductive groups over a global field as well as over a local field.  But
so far only the WZW model associated to the trivial group bundle has been
considered and the counterpart of the non-split reductive group over a
global field has been absent. The twisted WZW models fills this blank.

The second motivation comes from the geometric Langlands program over
$\C$ and its relation with quantum integrable systems. A geometric
analogue of the Langlands correspondence over $\C$ is described by using
the WZW model at the critical level, where the centers of the completed
enveloping algebras of affine Lie algebras are sufficiently large so that we 
can consider analogues of the infinitesimal characters of finite-dimensional
semisimple Lie algebras (\cite{hayashi}, \cite{go-wa}).  For introduction to
the original Langlands program we refer to \cite{borel:79} and
\cite{gelbart:84}.  For a  general formulation of the geometric Langlands
correspondence over $\C$ related to the non-twisted WZW models at the critical
level, see \cite{bei:91} and \cite{bei-dr:94} and for the analogue of the local
Langlands correspondence to affine Lie algebras at the critical level, see
\cite{fei-fre:92} and \cite{fre:91}.  The twisted WZW model at the critical
level shall give a geometric analogue of the Langlands correspondence for
a non-split reductive group over a global field.

To study this model at the critical level is important not only in this
context of the geometric Langlands program but also in the theory of the
quantum integrable spin chains. B.~Feigin, E.~Frenkel and N.~Reshetikhin
found in \cite{f-f-r:94} that the non-twisted WZW model at the critical level
on the Riemann sphere is closely related to a spin chain model called the
Gaudin model. See also \cite{fre:95}. Its Hamiltonian is described as an
insertion of a singular vector of the vacuum representation at a point
and the diagonalization problem turns out to be equivalent to a
description of a certain space of conformal blocks. This ``Gaudin'' model
is, however, merely a special case of the model introduced by M.~Gaudin
\cite{gau:73}, \cite{gau:76}, \cite{gau:83} as a quasi-classical limit of
the XYZ spin chain model. Let us call this general model the XYZ Gaudin
model, following \cite{skl-tak:96} where the diagonalization problem of
this model was studied by the algebraic Bethe Ansatz. In order to extend
the results of \cite{f-f-r:94}, we need the twisted WZW model on a
elliptic curve at the critical level as we shall see in
\secref{xyz-gaudin}. 

We remark that the non-twisted WZW models at the critical level on an elliptic
curve is related to quantum integrable systems on root systems.  In fact those
systems defined by the trigonometric (dynamical) $r$-matrices are described by
the non-twisted WZW model on a degenerated elliptic curve with only one
ordinary double point.  The non-twisted WZW model at the critical level on an
elliptic curve leads to a system called the Gaudin-Calogero model
(\cite{enr-rub:95}, \cite{nek:95}) which was defined as Hitchin's classical
integrable system (\cite{hitch:87}) on the moduli space of semistable
principal bundles on an elliptic curve.%
\footnote{This relation of the non-twisted WZW model and the Gaudin-Calogero
  model is due to B.~Enriquez and A.~Stoyanovsky. TT thanks Enriquez for
  communicating their unpublished result.} %
The reason why root systems appear in the non-twisted WZW models is explained
as follows.  Let $G$ be a complex semisimple group and $T$ its maximal torus.
Let $a$ and $b$ denote generators of the fundamental group $\pi_1(X)$ of an
elliptic curve $X$.  Then, for $g\in T$, the homomorphism from $\pi_1(X)$ into
$G$ sending $a$ and $b$ to $1$ and $g$ respectively induces a semistable
principal $G$-bundle on $X$.  This defines the covering by $T$ of the moduli
space of semistable $G$-bundles on $X$.  Furthermore the universal covering of
$T$ is identified with its Lie algebra, on which the root system structure
exists.  Namely, the root system appears as a covering space of the moduli
space of semistable principal $G$-bundles on $X$.

The third motivation is a geometric interpretation of Etingof's elliptic
KZ equations. As is well-known, the Knizhnik-Zamolodchikov equation is a
system of differential equations satisfied by matrix elements of products
of vertex operators (\cite{kni-zam:84}, \cite{tsu-kan:88}) and is a flat
connection over the family of pointed Riemann spheres.  Similarly from the
non-twisted WZW model over elliptic curves arises the elliptic
Knizhnik-Zamolodchikov-Bernard equations (KZB equations, for short), which
Bernard found in \cite{ber:88-1} by computing traces of products of vertex
operators twisted by $g\in G$.  The interpretation of the elliptic KZB
equations as flat connections on sheaves of conformal blocks, which are
defined without use of the traces, was found in \cite{fel-wie:96}.  Using
the same idea as \cite{ber:88-1}, Etingof computed in \cite{eti:94} a
twisted trace of a product of vertex operators and found that it obeys
linear differential equations of KZ type defined by the elliptic classical
$r$-matrices.  We call these equations the elliptic KZ equations. In the
present paper it is shown that the elliptic KZ equations also has an
interpretation as flat connections on sheaves of conformal blocks.

Let us explain now the contents of this paper.  In \secref{conf-block}, we
give a definition of the conformal coinvariants and the conformal blocks
of the twisted WZW model on an elliptic curve. The definition of the
non-trivial group bundle $\Gtw$ \eqref{def:Gtw} and the associated Lie
algebra bundle $\gtw$ \eqref{def:gtw} is given in \secref{bundles} and
their fundamental properties are studied. This bundle $\gtw$ was 
used by I.~Cherednik \cite{cher:83} for a algebro-geometric interpretation of
classical elliptic $r$-matrices. An important point is that the cohomology
groups of $\gtw$ vanish in all degrees.  Since the 1-cohomology
$H^1(X,\gtw)$ can be canonically identified with the tangent space of the
moduli space of $\Gtw$-torsors at the equivalence class consisting of
trivial ones, the trivial $\Gtw$-torsor can not be deformed.  Thus
non-trivial $\Gtw$-torsors do not appear in the twisted WZW model
associated to $\Gtw$.  Conformal coinvariants and conformal blocks of this
model are defined in \secref{def-conf-block}. We introduce correlation
functions of current and the energy-momentum tensor in \secref{correl},
following mostly \cite{tu-ue-ya:89}. An action of the Virasoro algebra on
the conformal coinvariants and the conformal blocks is defined in
\secref{virasoro-action}. 

This model at the critical level for $G=SL_2(\C)$ describes the XYZ Gaudin
model and the case for $G=SL_N(\C)$ is related to the higher rank
generalizations of the XYZ Gaudin model, as shown in \secref{xyz-gaudin}.

Away from the critical level, we can define a connection on the sheaves of
conformal coinvariants and blocks over the family of pointed elliptic
curves. \secref{sheaf-conf-block} and \secref{sheaf-of-virasoro} are a
sheaf version of \secref{conf-block} over a family of pointed elliptic
curves. By extending the tangent sheaf of the base space of the family
\eqref{eq:LieAlgExtOfTSByTDPr} and constructing its action on the sheaves
of conformal coinvariants and conformal blocks, we can introduce the
$D$-module structure on them in \secref{flat-connections}, which naturally
implies the existence of flat connections. The explicit formulae in
\secref{elliptic-kz} show that our connections are identical with
Etingof's elliptic KZ equations. This connection has modular invariance
which Etingof proved by his explicit expressions of the equations. We give
a geometric proof of this fact in \secref{automorphic}.

Useful properties of theta functions are listed in
\appref{theta-func}.  An algebro-geometric meaning of the extension
\eqref{eq:LieAlgExtOfTSByTDPr} is explained in
\appref{kodaira-spencer}. Higher-genus generalization of the theory is
discussed in \appref{higher-genus}.

\section{Spaces of conformal coinvariants and conformal blocks}
\label{conf-block}

In this section we define the space of conformal coinvariants and
conformal blocks associated to a twisted Group bundle.

\subsection{Group bundles and their associated Lie algebra bundles}
\label{bundles}

In this section we define a group bundle $\Gtw$ and an associated Lie algebra
bundle $\gtw$ on an elliptic curve $X=X_\tau = \C/(\Z+\tau\Z)$, where $\tau$
belongs to the upper half plane $\UHP := \{\Im\tau>0\}$. We fix a global
coordinate $t$ on $X$ which comes from that of $\C$.

Let $G$ be the Lie group $SL_N(\C)$ and $\g$ be its Lie algebra,
\begin{equation*}
    sl_N(\C) =
    \{ A \in M_N(\C) \, | \, \tr A = 0\}.
\end{equation*}
We fix an invariant inner product of $\g$ by
\begin{equation}
    (A | B) := \tr (AB) \quad \text{for $A, B \in \g$}.
\label{def:inner-product}
\end{equation}
Define matrices $\alpha$ and $\beta$ by
\begin{equation}
    \alpha :=  \begin{pmatrix}
         0 & 1 &        & 0 \\
           & 0 & \ddots &   \\
           &   & \ddots & 1 \\
         1 &   &        & 0
         \end{pmatrix},
    \quad
    \beta := \begin{pmatrix}
         1 &       &        & 0 \\
           & \eps  &        &   \\
           &       & \ddots &   \\
         0 &       &        & \eps^{N-1}
         \end{pmatrix},
\label{def:alpha-beta}
\end{equation}
where $\eps = \exp (2\pi i/ N)$.  Then we have $\alpha^N=\beta^N=1$
and $\alpha\beta=\eps\beta\alpha$.

We define the group bundle $\Gtw$ and its associated Lie algebra bundle
$\gtw$ by
\begin{align}
    \Gtw :=& (\C \times G)/{\sim},
\label{def:Gtw}
\\
    \gtw :=& (\C \times \g)/{\approx},
\label{def:gtw}
\end{align}
where the equivalence relations $\sim$ and $\approx$ are defined by
\begin{align}
   (t, g) &\sim (t+1, \alpha g \alpha^{-1}) 
          \sim (t+\tau, \beta g \beta^{-1}),
\label{def:eq-rel-Gtw}
\\
   (t, A) &\approx (t+1, \alpha A \alpha^{-1}) 
          \approx (t+\tau, \beta A \beta^{-1}).
\label{def:eq-rel-gtw}
\end{align}
(Because of $\alpha\beta = \eps\beta\alpha$, the group bundle $\Gtw$ is 
{\em not} a principal bundle.)  The fibers of $\Gtw$ are isomorphic to $G$ and
those of $\gtw$ are isomorphic to $\g$, but there are not canonical
isomorphisms.

The twisted Lie algebra bundle $\gtw$ has a natural connection, $\nabla_{d/dt}
= d/dt$, and is decomposed into a direct sum of line bundles:
\begin{equation}
  \gtw \simeqq \bigoplus_{(a,b)\ne(0,0)} L_{a,b},
\label{decomp:gtw}
\end{equation}
where the indices $(a,b)$ runs through $(\Z/N\Z)^N\setminus\{(0,0)\}$.  Here
the line bundle $L_{a,b}$ on $X$ is defined by
\begin{equation}
    L_{a,b} := (\C \times \C)/{\approx_{a,b}},
\label{def:L-a,b}
\end{equation}
where $\approx_{a,b}$ is an equivalence relation defined by
\begin{equation}
    (t, x) \approx_{a,b} (t+1,    \eps^a x)
           \approx_{a,b} (t+\tau, \eps^b x).
\label{def:approx-a,b}
\end{equation}
We regard $L_{a,b}$ as a line subbundle of $\gtw$ through the injection given
by
\begin{equation}
    L_{a,b} \owns (t, x) \mapsto (t, x J_{a,b}) \in \gtw,
\label{L->gtw}
\end{equation}
where $J_{a,b}$ is the element of $\g$ defined by
\begin{equation}
    J_{a,b} 
    := \beta^a \alpha^{-b}.
\label{def:J-ab}
\end{equation}
We remark that $\{\,J_{a,b}\mid(a,b)\in(\Z/N\Z)^N\setminus\{(0,0)\}\,\}$
is a basis of $\g=sl_N(\C)$.

The space of meromorphic sections of $L_{a,b}$ over $X$ pulled back
to $\C$ can be canonically identified with
\begin{equation}
  K_{a,b}
  =
  \{\, f \in{\cal M}({\C})
  \mid f(t+1)=\eps^a f(t),\; f(t+\tau)=\eps^b f(t) \,\}.
\label{def:K-a,b}
\end{equation}
Here ${\cal M}(\C)$ is the space of meromorphic functions on $\C$.
($K_{0,0}$ is the space of elliptic functions and corresponds to the
trivial line bundle on $X$.)  The mapping 
\begin{equation}
    f \mapsto (t, f(t)) \text{ modulo $\approx_{a,b}$}
\label{isom:K->L}
\end{equation}
gives a canonical isomorphism from $K_{a,b}$ onto 
$H^0(X, L_{a,b} \tensor \K_X)$. 

The Liouville theorem implies that the only holomorphic function in $K_{a,b}$
is zero when $(a,b)\ne(0,0)$.  This is equivalent to $H^0(X,L_{a,b})=0$.
Since $L_{a,b}^*\simeqq L_{-a,-b}$ and the canonical line bundle of $X$ is
trivial, it follows from the Serre duality that $H^1(X,L_{a,b})=0$.  Thus we
obtain a simple vanishing result $H^p(X,L_{a,b})=0$ and therefore from the
decomposition \eqref{decomp:gtw} we obtain the following result.

\begin{lem}
\label{lem:vanishing}
$H^0(X, \gtw) = H^1 (X, \gtw) = 0$.
\end{lem}

\begin{example}
For $N=2$, matrices $\alpha$ and $\beta$ are nothing but the Pauli
matrices $\sigma^1$ and $\sigma^3$. The Jacobian elliptic functions $\sn$,
$\cn$, and $\dn$ are meromorphic functions in $K_{1,0}$, $K_{1,1}$, and
$K_{0,1}$ and can be regarded as meromorphic sections of 
the line bundles $L_{1,0}$, $L_{1,1}$, and $L_{0,1}$ respectively. 
\end{example}

\begin{example}
For general $N$ and each $(a,b)\in(\Z/N\Z)^2\setminus\{(0,0)\}$, we define
the function $w_{a,b}$ by
\begin{equation}
  w_{a,b}(t)
  = w_{a,b}(\tau;t)
  :=
  \frac{\theta'_{[0,0]}}{\theta_{[a,b]}}
  \frac{\theta_{[a,b]}(t;\tau)}{\theta_{[0,0]}(t;\tau)}.
\label{def:w-a,b}
\end{equation}
(See \eqref{def:theta[a,b]} in \appref{theta-func} for the notation.) The
function $w_{a,b}(t)$ on $\C$ is uniquely characterized by the following
properties:
\begin{enumerate}
\item The function $w_{a,b}(t)$ is a meromorphic functions in $K_{a,b}$ and
  hence can be regarded as a global meromorphic section of $L_{a,b}$;
\item The poles of $w_{a,b}(t)$ are all simple and contained in $\Z + \Z\tau$;
\item The residue of $w_{a,b}(t)$ at $t=0$ is equal to $1$.
\end{enumerate}
Because of these properties, it will play an important role in
concrete computations in later sections. For convenience of those
computations, let us list several other properties of $w_{a,b}(t)$:

\begin{itemize}
\item
The Laurent expansion at $t=0$ is equal to 
\begin{equation}
    w_{a,b}(t)
    = \frac{1}{t} + \sum_{\nu=0}^\infty w_{a,b,\nu} t^\nu
    = \frac{1}{t} + w_{a,b,0} + w_{a,b,1} t + \cdots,
    \label{w-a,b:expand}
\end{equation}
where the coefficients are written in the following forms:
\begin{equation}
    w_{a,b,0} = \frac{\theta'_{[a,b]}}{\theta_{[a,b]}}, \quad
    w_{a,b,1} = \frac{\theta''_{[a,b]} }{2 \theta_{[a,b]} }
              - \frac{\theta'''_{[0,0]}}{6 \theta'_{[0,0]}}, \quad
    \ldots
    \label{w-a,b:coeff}
\end{equation}

\item
Formulae \eqref{theta:odd} and \eqref{theta:char-period} imply
\begin{equation}
    w_{-a,-b}(t) = - w_{a,b}(-t), \quad
    w_{a,b}(t) = w_{a',b'}(t)
    \text{ if $a\equiv a'$ and $b\equiv b'$ mod $N$}.
\label{w-a,b:property}
\end{equation}

\item
\lemref{jacobi} will be used in the following form:
\begin{equation}
    \sum_{(a,b)\neq(0,0)} w_{a,b,1} = 0,
\label{w-a,b,1:sum}
\end{equation}
where the summation is taken over all 
$(a,b)\in (\Z/N\Z)^2\setminus\{(0,0)\}$. 
\end{itemize}
\end{example}

\subsection{Definition of the spaces of conformal coinvariants and conformal
  blocks} 
\label{def-conf-block}

In this section we define a conformal block associated to the twisted Lie
algebra bundle $\gtw$ defined in \secref{bundles}.

First let us introduce notation of sheaves.  As usual, the structure sheaf on
$X=X_\tau$ is denoted by $\O_X$ and the sheaf of meromorphic functions on $X$
by $\K_X$. A stalk of a sheaf $\F$ on $X$ at a point $P \in X$ is denoted by
$\F_P$.  When $\F$ is a $\O_X$-module, we denote its fiber $\F_P/\m_P\F_P$ by
$\F|_P$, where $\m_P$ is the maximal ideal of the local ring $\O_{X,P}$.
Denote by $\F^\wedge_P$ the $\m_P$-adic completion of $\F_P$.

We shall use the same symbol for a vector bundle and for a locally free
coherent $\O_X$-module consisting of its local holomorphic sections.  For
instance, the invertible sheaf associated to the line bundle $L_{a,b}$ is also
denoted by the same symbol $L_{a,b}$. Denote by $\Omega_X^1$ the sheaf of
holomorphic 1-forms on $X$, which is isomorphic to $\O_X$ since $X$ is an
elliptic curve.  The fiberwise Lie algebra structure of the bundle $\gtw$
induces that of the associated sheaf $\gtw$ over $\O_X$.  Define the invariant
$\O_X$-inner product on $\gtw$ by
\begin{equation}
  (A|B) := \frac{1}{2N} \tr_{\gtw}(\ad A \ad B) \in \O_X
  \quad
  \text{for $A,B\in\gtw$},
  \label{def:inner-prod-gtw}
\end{equation}
where the symbol $\ad$ denotes the adjoint representation of the $\O_X$-Lie
algebra $\gtw$.  Then the inner product on $\gtw$ is invariant under the
translations with respect to the connection 
$\nabla:\gtw\to\gtw\tensor_{\O_X}\Omega_X^1$:
\begin{equation}
  d(A|B) = (\nabla A|B) + (A|\nabla B) \in \Omega_X^1
  \quad
  \text{for $A,B\in\gtw$}.
  \label{compati-conn-innerprod}
\end{equation}
Under the trivialization of $\gtw$ defined by the construction
\eqref{def:gtw}, the connection $\nabla$ and 
the inner product $(\cdot|\cdot)$ on $\gtw$ respectively coincide with the
exterior derivative by $t$ and the inner product defined by
\eqref{def:inner-product}.

For any point $P$ on $X$ with $t(P)=z$, we put
\begin{equation*}
    \g^P := (\gtw \tensor_{\O_X} \K_X)^{\wedge}_P,
\end{equation*}
which is a topological Lie algebra non-canonically isomorphic to the loop
algebra $\g((t-z))$. Its subspace $\g^P_+:=(\gtw)_P^{\wedge}\simeqq\g[[t-z]]$
is a maximal linearly compact subalgebra of $\g^P$ under the $(t-z)$-adic
linear topology.

Let us fix mutually distinct points $Q_1, \ldots, Q_L$ on $X$ whose
coordinates are $t = z_1, \ldots, z_L$ and put $D:=\{Q_1,\dots,Q_L\}$.  We
shall also regard $D$ as a divisor on $X$ (i.e., $D=Q_1+\cdots+Q_L$).  Denote
$X\setminus D$ by $\dot X$.  The Lie algebra
$\g^D:=\bigoplus_{i=1}^L\g^{Q_i}$ has the natural 2-cocycle defined by
\begin{equation}
  \caff(A,B) := \sum_{i=1}^L \Res\limits_{t = z_i} (\nabla A_i | B_i),
\label{def:cocycle}
\end{equation}
where $A=(A_i)_{i=1}^L, B=(B_i)_{i=1}^L \in \g^D$ and $\Res_{t=z}$ is the
residue at $t=z$. (The symbol ``$\caff$'' stands for ``Cocycle defining the
Affine Lie algebra''.)  We denote the central extension of $\g^D$ with respect
to this cocycle by $\hat\g^D$:
\begin{equation*}
    \hat\g^D := \g^D \oplus \C \khat,
\end{equation*}
where $\khat$ is a central element.  Explicitly the bracket of $\hat\g^D$ is
represented as
\begin{equation}
  [A, B] = ([A_i, B_i]^0)_{i=1}^L \oplus \caff(A,B) \khat
  \quad
  \text{for $A,B\in\g^D$,}
  \label{def:aff-alg-str}  
\end{equation}
where $[A_i, B_i]^0$ are the natural bracket in $\g^{Q_i}$.  The Lie algebra
$\hat\g^P$ for a point $P$ is non-canonically isomorphic to the affine Lie
algebra $\hat\g$ of type $A^{(1)}_{N-1}$ (a central extension of the loop
algebra $\g((t-z))=sl_N\bigl(\C((t-z))\bigr)$).  If $P=Q_i$ for $i=1,\dots,L$,
then $\hat\g^P=\hat\g^{Q_i}$ can be regarded as a subalgebra of $\hat\g^D$.
Put $\g^P_+:=(\gtw)^\wedge_P$ as above.  Then $\g^{Q_i}_+$ can be also
regarded as a subalgebra of $\hat\g^{Q_i}$ and $\hat\g^D$.

Let $\gDpr$ be the space of global meromorphic sections of $\gtw$ which are
holomorphic on $\dot X$:
\begin{equation*}
    \gDpr := \Gamma(X, \gtw(\ast D)).
\end{equation*}
There is a natural linear map from $\gDpr$ into $\g^D$ which maps a
meromorphic section of $\gtw$ to its germ at $Q_i$'s. As in the
non-twisted case (e.g., \S2.2 of \cite{tu-ue-ya:89}), the residue theorem
implies that this linear map is extended to a Lie algebra injection from
$\gDpr$ into $\hat\g^D$, which allows us to regard $\gDpr$ as a subalgebra
of $\hat\g^D$.

\begin{defn}
\label{def:CC-CB}
The space of {\em conformal coinvariants} $\CC(M)$ and that of 
{\em conformal blocks} $\CB(M)$ associated to $\hat\g^{Q_i}$-modules $M_i$
with the same level $\khat = k$
are defined to be the space of coinvariants of $M := \bigotimes_{i=1}^L M_i$
with respect to $\gDpr$ and its dual:
\begin{equation}
    \CC(M) := M/\gDpr M, \quad
    \CB(M) := (M/\gDpr M)^\ast.
\label{def:conf-block}
\end{equation}
(In \cite{tu-ue-ya:89}, $\CC(M)$ and $\CB(M)$ are called the space of {\em
covacua} and that of {\em vacua} respectively.)  In other words, the space
of conformal coinvariants $\CC(M)$ is generated by $M$ with relations
\begin{equation}
    A_\prin v \equiv 0
\label{ward-cc}
\end{equation}
for all $A_\prin \in \gDpr$ and $v\in M$,
and a linear functional $\Phi$ on $M$ belongs to the space of conformal
blocks $\CB(M)$ if and only if it satisfies that
\begin{equation}
    \Phi(A_\prin v) = 0
    \quad\text{for all $A_\prin \in \gDpr$ and $v \in M$.}
\label{ward-cb}
\end{equation}
These equations \eqref{ward-cc} and \eqref{ward-cb} are called  the
{\em Ward identities}.
\end{defn}

The most important conformal blocks for our purpose are constructed from
Weyl modules (or generalized Verma modules) which are determined from the
following data:
\begin{itemize}

\item A parameter $k\in\C$ which is called the level of the model.

\item Finite-dimensional irreducible representations $V_i$ of the fiber Lie
  algebra $\gtw|_{Q_i}$ isomorphic to $\g$.

\end{itemize}

Put $\g^{Q_i}_+:=(\gtw)^\wedge_{Q_i}$ and
$\hat\g^{Q_i}_+:=\g^{Q_i}_+\oplus\C\khat$, which are subalgebras of
$\hat\g^{Q_i}$.
The subalgebra $\hat\g^{Q_i}_+:=\g^{Q_i}_+\oplus\C\khat$ of $\hat\g^{Q_i}$
acts on $V_i$ through the linear map $\khat \mapsto k\,\id_{V_i}$ and the
natural projection $\g^{Q_i}_+\to\gtw|_{Q_i}$, $A\mapsto A(Q_i)$.  The
$\hat\g^{Q_i}$-module induced from $V_i$ is called a {\em Weyl
  module} or a {\em generalized Verma module}:
\begin{equation}
  \Weyl_{k}(V_i)
  := \Ind_{\hat\g^{Q_i}_+}^{\hat\g^{Q_i}} V_i
  = U(\hat\g^{Q_i}) \tensor_{U(\hat\g^{Q_i}_+)} V_i
\label{def:weyl-mod}
\end{equation}
See \cite{kazh-lus:93} for properties of Weyl modules.

The space of conformal coinvariants and that of conformal blocks associated to
the data $(Q,V)=(\{Q_i\}, \{V_i\})$ are defined to be the space of
conformal coinvariants and that of conformal blocks associated to the
$\hat\g^D$-module
\begin{equation*}
  \Weyl_k(V) :=
  \bigotimes_{i=1}^L \Weyl_{k}(V_i),
\end{equation*}
on which the center $\khat$ acts as multiplication by $k$.  Namely we define
them as follows:
\begin{align}
    \CC_k(Q,V) = \CC_k(\{Q_i\}, \{V_i\})
    :=& \CC(\Weyl_k(V)) = \Weyl_k(V)/ \gDpr \Weyl_k(V),
\label{def:coinvariants:weyl}
\\
    \CB_k(Q,V) = \CB_k(\{Q_i\}, \{V_i\})
    :=& \CB(\Weyl_k(V)) = (\Weyl_k(V)/ \gDpr \Weyl_k(V))^\ast.
\label{def:conf-block:weyl}
\end{align}
Hereafter we use the word ``conformal coinvariants'' and ``conformal
block'' for this kind of conformal coinvariants and conformal
blocks, namely those associated to Weyl modules, unless
otherwise stated.

It is easy to see that the spaces of conformal coinvariants and conformal
blocks are determined by the finite-dimensional part
$V=\bigotimes_{i=1}^L V_i$ of $\Weyl_k(V)$ as is the case with the space
of conformal coinvariants and conformal blocks on $\bbbP^1(\C)$
(e.g., Lemma 1 of \cite{f-f-r:94}), because of the cohomology
vanishing. In fact, \lemref{lem:vanishing} implies a decomposition,
\begin{equation}
  \hat\g^D = \gDpr \oplus \hat\g^D_+,
\label{decomp:g-hat}
\end{equation}
where $\hat\g^D_+ = \bigoplus_{i=1}^L \g^{Q_i}_+ \oplus \C \khat$. Hence
we have, as left $\gDpr$-modules, 
\begin{equation}
  \Weyl_k(V)
  = \Ind_{\hat\g^D_+}^{\hat\g^D} V
  = U(\hat\g^D) \tensor_{U(\hat\g^D_+)} V
  = U(\gDpr) \tensor_\C V,
\label{decomp:bbbV}
\end{equation}
where $V:=\bigotimes_{i=1}^L V_i$ and the action of $\hat\g^D_+$ is defined by
the mapping $\khat \mapsto k\cdot\id$ and the natural projection 
$\g^D_+ \to \prod_{i=1}^L(\gtw|_{Q_i})$.  Therefore, due to the Ward identity
\eqref{ward-cb} and the definition of the Weyl module \eqref{def:weyl-mod},
the space of conformal coinvariants is canonically isomorphic to the tensor
product of $\gtw|_{Q_i}$-modules by the natural inclusion map
$V=\bigotimes_{i=1}^L V_i\injto\Weyl_k(V)$, 
$\bigotimes_{i=1}^L v_i \mapsto \bigotimes_{i=1}^L(1\tensor v_i)$.  (In the
following we shall identify $v_i\in V_i$ with $1\tensor v_i\in\Weyl_k(V_i)$.)

\begin{prop}
\label{restriction}
The inclusion map $V\injto\Weyl_k(V)$ and the induced restriction map
$\Weyl_k(V)^\ast\onto V^\ast$ induce the following isomorphisms
respectively:
$$
    \CC_k(Q,V) 
    \isoot
    V = \bigotimes_{i=1}^L V_i
    \quad\text{and}\quad
    \CB_k(Q,V) 
    \isoto
    V^\ast = \bigotimes_{i=1}^L V_i^\ast.
$$
\end{prop}

For any point $P\in X=X_\tau$, Let us denote the $1$-dimensional trivial
representation of $\gtw|_P$ by $\C_P = \C u_P$.  Then the 
proposition above readily leads to the following corollary.

\begin{cor}
Let $P$ be a point of $X$ distinct from $Q_i$'s. Then the canonical
inclusion $\Weyl_k(V) \injto \Weyl_k(\C_P) \tensor \Weyl_k(V)$, 
$v \mapsto u_P \tensor v$, induces an isomorphism
\begin{equation}
  \CB_k(\{P,Q_i\},\{\C_P, V_i\}) \isoto \CB_k(\{Q_i\},\{V_i\}).
\label{propagation}
\end{equation}
\end{cor}

The property above is called {\em propagation of vacua} in
\cite{tu-ue-ya:89}. In our case the proof is far simpler due to
\propref{restriction}, as in the case of $\bbbP^1$. 
(cf.\ \S3 of \cite{f-f-r:94}.) 

\subsection{Correlation functions}
\label{correl}

The current and the energy-momentum correlation functions are defined
as in \S2 of \cite{tu-ue-ya:89}, but we must take twisting into account
and use the decomposition \eqref{decomp:gtw}.

First we consider the current correlation functions.  Let $\Phi$ be a
conformal block in $\CB_k(Q,V)$ and $v$ a vector in $\Weyl_k(V)$.  There
exists a unique $\omega_i\in
(L_{a,b}^\ast\tensor_{\O_X}\K_X\tensor_{\O_X}\Omega_X^1)_{Q_i}^{\wedge}$ with
the property that
\begin{equation}
    \Res_{t=z_i}\langle f_i, \omega_i \rangle 
    := \Phi(\rho_i(f_i J_{a,b}) v)
    \quad
    \text{for all $f_i\in(\K_X)_{Q_i}^{\wedge}$},
\label{def:omega-i:res}
\end{equation}
where $\langle \cdot, \cdot \rangle$ is the canonical pairing of
$(L_{a,b})_P^{\wedge}$ and $(L_{a,b}^\ast)_P^{\wedge}$, $f_i J_{a,b}$ can be
regarded as an element of $\g^{Q_i}$ by means of \eqref{L->gtw} and its
action on the $i$-th component of $v$ (namely, the $\Weyl_k(V_i)$-component)
is denoted by $\rho_i(f_i J_{a,b})$.  Thus we obtain a linear functional
\begin{equation*}
  \sum_{i=1}^L\Res_{Q_i}\langle\cdot,\omega_i\rangle : 
  \bigoplus_{i=1}^L (L_{a,b}\tensor_{\O_X}\K_X)_{Q_i}^{\wedge} \to \C
\end{equation*}
which maps $(f_i)_{i=1}^L\in(L_{a,b}\tensor_{\O_X}\K_X)_{Q_i}^{\wedge}$ to 
$\sum_{i=1}^L\Res_{Q_i}\langle f_i,\omega_i\rangle\in\C$.

The Ward identity \eqref{ward-cb} implies that
\begin{equation*}
  \sum_{i=1}^L \Res_{Q_i} \langle f_{Q_i}, \omega_i \rangle
  = \sum_{i=1}^L \Phi(\rho_i(f_{Q_i} J_{a,b}) v)
  = 0
\end{equation*}
for any meromorphic section $f \in H^0(X,L_{a,b}(\ast D))$, where
$f_{Q_i}$ is the germ of $f$ at $Q_i$. 
Since $H^0(X,L_{a,b}(\ast D))$ and
$H^0(X,L_{a,b}^\ast\tensor_{\O_X}\Omega_X^1(\ast D))$ are orthogonal
complements to each other under the residue pairing of 
$\bigoplus_{i=1}^L (L_{a,b}\tensor_{\O_X}\K_X)_{Q_i}^{\wedge}$ and
$\bigoplus_{i=1}^L
(L_{a,b}^\ast\tensor_{\O_X}\K_X\tensor\Omega_X^1)_{Q_i}^{\wedge}$
(cf.\ \cite{tate:68} or Theorem 2.20 of \cite{iwa:52}), 
we have a meromorphic 1-form $\omega$ with values
in $L_{a,b}^\ast$ such that the germ of $\omega$ at $Q_i$ gives $\omega_i$ and
is holomorphic outside of $\{Q_1, \ldots, Q_L\}$:
\begin{equation}
  \omega \in H^0(X, L_{a,b}^\ast\tensor_{\O_X}\Omega_X^1(\ast D)),\quad
  (\omega)_{Q_i} = \omega_i.
\label{def:omega}
\end{equation}

In order to define the correlation functions, we need explicit expression
of the action of $(\gtw)^\wedge_P$ which we can identify with the affine
Lie algebra by fixing a trivialization of $\gtw$ around $P$.

Let $P$ be any point of $X$ and $z(P)$ a point of $\C$ whose image in
$X = \C/(\Z+\tau\Z)$ is equal to $P$. 
The description \eqref{def:gtw} of $\gtw$ naturally determines a
local trivialization of $\gtw$ at $P$, once we fix the coordinate $t=z(P)$ of
$P$.  By means of this trivialization, we fix isomorphisms $\g^P \simeqq
\hat\g$, $\gtw|_P \simeqq \g$, and so on.
The induced trivialization of $L_{a,b}$ at $P$ is the same as
the trivialization defined by the isomorphism \eqref{isom:K->L}:
\begin{equation}
  (L_{a,b})_P \isoto \O_{\C,z(P)} \isoot \O_{X,P}
  \quad
  \text{for $(a,b)\in(\Z/N\Z)^2$},
\label{isom:L->O<-O}
\end{equation}
which corresponds $\nabla$-flat sections of $L_{a,b}$ to constant functions on
$X$.   The decomposition \eqref{decomp:gtw}
of $\gtw$ induces a decomposition of its stalk at $P$ and is consistent with
the trivializations above:
\begin{equation}
  \gtw_P 
  = \bigoplus J_{a,b}(L_{a,b})_P
  \simeqq \bigoplus J_{a,b}\O_{X,P}
  = \g\tensor\O_{X,P},
\label{triv:gtw=g*O}
\end{equation}
where the indices $(a,b)$ run through $(\Z/N\Z)^2\setminus\{(0,0)\}$.

Let $\xi$ be a local coordinate at $P$.  For $A\in\g$ and $m\in\Z$, we denote
by $A[m]$ the element of $(\gtw)_P^{\wedge}\subset\hat\g^P$ which is
represented by $A\xi^m$ under the trivialization \eqref{triv:gtw=g*O}.  Since
$A[0]$ is $\nabla$-flat for $A\in\g$ (i.e., $\nabla A[0]=0$), the bracket of
$\hat\g^P$ is represented as:
\begin{equation*}
  [A[m], B[n]] = [A,B][m+n] + (A|B) m\delta_{m+n,0} \khat
  \quad
  \text{for $A,B\in\g$ and $m,n\in\Z$},
\end{equation*}
which coincides with the usual commutation relation of the affine Lie algebra.

\begin{lem}
\label{lem:current-correl}
Under the situation above, let $P$ be in $\dot X$ (i.e., distinct from
$Q_i$'s) and $\tilde\Phi$ the conformal block in
$\CB_k(\{P,Q_i\},\{\C_P,V_i\})$ corresponding to $\Phi$ by the
isomorphism \eqref{propagation}.  Then we have the following:
\par\noindent{\rm (i)}\enspace
Take $x \in L_{a,b}|_P$ and let $f_x$ be an element of
$(L_{a,b})_P^{\wedge}$ with a principal part $x/\xi$ 
(i.e., $f_x = (x/\xi + \text{{{\rm regular}}})$).  Then 
$\tilde\Phi ( f_x J_{a,b}u_P \tensor v)\,d\xi$ 
does not depend on the choice of \/$\xi$ and $f_x$. Thus we can define
$\tilde\Phi(J_{a,b}[-1]u_P\tensor v)\,d\xi
\in(L_{a,b}^\ast\tensor_{\O_X}\Omega_X^1)|_P$ by
\begin{equation*}
  \langle x, \tilde\Phi(J_{a,b}[-1]u_P \tensor v) \rangle \,d\xi
  :=
  \tilde\Phi ( f_x J_{a,b}u_P \tensor v )\,d\xi.
\end{equation*}
\noindent{\rm (ii)}\enspace
The following equation holds at any point $P$ in 
$\dot X = X\setminus\{Q_1,\ldots, Q_L\}$:
\begin{equation}
  \omega(P) = \tilde\Phi(J_{a,b}[-1]u_P \tensor v)\,d\xi
  \in (L_{a,b}^\ast\tensor_{\O_X}\Omega_X^1)|_P.
\label{current-correl}
\end{equation}
\end{lem}
\begin{pf}
The statement (i) can be shown by the same argument as the proof of Claim
1 of Theorem 2.4.1 of \cite{tu-ue-ya:89}.
Using the Riemann-Roch theorem (or the function $w_{a,b}(t)$ defined by
\eqref{def:w-a,b}), we can choose $f_x$ in (i) from 
$H^0(X,L_{a,b}(\ast D+P))$.  Then we have
\begin{align*}
    & \langle x, \tilde\Phi(J_{a,b}[-1]u_P \tensor v) \rangle \,d\xi
    = \tilde\Phi ( f_x J_{a,b}u_P \tensor v )\,d\xi \\
    &= - \sum_{i=1}^L \tilde\Phi(u_P\tensor\rho_i(f_x J_{a,b})v)\,d\xi \\
    &= - \sum_{i=1}^L \Phi(\rho_i(f_x J_{a,b}) v)\,d\xi 
     = - \sum_{i=1}^L \Res_{Q_i} \langle f_x, \omega\rangle\,d\xi \\
    &= \Res_P \langle f_x, \omega \rangle \,d\xi 
     = \Res_P \frac{\langle x, \omega \rangle}{\xi}\,d\xi
     = \langle x, \omega|_P \rangle.
\end{align*}
Here we have used the Ward identity \eqref{ward-cb} and the residue theorem.
Thus we have proved the equation \eqref{current-correl}.
\end{pf}

\begin{defn}
We call this 1-form $\omega = \tilde\Phi(J_{a,b}[-1]u_P \tensor v)\,d\xi$
a {\em correlation function of the current} $J_{a,b}(\xi)$ and $v$ under
$\Phi$, or a {\em current correlation function} for short, and denote it by
$\Phi(J_{a,b}(P)v)\, dP$ or $\Phi (J_{a,b}(\xi)v)\, d\xi$ 
when we fix a local coordinate $\xi$.
\end{defn}

We now proceed to the definition of the energy-momentum correlation
functions.

\begin{lem}
\label{lem:sug-correl}
Let $P$ be in $\dot X$ and $\xi$ a local coordinate defined on an open
neighborhood $U$ of $P$. Then the following expression gives a holomorphic
section of $\Omega_X^2(\ast D)=(\Omega_X^1)^{\tensor 2}(\ast D)$ on
sufficiently small $U$:
\begin{multline}
  \Phi(S(P) v) \,(d\xi(P))^2 
  \\
  :=
  \frac{1}{2}
  \lim_{P'\to P} 
  \Biggl(
    \sum_{(a,b)} 
    \Phi( J_{a,b}(\xi(P)) \, J^{a, b}(\xi(P')) v )\,
  \\
    - \frac{k \dim \g}{(\xi(P) - \xi(P'))^2}
  \Biggr)
  \,d\xi(P)\,d\xi(P').
\label{def:sug-correl}
\end{multline}
Here the indices $(a,b)$ run through 
$(\Z/N \Z)^2 \smallsetminus \{(0,0)\}$, and $J^{a,b}$ is 
the dual basis of $J_{a,b}$ with respect to $(\cdot|\cdot)$, namely,
$J^{a,b} = J_{-a,-b}/N$.
\end{lem}
\begin{pf}
The argument in the proof of \lemref{lem:current-correl} and the Hartogs
theorem of holomorphy show that the current correlation function
$\Phi(J_{a,b}(\xi)J^{a,b}(\zeta) v)\,d\xi\,d\zeta$ defines a global section on
$X\times X$ of the sheaf
\(
  \F :=
  (L_{ a, b}^\ast\tensor_{\O_X}\Omega_X^1(\ast D)) \boxtimes
  (L_{-a,-b}^\ast\tensor_{\O_X}\Omega_X^1(\ast D)) (\ast \Delta)
\), %
where $\Delta$ is the diagonal divisor of $X \times X$. 
We define 
\begin{equation*}
\Phi(S(P,P') v)\,d\xi(P)\,d\xi(P') \in H^0(U \times U,\F)
\end{equation*}
by
\begin{align*}
  & \Phi(S(P,P')v)\,d\xi(P)\,d\xi(P')
  \\
  & :=
  \left(
    \sum_{a,b}
    \Phi(J_{a,b}(\xi(P))\, J^{a,b}(\xi(P')) v)
    - \frac{k \dim \g}{(\xi(P) - \xi(P'))^2}
  \right)
  \,d\xi(P)\,d\xi(P').
\end{align*}

In order to show \lemref{lem:sug-correl}, first take a local coordinate $\xi$
on a sufficiently small neighborhood $U$ of $P$ with $\xi(P)=0$ and a local
trivialization of $L_{a,b}$ on $U$ and choose a meromorphic section 
$f\in H^0(X, L_{a,b}\tensor_{\O_X}\K_X)$ whose Laurent expansion has the form
\begin{equation}
    f(\xi) = \xi^{-1} + \text{(regular at $\xi=0$)}.
\label{local-coord}
\end{equation}
The inclusion \eqref{L->gtw} and the Ward identity \eqref{ward-cb} imply that
$\Phi(J_{a,b}(\xi(P))J^{a,b}(\xi(P'))v)$ is equal to
\begin{equation}
\begin{aligned}
  & \Tilde{\Tilde{\Phi}}
  (J_{a,b}[-1]u_P \tensor J^{a,b}[-1]u_{P'} \tensor v)
  \\
  & =
  - \tilde\Phi((f J_{a,b})_{P'}J^{a,b}[-1]u_{P'} \tensor v)
  - \sum_{i=1}^L
    \tilde\Phi(J^{a,b}[-1]u_{P'} \tensor \rho_i((f J_{a,b})_{Q_i}) v),
  \hphantom{--} 
\end{aligned}
\label{tmp:em1}
\end{equation}
where
$\Tilde{\Tilde{\Phi}} \in \CB_k(\{P, P', Q_i\}, \{\C_P, \C_{P'}, V_i\})$,
and $\tilde \Phi \in \CB_k(\{P', Q_i\}, \{\C_{P'}, V_i\})$
correspond to $\Phi$ through the isomorphism \eqref{propagation}.
The second term of \eqref{tmp:em1} is regular as a function of $P$ at $P'$
as shown in \lemref{lem:current-correl} and the first term is rewritten as
\begin{equation}
  - \tilde\Phi((f J_{a,b})_{P'} J^{a,b}[-1]u_{P'} \tensor v) 
  =
  \frac{k}{(\xi(P) - \xi(P'))^2} \Phi(v) +
  \text{(regular at $P' = P$).}
\label{tmp:em2}
\end{equation}
(Details of the computation is the same as that of the proof of the
assertion (4) of Theorem 2.4.1 of \cite{tu-ue-ya:89}. 
Note that $[J_{a,b}, J^{a,b}] = 0$.)
Equations \eqref{tmp:em1} and \eqref{tmp:em2} mean that
$\Phi(S(P,P') v)\,d\xi(P)\,d\xi(P')$ is a holomorphic section of 
$(L_{ a, b}^\ast\tensor_{\O_X}\Omega_X^1(\ast D)) \boxtimes
 (L_{-a,-b}^\ast\tensor_{\O_X}\Omega_X^1(\ast D))$ on $U \times U$ for
the coordinate $\xi$ that satisfies \eqref{local-coord}.

Restricting it to the diagonal of $U\times U$, we obtain a local holomorphic
section
\begin{equation*}
    \Phi(S(P) v)\,(d\xi(P))^2 \in H^0(U, \Omega_X^2(\ast D)).
\end{equation*}
Note that since $L_{a,b}^\ast = L_{-a,-b}$, the factors $L_{a,b}^\ast$ and
$L_{-a,-b}^\ast$ cancel out on the diagonal.  Thanks to this fact, the
trivialization of $L_{a,b}$'s which we implicitly fixed in the argument above
does not affect the result.
\end{pf}

\begin{defn}
Put $\kappa:=k+\hvee$, where $\hvee$ is the dual Coxeter number of $\g$ 
(i.e., $\hvee = N$ because $\g=sl_N$). We call $\Phi(S(\xi) v) \, (d\xi)^2$ in
\eqref{def:sug-correl} a {\em correlation function of the Sugawara tensor}
$S(\xi)$ and $v$ under $\Phi$, or a
{\em Sugawara correlation function} for short, and
\begin{equation}
    \Phi(T(\xi) v) \, (d\xi)^2 := \kappa^{-1}\Phi(S(\xi) v) \,(d\xi)^2
\label{def:em-correl}
\end{equation}
a {\em correlation function of the energy-momentum tensor} $T(\xi)$ and $v$
under $\Phi$, or a {\em energy-momentum correlation function} for short 
when $\kappa\ne0$.
\end{defn}

Let us calculate the coordinate transformation law of the correlation function
of the Sugawara tensor.  Let $\zeta$ be another coordinate on $U$.  The
differential $d\xi(P) d\xi(P')/ (\xi(P) - \xi(P'))^2$ transforms under the
coordinate change $\xi \mapsto \zeta = \zeta(\xi)$ as
\begin{equation}
  \frac{d\zeta(P)\, d\zeta(P')}{(\zeta(P) - \zeta(P'))^2}
  =
  \frac{d\xi(P)\, d\xi(P')}{(\xi(P) - \xi(P'))^2}   
  + \frac{\{\zeta, \xi\}(P')}{6} \,d\xi(P)\,d\xi(P') 
  + O(\xi(P) - \xi(P')),
  \label{transform:2-diff}
\end{equation}
where $\{\zeta, \xi\} = \zeta'''/\zeta' - 3/2 (\zeta''/\zeta')^2$
($\zeta'=d\zeta/d\xi$) is the Schwarzian derivative of $\zeta = \zeta(P)$ with
respect to $\xi = \xi(P)$.  Hence the correlation function of the Sugawara
tensor transforms with respect to a coordinate change 
$\xi \mapsto \zeta=\zeta(\xi)$ not as 2-differentials but as
\begin{equation}
    \Phi(S(\zeta) v) d\zeta^2
    =
    \Phi(S(\xi) v) d\xi^2 + \frac{k \dim\g}{12} \{\zeta, \xi\} \Phi(v).
\label{transform:sug-correl}
\end{equation}
This means that the family $\{\Phi(S(\xi) v) d\xi^2\}$ defines a meromorphic
projective connection on $X$.  (For the notion of projective connections on
Riemann surfaces, see \cite{Gun}.) The Schwarzian derivative $\{\zeta,\xi\}$
vanishes identically if and only if $\zeta$ is a fractional linear
transformation of $\xi$ (i.e., $\zeta=(a\xi+b)/(c\xi+d)$).  From this fact
it follows that $\{\Phi(S(\xi)v)\,d\xi^2\}$ behaves like a 2-differential
under fractional linear coordinate changes. 

For later use, we compute the local expression of the energy-momentum
correlation function around $Q_i$. We take a holomorphic local chart
$(U,\zeta)$ with $Q_i\in U$ and $\zeta(Q_i)=0$, a local trivialization of
$L_{a,b}$ by \eqref{isom:L->O<-O}, and a trivialization of $\gtw_{Q_i}$ by
\eqref{triv:gtw=g*O}. Under these trivializations, we have, due to
\eqref{def:omega},
\begin{multline}
  \Phi(J_{a,b}(P_1) J^{a,b}(P_2) v) \,dP_1\,dP_2
  \\
  =
  \Biggl(
    \sum_{i=1}^L \sum_{n,m \in \Z} \zeta_1^{-m-1} \zeta_2^{-n+m-1}
    \Phi(\rho_i(\NP J^{a,b}[n-m]J_{a,b}[m] \NP) v) 
  \\
  + \frac{k \Phi(v)}{(\zeta_2 - \zeta_1)^2}
  \Biggr) \,d\zeta_1\,d\zeta_2
\label{tmp:em3}
\end{multline}
if $|\zeta_1| > |\zeta_2|$. Here $P_1,P_2\in U$, $\zeta_1=\zeta(P_1)$,
$\zeta_2=\zeta(P_2)$, and $\NP\quad\NP$ denotes the normal ordered product
defined by
\begin{equation}
  \NP A[m] B[n] \NP = 
  \begin{cases}
    A[m] B[n],                         & \text{if $m < n$,} \\
    {\tfrac12}(A[m] B[n] + B[n] A[m]), & \text{if $m = n$,} \\
    B[n] A[m],                         & \text{if $m > n$.}
  \end{cases}
\label{def:normal-order}
\end{equation}
Using \eqref{tmp:em3}, we obtain an expression of the correlation
function \eqref{def:sug-correl} around the point $Q_i$:
\begin{equation}
    \Phi(S(\zeta) v) (d\xi)^2
    = 
    \sum_{m \in \Z} \zeta^{-m-2} \Phi(\rho_i(S[m]) v) \,(d\zeta)^2,
\label{sug-correl:local}
\end{equation}
where $S[m]$ are the {\em Sugawara operators} defined by:
\begin{equation}
    S[m] = 
    \frac12
    \sum_{a,b} \sum_{n\in\Z} \NP J_{a,b}[m-n]\, J^{a,b}[n]\NP,
\label{def:sugawara}
\end{equation}
which satisfy the following commutation relations:
\begin{align}
  & [S[m], A[n]] = - \kappa n A[m+n]
  \quad\text{for $A\in\g$},
  \label{[Sm,An]}
  \\
  & [S[m], S[n]]
  = \kappa \left(
    (m-n) S[m+n] + \frac{k\dim\g}{12}(m^3 - m)\delta_{m+n,0}\,\id
    \right) 
\end{align}
In particular the Sugawara operators $S[m]$ commute with $\g^P$ if $\kappa=0$
(i.e., the level $k$ is critical).

When $\kappa=k+\hvee\ne0$, the usual Virasoro operators are defined by
normalizing $S[m]$:
\begin{equation}
    T[m] := \kappa^{-1} S[m],
\label{def:virasoro}
\end{equation}
which satisfy the well-known commutation relations:
\begin{align}
  & [T[m], A[n]] = - n A[m+n]
  \quad\text{for $A\in\g$},
  \\
  & [T[m], T[n]] 
  = (m-n) T[m+n] + \frac{c_k}{12}(m^3 - m)\delta_{m+n,0}\,\id,
  \label{comm-rel:virasoro}
\end{align}
where $c_k = k\dim\g/\kappa$. 

Later we fix the local coordinate at $Q_i$ to $\xi_i=t-z_i$ and the one at
$P\in X$ to $\xi=t-z$ with $t(P)=z$, where $t$ is the global coordinate of
$\C$ (cf.\ \eqref{def:gtw}) and $z_i=t(Q_i)$.

\begin{lem}
\label{global-correl}
In this coordinate $\Phi(S(\xi)v)\,(d\xi)^2$ and
$\Phi(T(\xi)v)(d\xi)^2$ can be extended to a global 2-differentials
$\Phi(S(t)v)(dt)^2$ and $\Phi(T(t)v)\,(dt)^2$.
\end{lem}
\begin{pf}
Under fractional linear coordinate changes, $\{\Phi(S(\xi)v)\,(d\xi)^2\}$
behaves like a 2-differential due to \eqref{transform:sug-correl}.  Since the
coordinate changes between two of $\xi$'s are merely translations, if $X$ is
covered by these coordinates, then $\{\Phi(S(\xi)v)\,(d\xi)^2\}$ gives
a meromorphic 2-differentials on $X$.
\end{pf}

\begin{rem}
Using Weierstra\ss' $\wp$-function, we can prove the lemma above in a more
explicit manner.  In fact, since
\begin{equation*}
    \wp(t_1 - t_2)\, dt_1 dt_2
    = \left( \frac{1}{(t_1 - t_2)^2} + O((t_1-t_2)^2) \right)\,dt_1\,dt_2
\end{equation*}
is a global meromorphic 2-form on $X \times X$ with a pole along the
diagonal, we can equivalently replace the definition
\eqref{def:sug-correl} by
\begin{multline}
  \Phi( S(P) v)\, (dt(P))^2
  \\
  :=
  \frac{1}{2}
  \lim_{P'\to P} 
  \Biggl(
    \sum_{a,b} \Phi( J_{a,b}(t(P)) J^{a,b}(t(P')) v)
  \\
  - k\dim\g\cdot\wp(t(P) - t(P'))
  \Biggr)
  \,dt(P)\, dt(P').
\label{def:global-correl}
\end{multline}
This definition is meaningful globally on $X$ and coincides with that of
the proof above.
\end{rem}

\subsection{Action of the Virasoro algebra}
\label{virasoro-action}

When the level $k$ is not $-\hvee$, the Lie algebras of formal meromorphic
vector fields at $Q_i$ are projectively represented on $\Weyl_k(V)$ through
the energy-momentum tensor.

We denote by $\T_X$ the tangent sheaf of $X=X_\tau$ (i.e., the sheaf of vector
fields on $X$). Let us fix a local coordinate at $Q_i$ to $\xi_i=t-z_i$, and
denote by $\T^D$ the direct sum of Lie algebra of formal meromorphic vector
fields at $Q_i$ for $i=1,\ldots,L$:
\begin{equation}
    \T^D := \bigoplus_{i=1}^L \T^{Q_i}, \quad
    \T^{Q_i} := (\T_X \tensor_{\O_X} \K_X)_{Q_i}^\wedge
    = \bigoplus_{i=1}^L \C((\xi_i)) \pd{}{\xi_i}.
\label{def:T-D}
\end{equation}
The {\em Virasoro algebra} $\Vir^D$ at $D$ is defined to be the central
extension of $\T^D$ by $\C \chat$:
\begin{equation}
  \Vir^D := \T^D \oplus \C \chat,
  \label{def:Vir-D}
\end{equation}
whose Lie algebra structure is defined by
\begin{multline}
  \left[
    \bigl( \theta_i(\xi_i)\partial_{\xi_i} \bigr)_{i=1}^L,
    \bigl(   \eta_i(\xi_i)\partial_{\xi_i} \bigr)_{i=1}^L
  \right]
  \\
  =
  \Bigl(
     (\theta_i(\xi_i)\eta'_i(\xi_i) - \eta_i(\xi_i)\theta'_i(\xi_i) )
     \partial_{\xi_i}
  \Bigr)_{i=1}^L
  \oplus
  \frac{\chat}{12}
  \sum_{i=1}^L
  \Res\limits_{\xi_i=0}(\theta_i'''(\xi_i)\eta_i(\xi_i)\,d\xi_i),  
  \label{def:Vir-bracket}
\end{multline}
where $\theta_i(\xi_i), \eta_i(\xi_i) \in \C((\xi_i))$ and 
$\partial_{\xi_i}:=\partial/\partial \xi_i$. 
When $L=1$, this Virasoro algebra $\Vir^D$ is the usual one defined as a
central extension of the Lie algebra of vector fields on a circle.

The action of $\theta_i = \theta_i(\xi_i)\partial_{\xi_i}\in \T^{Q_i}$
and that of $\theta = (\theta_i(\xi_i)\partial_{\xi_i})_{i=1}^L\in\T^D$ 
on $v \in \Weyl_k(V)$ and $\Phi \in (\Weyl_k(V))^\ast$ are given by 
\begin{alignat}{2}
    & T_i\{\theta_i\}v
    := - \sum_{m\in\Z} \rho_i(\theta_{i,m}T[m]) v,
    & \quad
    T\{\theta\}(v)
    & := \sum_{i=1}^L T_i\{\theta_i\}v,
\label{def:T(vec-field)}
\\
    & (T^\ast_i\{\theta_i\}\Phi)(v)
    := - \Phi(T_i\{\theta_i\}v),
    & \quad
    T^\ast\{\theta\}\Phi
    & := \sum_{i=1}^L T^\ast_i\{\theta_i\}\Phi,
\label{def:T*(vec-field)}
\end{alignat}
where
$\theta_i(\xi_i)\partial_{\xi_i} =
\sum_{m\in\Z}\theta_{i,m}\xi_i^{m+1}\partial_{\xi_i}$ is the
Laurent expansion (cf.\ \eqref{def:sugawara} and \eqref{def:virasoro}).
Note that we have
\begin{equation}
  (T^\ast_i\{\theta_i\}\Phi)(v)
  = \Res\limits_{\xi_i=0} 
    \langle \Phi(T(\xi_i)v)\,(d\xi_i)^2,
            \theta_i(\xi_i)\partial_{\xi_i} \rangle,
  \label{res-repre-T-eta}
\end{equation}
where $\langle \cdot, \cdot \rangle$ is the contraction of a
2-differential and a tangent vector. 

\begin{prop}
\label{prop:vec-field-action}
Let $\theta =(\theta_i\partial_{\xi_i})_{i=1}^L$ and 
$\eta =(\eta_i\partial_{\xi_i})_{i=1}^L$ be elements of $\T^D$.
Then the operators $T\{\theta\}$, $T\{\eta\}$ acting on $\Weyl_k(V)$ and
the operators $T^\ast\{\theta\}$, $T^\ast\{\eta\}$ acting on
$(\Weyl_k(V))^\ast$ satisfy
\begin{align} 
  & [T\{\theta\}, T\{\eta\}]
  =
  T\{[\theta, \eta]\}
  +
  \frac{c_k}{12}
  \sum_{i=1}^L
  \Res\limits_{\xi_i=0}(\theta_i'''(\xi_i)\eta_i(\xi_i)\,d\xi_i)\,\id,
  \label{comm-rel:T(vec-field)}
  \\
  & [T^\ast\{\theta\}, T^\ast\{\eta\}]
  =
  T^\ast\{[\theta, \eta]\}
  -
  \frac{c_k}{12}
  \sum_{i=1}^L
  \Res\limits_{\xi_i=0}(\theta_i'''(\xi_i)\eta_i(\xi_i)\,d\xi_i)\,\id,
  \label{comm-rel:T*(vec-field)}
\end{align}
where $c_k=k\dim\g/\kappa$ and $\kappa=k+\hvee$.  Namely, the definitions
\eqref{def:T(vec-field)} and \eqref{def:T*(vec-field)} define representations
of $\Vir^D$ on $\Weyl_k(V)$ with central charge $c_k$ and on
$(\Weyl_k(V))^\ast$ with central charge $-c_k$ respectively.
\end{prop}
This is a direct consequence of the definitions \eqref{def:T(vec-field)},
\eqref{def:T*(vec-field)} and the Virasoro commutation relations
\eqref{comm-rel:virasoro}.

As shown in \lemref{global-correl}, we have a global meromorphic
2-differential $\Phi(T(t)v)(dt)^2$ for $\Phi\in\CB_k(Q,V)$ and
$v\in\Weyl_k(V)$.  Therefore \eqref{res-repre-T-eta} and the residue
theorem imply the following lemma.

\begin{lem}
\label{global-vec-field-action}
Let $\Phi$ be a conformal block in $\CB_k(Q,V)$, $\theta(t)\partial_t$ in
$H^0(X,\T_X(\ast D))$, and $\theta_i(\xi_i) \partial_{\xi_i}$ the Laurent
expansion of $\theta(t)\partial_t$ in $\xi_i$ for each $i=1,\dots,L$.  Denote
$(\theta_i(\xi_i)\partial_{\xi_i})_{i=1}^L$ by $\theta$.  Then
$T^\ast\{\theta\} \Phi = 0$.
\end{lem}

\section{Critical level and the XYZ Gaudin model}
\label{xyz-gaudin}

In this section we restrict ourselves to the case $k = -\hvee = -N$, namely
the case when the level is critical.

We showed in \secref{def-conf-block} that the conformal block is determined by
its finite-dimensional part. (See \propref{restriction}.)  We shall see in
this section that the Sugawara tensor is expressed by integrals of motion of
the XYZ Gaudin model on this finite-dimensional space. Indeed, it shall be
shown that determining certain spaces of conformal blocks at the critical
level is equivalent to solving the XYZ Gaudin model.

First let us recall the definition of the XYZ Gaudin model, generalizing
the definition in \cite{gau:73}, \cite{skl-tak:96} to $sl_N$ case. Keeping
in mind that we will show the relation of the conformal field theory and
the XYZ Gaudin model, we use the same notation for $sl_N$-modules, points
$Q_i$ on a elliptic curve etc.\ as in the previous section
\secref{conf-block}, and fix local coordinates at each point $Q_i$ to
$\xi_i=t-z_i$, where $t$ is the global coordinate of $\C$. 

The Hilbert space of the model is a tensor product of the finite-dimensional
irreducible representation spaces of $sl_N(\C)$:
$V := \bigotimes_{i=1}^L V_i$. The
generating function $\hat\tau(u)$ of the integrals of motion of the model
is defined as the trace of square of the quasi-classical limit $\T(u)$ 
of the monodromy matrix of the spin chain model associated with the
Baxter-Belavin's elliptic $R$-matrix:
\begin{align}
  & \T(u) :=
  \sum_{i=1}^L  \sum_{(a,b)\neq(0,0)}
  w_{a,b}(u-z_i) J_{a,b} \tensor \rho_i(J^{a,b}),
  \\ &
  \begin{aligned}
    \hat\tau(u) :=& \frac12 \tr(\T(u))^2
    \\    
    =&
    \frac12 \sum_{i,j=1}^L \sum_{(a,b)\neq(0,0)}
    w_{-a,-b}(u - z_i) w_{a,b}(u - z_j) \rho_i(J_{a,b}) \rho_j(J^{a,b}),
  \end{aligned}
  \label{def:tau-hat}
\end{align}
where the indices of the summations over $(a,b)$ run through $a = 0,
\ldots, N-1$, $b = 0, \ldots, N-1$, $(a,b)\neq(0,0)$, and $w_{a,b}$ are
functions defined by \eqref{def:w-a,b}.  As before, $\rho_i$ is the
representation of $\g$ on the $i$-th factor $V_i$ of $V$.

The integrals of motion are encoded here in the following way:
\begin{equation}
    \hat\tau(u)
    = \sum_{i=1}^L C_i \wp(u - z_i)
    + \sum_{i=1}^L H_i \zeta(u - z_i) + H_0,
\label{tau-hat}
\end{equation}
where $\zeta$ and $\wp$ are Weierstra{\ss}' $\zeta$ and $\wp$ functions,
$C_i$ is the Casimir operator of $\g$ acting on $V_i$, i.e.,
\begin{equation}
    C_i = \frac12 \sum_{(a,b)\neq(0,0)}\rho_i(J_{a,b}) \rho_i(J^{a,b}),
\label{def:casimir}
\end{equation}
and $H_i$ ($i = 1, \ldots, L$) and $H_0$ are integrals of motion.
Operators $H_i$ satisfy $\sum_{i=1}^L H_i = 0$, and hence there are $L$
independent integrals of motion.

\begin{example}
\label{ex:gaudin:N=2}
When $N=2$, the Casimir operator $C_i$ is equal to $l_i(l_i+1)\id_{V_i}$,
where $l_i = (\dim V_i - 1)/2$, and $H_i$ are expressed as 
(cf.\ \cite{skl-tak:96}):
\begin{equation}
\begin{aligned}
    H_i =&
    \sum_{j \neq i} \sum_{(a,b)= (0,1),(1,1),(1,0)}
    w_{a,b}(z_i - z_j) \rho_i(J_{a,b}) \rho_j(J^{a,b}),
\\
    H_0 =& \frac12 \sum_{i=1}^L \sum_{(a,b)=(0,1),(1,1),(1,0)} \Bigl(
            - e_{a,b} \rho_i(J_{a,b}) \rho_i(J^{a,b}) \\
         &  +
              \sum_{j \neq i} w_{a,b}(z_i - z_j)
              \left(
              \zeta\Bigl(z_i - z_j +\frac{\omega_{a,b}}{2}\Bigr)
            - \zeta\Bigl(\frac{\omega_{a,b}}{2}\Bigr)
              \right) \rho_i(J_{a,b}) \rho_j(J^{a,b})
            \Bigr),
\end{aligned}
\label{def:H-i}
\end{equation}
where $\omega_{a,b} = a\tau + b$ and $e_{a,b} =
\wp(\omega_{a,b}/2)$. 
\end{example}

We interpret this system as a twisted WZW model at the critical level. Let
us come back to the situation in \secref{conf-block} and put $u=t(P)$.
The $sl_N$-module $V_i$ is assigned to the point $Q_i$ and regarded as
$\gtw|_{Q_i}$ module by the trivialization \eqref{triv:gtw=g*O}.
Assign the vacuum module $\Weyl_k(\C_P)$ at $P$ ($k=-\hvee=-N$). As before
$\tilde\Phi \in \CB_k(\{P, Q_i\}, \{\C_P, V_i\})$
corresponds to a conformal block $\Phi \in \CB_k(\{Q_i\}, \{V_i\})$ through the
isomorphism \eqref{propagation}. The correlation function 
of the Sugawara tensor 
$\tilde\Phi(S(t) (u_P \tensor v))\,(dt)^2$ has an expansion
\eqref{sug-correl:local} at $Q_i$ and at $P$,
\begin{equation}
\begin{split}
    \tilde\Phi(S(t)(u_P \tensor v))\,(dt)^2
    &= 
    \sum_{m\in \Z}
    \tilde\Phi(u_P \tensor \rho_i(S[m]) v) (t-z_i)^{-m-2}\,(dt)^2 \\
    &=
    \sum_{m\in \Z}
    \tilde\Phi(S[m] u_P \tensor v) (t-u)^{-m-2}\,(dt)^2,
\end{split}
\label{sug-correl-at-Q-P}
\end{equation}
respectively. The right-hand side of \eqref{sug-correl-at-Q-P} is
$\sum_{m \leq -2} \tilde\Phi(S[m] u_P \tensor v) (t-u)^{-m-2} (dt)^2$,
since $S[m] u_P = 0$ for all $m \geqq -1$. Hence, evaluating
\eqref{sug-correl-at-Q-P} at $t=u$, we have
\begin{equation}
    \sum_{m\in \Z} \Phi(\rho_i(S[m]) v) (u-z_i)^{-m-2}
    =
    \tilde\Phi(S[-2] u_P \tensor v).
\label{sug-correl-at-P}
\end{equation}

\begin{lem}
\label{gaudin-as-sug-correl}
Let $v$ be a vector in $V=\bigotimes_{i=1}^L V_i$. Then we have
\begin{equation*}
  S[-2]u_P \tensor v \equiv u_P \tensor \hat\tau(u)v
  \quad\text{in $\CC_k(\{P,Q_i\},\{\C_P,V_i\})$.}
\end{equation*}
Hence the right-hand side of \eqref{sug-correl-at-P} is equal to
$\Phi(\hat\tau(u)v)$.
\end{lem}
\begin{pf}
First note that
\begin{equation*}
    S[-2] u_P =
    \frac12 \sum_{(a,b)\neq(0,0)} J_{a,b}[-1] J^{a,b}[-1] u_P.
\end{equation*}
The key step is to exchange the operators $J_{a,b}[-1]$ and $J^{a,b}[-1]$
with operators acting on $v$ by
using the Ward identity \eqref{ward-cc}. Recall that the functions
$w_{a,b}(t-u)$ \eqref{def:w-a,b} in $K_{a,b}$ define meromorphic sections
\begin{equation}
    J_{a,b,P}(t) := w_{ a, b}(t-u) J_{a,b}, \quad
    J^{a,b,P}(t) := w_{-a,-b}(t-u) J^{a,b},
\label{def:J-a,b}
\end{equation}
of $\gtw$ through the inclusion \eqref{isom:K->L} and \eqref{L->gtw}.  These
sections belong to $\pr_\ast\gtw(P)$.  Since $J_{a,b,P}(t)$ has a Laurent
expansion
\begin{equation}
    J_{a,b,P}(t) 
    = \frac{J_{a,b}}{t-u} + w_{a,b,0} J_{a,b} 
                          + w_{a,b,1} J_{a,b} \cdot (t-u) + O((t-u)^2)
\label{J-a,b:expand}
\end{equation}
at $P$ (see \eqref{w-a,b:expand}), and $J^{a,b}(t)$ has a similar
expansion, we have
\begin{equation}
    J_{a,b}[-1] J^{a,b}[-1] u_P =
    (J_{a,b,P}(t))_P (J^{a,b,P}(t))_P u_P 
    - k w_{a,b,1} u_P,
\label{JJu-0}
\end{equation}
where $k=-\hvee=-N$.  Summing up \eqref{JJu-0} for $(a,b)$ and using
\eqref{w-a,b,1:sum}, we obtain
\begin{equation}
    S[-2] u_P = \frac12 
    \sum_{(a,b) \neq (0,0)} (J_{a,b,P}(t))_P \, (J^{a,b,P}(t))_P u_P.
\label{Su-0}
\end{equation}
Substituting \eqref{Su-0} into $S[-2] u_P \tensor v$ and
swapping $J_{a,b,P}(t)$ and then $J^{a,b,P}(t)$ by the Ward identity
\eqref{ward-cc}, we obtain
\begin{equation}
    S[-2] u_P \tensor v
    \equiv
    u_P \tensor
    \frac12 \sum_{i,j=1}^L \sum_{(a,b)\neq(0,0)}
    w_{a,b}(z_i - u) w_{-a,-b}(z_j - u)
    \rho_i(J_{a,b}) \rho_j(J^{a,b}) v,
\label{S=tau}
\end{equation}
which proves the lemma because of \eqref{w-a,b:property}.
\end{pf}

\begin{cor}
\label{commutativity-gaudin}
\( %
    [\hat\tau(u), \hat\tau(u')] = 0
\) %
for any $u$ and $u'$. In particular, $H_i$ $(i=0, 1, \ldots, L)$ commute
with each other.
\end{cor}
\begin{pf}
Since the Sugawara operators $S[m]$ commute with the affine Lie algebra
at the critical level due to \eqref{[Sm,An]}, we have
\begin{equation*}
  A[n]S[m]u_P = 0
  \quad
  \text{for $A\in\g$ and $n\geqq 0$.}
\end{equation*}
Hence we can find the following formula in the similar way as the proof of
\lemref{gaudin-as-sug-correl}:
\begin{equation*}
  S[-2]u_{P'}\tensor S[m]u_P\tensor v
  \equiv
  u_{P'}\tensor S[m]u_P\tensor \hat\tau(u')v,
\end{equation*}
where $t(P)=u$, $t(P')=u'$, and $v\in V$.  Using this formula and
\lemref{gaudin-as-sug-correl}, we obtain
\begin{align*}
  u_{P'}\tensor u_P\tensor \hat\tau(u)\hat\tau(u')v
  & \equiv u_{P'}\tensor S[-2]u_P\tensor \hat\tau(u')v
  \equiv S[-2]u_{P'}\tensor S[-2]u_P\tensor v
  \\
  & \equiv S[-2]u_{P'}\tensor u_P\tensor \hat\tau(u)v
  \equiv u_{P'}\tensor u_P\tensor \hat\tau(u')\hat\tau(u)v.
\end{align*}
This proves the corollary in view of \propref{restriction}.
\end{pf}

Once the correspondence of the Hamiltonians of the XYZ Gaudin model and
the correlation functions of the twisted WZW model is established, the
eigenvalue problem of the XYZ Gaudin model is rewritten in terms of the
conformal block of the twisted WZW model, as is the case with the (XXX)
Gaudin model. (See \cite{fre:95}.) We sketch below how it goes,
restricting ourselves to $sl_2$ case. For general $sl_N$ case, we should
introduce higher order Sugawara operators, whose constructions are found
in \cite{hayashi} and \cite{go-wa}.

Let us introduce a meromorphic (single-valued) function on $X$ of the
form
\begin{equation}
    q(t)
    = \sum_{i=1}^L l_i(l_i + 1) \wp(t - z_i)
    + \sum_{i=1}^L \mu_i \zeta(t - z_i) + \mu_0,
\label{def:q}
\end{equation}
where $l_i = (\dim V_i - 1)/2$ (cf.\ \exref{ex:gaudin:N=2}), 
$\mu_i$ and $\mu_0$ are parameters satisfying
$\sum_{i=1}^L \mu_i = 0$. Let $q_i(t-z_i) = \sum_{n\in\Z}
q_{i,n}(t-z_i)^{-n-2}$ be the Laurent expansion of $q(t)$ at $Q_i$.
Denote by $K^{q_i}(V_i)$ the submodule of $\Weyl_{-2}(V_i)$ generated by
the vectors $(S[m]-q_{i,m})v_i$ for $v_i\in V_i$, $m\in\Z$ and 
put $\Weyl^{q_i}(V_i):=\Weyl_{-2}(V_i)/K^{q_i}(V_i)$.

\begin{thm}
\label{cft=gaudin}
The space of conformal coinvariants and that of conformal blocks associated to
the module $\Weyl^q(V) := \bigotimes_{i=1}^L \Weyl^{q_i}(V_i)$ are isomorphic
to the quotient of $V := \bigotimes_{i=1}^L V_i$ by the subspace $J^\mu(V)$
spanned by vectors of the form $(H_i - \mu_i)v$ for $i = 0,1,\ldots,L$ and 
$v\in V$ and its dual: 
\begin{equation*}
  \CC_k(\Weyl^q(V)) \simeqq V/J^\mu(V),
  \quad
  \CB_k(\Weyl^q(V)) \simeqq  (V/J^\mu(V))^\ast.
\end{equation*}
\end{thm}
\begin{pf}
We prove the statement for the conformal blocks. The statement for the
space of conformal coinvariants follows from this since it is
finite-dimensional and dual to the space of conformal blocks.

Let $\Phi$ be any linear functional on 
$\Weyl_{-2}(V) = \bigotimes_{i=1}^L \Weyl_{-2}(V_i)$. A necessary and
sufficient condition for $\Phi$ to be a conformal block in $\CB_k(\Weyl^q(V))$
is that it vanishes on $\gDpr \Weyl_{-2}(V)$ and on the subspaces 
\[ %
  K^q(V) :=
  \sum_{i=1}^L
  \Weyl_{-2}(V_1) \tensor \cdots
  \tensor K^{q_i}(V_i) \tensor \cdots
  \tensor \Weyl_{-2}(V_L).
\] %
First we show that this condition
implies $\Phi((H_i - \mu_i)v) = 0$ for $i = 0, 1, \ldots, L$ and 
$v \in V$. 

The assumption is encapsulated in the following expression by a generating
function,
\begin{equation}
    \sum_{m\in\Z} \Phi(\rho_i(S[m] - q_{i,m}) v) (u-z_i)^{-m-2} = 0,
\label{em-tensor-vanish}
\end{equation}
which means $\Phi(\hat\tau(u)v) = q(u) \Phi(v)$ because of
\eqref{sug-correl-at-P} and \lemref{gaudin-as-sug-correl}. Thus
\eqref{tau-hat} and \eqref{def:q} shows that
$\Phi(H_i v) = \mu_i \Phi(v)$.

We prove the converse statement next. Assume that $\Phi$ vanishes on
the subspace $J^\mu(V)$. 
Let $v$ be an arbitrary vector in $\Weyl_{-2}(V)$. 
We want to show that $\Phi(\rho_i(S[m] - q_{i,m})v)$ vanishes 
for any $m$ and $i$, but for this purpose we may assume $v\in V$
without loss of generality. Indeed
any $v$ can be written in the form $v = g_\prin v^0$ by
the decomposition \eqref{decomp:g-hat}, where
$g_\prin\in U(\gDpr)$, $v^0 \in V$, and therefore
$$
    \rho_i(S[m] - q_{i,m})v 
    =
    g_\prin \rho_i(S[m] - q_{i,m}) v^0,
$$
since $S[m]$ belongs to the center of $U_{-2}(\widehat{sl}(2))$. The Ward
identity \eqref{ward-cb} implies 
that $\Phi(\rho_i(S[m] - q_{i,m})v) = 0$ 
if   $\Phi(\rho_i(S[m] - q_{i,m}) v^0) = 0$.

For $v \in V$, we can prove $\Phi(\rho_i(S[m] - q_{i,m})v) = 0$ 
by tracing back the first part of this proof.
\end{pf}

\section{Sheaves of conformal coinvariants and conformal blocks}
\label{sheaf-conf-block}

So far we have fixed the modulus $\tau$ of an elliptic curve and marked
points on it.
In this section we introduce sheaves of conformal coinvariants and
conformal blocks on a family of pointed elliptic curves.


\subsection{Family of pointed elliptic curves and Lie algebra bundles}
\label{family-of-curves}

In this subsection we construct a family of elliptic curves with marked
points, a group bundle, and the associated Lie algebra bundle over this
family.  The fiber at a point of the base space of the family gives the group
bundle $\Gtw$ and the Lie algebra bundle $\gtw$ on a pointed elliptic curve
defined in \secref{bundles}.

Recall that $\UHP$ denotes the upper half plane.  We define $\Xtilde$ and $S$
by
\begin{align*}
  & S :=
  \{\,(\tau; z)=(\tau;z_1,\ldots,z_L)\in\UHP\times\C^L \mid
      z_i - z_j \not\in \Z + \tau \Z \ \text{if}\ i \ne j \,\},
  \\
  & \Xtilde := S\times\C.
\end{align*}
Let $\tilde\pr=\pr_{\Xtilde/S}$ be the projection from $\Xtilde$ onto $S$
along $\C$ and $\tilde q_i$ the section of $\tilde\pr$ given by $z_i$:
\begin{equation*}
  \tilde q_i(\tau;  z) := (\tau;  z; z_i) \in \Xtilde
  \quad
  \text{for $(\tau;  z)=(\tau; z_1,\ldots,z_L)\in S$}.
\end{equation*}

A {\em family of $L$-pointed elliptic curves} $\pr:\X\onto S$ is
constructed as follows. Define the action of $\Z^2$ on $\Xtilde$ by
\begin{equation}
  (m,n)\cdot(\tau; z;t) := (\tau; z;t + m \tau + n)
  \quad
  \text{for $(m,n)\in\Z^2$, $(\tau; z;t)\in\Xtilde$.}
  \label{def:action-Z2-on-Xtilde}
\end{equation}
Let $\X$ be the quotient space of $\Xtilde$ by the action of $\Z^2$:
\begin{equation}
  \X := \Z^2\backslash\Xtilde.
  \label{def:X}
\end{equation}
Let $\pr_{\Xtilde/\X}$ be the natural projection from $\Xtilde$ onto $\X$ and
$\pr=\pr_{\X/S}$ the projection from $\X$ onto $S$ induced by $\tilde\pr$. 
We put
\begin{equation*}
  q_i:=\pr_{\Xtilde/\X}\circ\tilde q_i, \quad
  Q_i := q_i(S), \quad
  D := \bigcup_{i=1}^L Q_i, \quad
  \dot\X := \X \setminus D.
\end{equation*}
Here $q_i$ is the section of $\pr$ induced by $\tilde q_i$ and $D$ is also
regarded as a divisor $\sum_{i=1}^L Q_i$ on $\X$.  
The fiber of $\pr$ at $(\tau, z)=(\tau;z_1,\ldots,z_L)\in S$
is an elliptic curve with modulus $\tau$ and marked points $z_1,\ldots,z_L$.

A group bundle $\Gtw_\X$ and a Lie algebra bundle $\gtw_\X$ on $\X$ are
defined as follows.  Due to the definition of $\X$, the Galois group of the
covering $\pr_{\Xtilde/\X} : \Xtilde \onto \X$ is naturally identified with
$\Z^2$.  Its natural right action on $\Xtilde$ is given by 
$(\tau; z;t)\cdot(m,n) := (-m,-n)\cdot(\tau; z;t)$. Then the covering
$\pr_{\Xtilde/\X} : \Xtilde \onto \X$ is regarded as a principal $\Z^2$-bundle
on $\X$.  The actions of the Galois group $\Z^2$ on $G$ and $\g$ are defined
by 
\begin{alignat}{2}
  (m,n)\cdot g &:= (\beta^m\alpha^n) g (\beta^m\alpha^n)^{-1}
  \quad &
  & \text{for $g\in G$ and $(m,n)\in\Z^2$},
  \label{def:Gal-action-on-G}
  \\
  (m,n)\cdot A &:= (\beta^m\alpha^n) A (\beta^m\alpha^n)^{-1}
  \quad &
  & \text{for $A\in\g$ and $(m,n)\in\Z^2$}.
  \label{def:Gal-action-on-g}
\end{alignat}
These actions produces the associated group bundle $\Gtw_\X$ and the
associated Lie algebra bundle $\gtw_\X$ on $\X$:
\begin{equation}
  \Gtw_\X := \Xtilde\times^{\Z^2} G,
  \quad
  \gtw_\X := \Xtilde\times^{\Z^2} \g.
  \label{def:GtwX-gtwX}
\end{equation}
Their fibers at a point $(\tau; z)\in S$ can be identified with $\Gtw$ and
$\gtw$ in \secref{bundles}.

We denote the $\O_\X$-Lie algebra associated to the Lie algebra bundle
$\gtw_\X$ by the same symbol $\gtw_\X$, as mentioned in
\secref{def-conf-block}.  The sheaf $\gtw_\X$ can be written in the following
form:
\begin{equation}
  \gtw_\X
  =
  \{\, A \in (\pr_{\Xtilde/\X})_\ast(\g\tensor\O_\Xtilde) \mid
       A(p\cdot\tilde x) = p\cdot A(\tilde x)
       \ \text{for } \tilde x\in\Xtilde,\ p\in\Z^2 \,\}.
  \label{ident:gtwX}
\end{equation}
Hence if we take an open subset $U'$ of $\Xtilde$ which does not intersect
$(m,n)\cdot U'$ for any $(m,n)\in\Z^2\setminus\{(0,0)\}$ and denote
by $U$ the image of $U'$ on $\X$, then the restriction of $\gtw_\X$ on $U$ can
be canonically identified with $\g\tensor\O_U$:
\begin{equation}
  \gtw_\X|_U \simeqq (\pr_{U'/U})_\ast(\g\tensor\O_{U'}) = \g\tensor\O_U,
  \label{triv:gtwX}
\end{equation}
where $\pr_{U'/U}$ is the natural biholomorphic projection $U'\isoto U$.

Denote by $\Omega_\X^1$ the sheaf of 1-forms on $\X$ and by
$\Omega_{\X/S}^1$ the sheaf of relative differentials on $\X$ over $S$.  It
follows from the definition of $\gtw_\X$ that $\gtw_\X$ possesses a natural
connection $\nabla:\gtw_\X\to\gtw_\X\tensor_{\O_\X}\Omega_\X^1$,
which is induced by the trivial connection 
$\id\tensor d: \g\tensor\O_\Xtilde\to\g\tensor\Omega_\Xtilde^1$ through the
identification \eqref{ident:gtwX}.  This means that, under the trivialization
\eqref{triv:gtwX}, the connection $\nabla$ is identified with $\id\tensor d_U$
where $d_U$ is the exterior derivation on $U$. The relative connection
$\nabla_{\X/S}$ along the fibers is defined to be the composite of the
connection $\nabla$ and the natural homomorphism 
$\Omega_\X^1 \to \Omega_{\X/S}^1$.  Under the trivialization \eqref{triv:gtwX}
and the coordinate $(\tau; z;t)$, the relative connection $\nabla_{\Xtilde/S}$
is equal to the exterior derivation by $t$.

Define the invariant $\O_X$-inner product on $\gtw_\X$ by
\begin{equation}
  (A|B)
  := \frac{1}{2\hvee} \tr_{\gtw_\X}(\ad A \ad B)
  = \frac{1}{2N} \tr_{\gtw_\X}(\ad A \ad B) \in \O_X
  \quad
  \text{for $A,B\in\gtw_\X$},
  \label{def:inner-prod-gtwX}
\end{equation}
where the symbol $\ad$ denotes the adjoint representation of the $\O_\X$-Lie
algebra $\gtw_\X$.  Under the trivialization \eqref{triv:gtwX}, the inner
product on $\gtw_\X$ is equal to the inner product defined by
\eqref{def:inner-product} and hence it is invariant under the translation
along the connection $\nabla$.

Recall that $\eps = \exp(2\pi i/N)$.  For $(a,b)\in(\Z/N\Z)^2$, the
1-dimensional representation $(m,n)\mapsto \eps^{bm+an}$ of $\Z^2$ defines
the associated flat line bundle $L_{a,b}$ on $\X$.  We obtain the
decomposition of $\gtw_\X$ into line bundles:
\begin{equation}
  \gtw_\X = \bigoplus_{(a,b)\ne(0,0)} J_{a,b} L_{a,b}.
  \label{decomp:gtwX}
\end{equation}
This is a sheaf version of \eqref{decomp:gtw}.  

\begin{lem} \quad
\label{lem:sheaf-vanishing}
  $R^p \pr_\ast \gtw_\X = 0$ for all $p$.
\end{lem}
\begin{pf}
Since $L_{a,b}^\ast\tensor_{\O_\X}\Omega_{\X/S}^1$ is isomorphic to
$L_{-a,-b}$, it follows that 
$R^1\pr_\ast L_{a,b}\simeqq\HOM_{\O_S}(\pr_\ast L_{-a,-b},\O_S)$ by the Serre
duality.  Therefore, because of the decomposition \eqref{decomp:gtwX}, it is
enough to show $\pr_\ast L_{a,b}=0$ for $(a,b)\ne(0,0)$.
Let $U$ be any open subset of $S$ and put $V:=\pr^{-1}(U)$.  For each
$s=(\tau;z)\in U$, the restriction $L_{a,b}|_{X_s}$ of $L_{a,b}$ on the fiber
$X_s:=\pr^{-1}(s)$ can be identified with the line bundle $L_{a,b}$ on
$X_\tau$ define in \secref{bundles}.  Hence we obtain
$H^0(X_s,L_{a,b}|_{X_s})=0$ for each $s\in S$ In particular, for every 
$f\in H^0(V,L_{a,b})=H^0(U,\pr_\ast L_{a,b})$, the restriction $f|_{X_s}$ of
$f$ on the fiber vanishes for each $s\in S$ and hence $f$ vanishes itself.
This means that $H^0(U,\pr_\ast L_{a,b})=0$.  We have proved the lemma.
\end{pf}

\subsection{Sheaf of affine Lie algebras}
\label{definition-of-sheaf-aff-alg}

In this section we define a sheaf version of the Lie algebras $\hat\g^D$,
$\gDpr$, etc.\ on the base space $S$ of the family.

For an $\O_\X$-module $\F$ and a closed analytic subset $W$ of $\X$, the
restriction $\F|_{W}$ of $\F$ on $W$ and the completion
$\hat\F_{|W}=(\F)_W^{\wedge}$ of $\F$ at $W$ are defined by
\begin{equation}
  \F|_{W} := \F/I_W\F,
  \quad
  \widehat\F_{|W} = (\F)_W^{\wedge} := \projlim_{n\to\infty}(\F/I_W^n\F),
  \label{def:restrection-completion}
\end{equation}
where $I_W$ is the defining ideal of $W$ in $\X$.

We define the $\O_S$-Lie algebras $\g_S^{Q_i}$, $\g_{S,+}^{Q_i}$, $\g_S^D$,
$\g_{S,+}^D$, and $\gDPr$ as follows:
\begin{equation}
\begin{alignedat}{2}
  & \g_S^{Q_i} 
  := \pr_\ast (\gtw_\X(\ast Q_i))_{Q_i}^{\wedge},
  \quad &
  & \g_{S,+}^{Q_i}
  := \pr_\ast (\gtw_\X)_{Q_i}^{\wedge},
\\
  & \g_S^D
  := \pr_\ast (\gtw_\X(\ast D))_D^{\wedge}
  = \bigoplus_{i=1}^L \g_S^{Q_i},
  \quad &
  & \g_{S,+}^D
  := \pr_\ast (\gtw_\X)_D^{\wedge}
  = \bigoplus_{i=1}^L \g_{S,+}^{Q_i},
\\
  & \gDPr
  := \pr_\ast (\gtw_\X(\ast D)).
  & &
\end{alignedat}
\label{def:Os-Lie-alg}
\end{equation}
The 2-cocycle of $\g_S^D$ is defined by
\begin{equation}
    \caff(A,B)
    := \sum_{i=1}^L \Res_{Q_i} (\nabla A_i | B_i)
     = \sum_{i=1}^L \Res_{Q_i} (\nabla_{\X/S} A_i | B_i),
  \label{def:localized-cocycle}
\end{equation}
where $A=(A_i)_{i=1}^L, B=(B_i)_{i=1}^L \in \g_S^D$
and $\Res_{Q_i}$ is the residue along $Q_i$.  Using the 2-cocycle
$\caff(\cdot,\cdot)$, we define a central extension $\hat\g_S^D$ of
$\g_S^D$:
\begin{equation}
  \hat\g_S^D := \g_S^D \oplus \O_S \khat,
  \label{def:sheaf-aff-alg}
\end{equation}
where its Lie algebra structure is defined by the formula similar to
\eqref{def:aff-alg-str}.  Put $\hat\g_S^{Q_i}:=\g_S^{Q_i}\oplus\O_S\khat$,
which is a $\O_S$-Lie subalgebra of $\hat\g_S^D$.  We call $\hat\g_S^D$ 
(resp.\ $\hat\g_S^{Q_i}$) the {\em sheaf of affine Lie algebras at} $D$
(resp.\ $Q_i$).  

The diagonal embedding of $\gDPr$ into $\g_S^D$ is
defined to be the mapping which sends $A\in\gDPr$ to
$(A_i)_{i=1}^L\in\g_S^D$, where each $A_i$ is the image of $A$ given by
the natural embedding 
$\gtw_\X(\ast D)\injto (\gtw_\X(\ast Q_i))_{Q_i}^{\wedge}$.
We identify $\gDPr$ with its image in $\g_S^D$ and $\hat\g_S^D$.
For $A,B\in\gDPr$, we can regard $(\nabla_{\X/S}A|B)$  as an element of
$\pr_*\Omega_{\X/S}^1(\ast D)$. Hence, using the residue theorem, we obtain
that $\caff(A,B)=0$. Thus $\gDPr$, as well as $\g_{S,+}^{Q_i}$, is an
$\O_S$-Lie subalgebra of $\hat\g_S^D$.  

Put $\hat\g^D_{S,+}:=\g^D_S \oplus \O_S\khat$.  Then
\lemref{lem:sheaf-vanishing} implies the sheaf version of
\eqref{decomp:g-hat}.

\begin{lem}
\label{lem:decomp-hatgSD}
  $\hat\g_S^D = \gDPr \oplus \hat\g_{S,+}^D$.
\end{lem}
\begin{pf}
We can calculate $R^p\pr_\ast\gtw_\X$ for $p=0,1$ as the kernel and the
cokernel of $\gDPr\oplus\g^D_{S,+}\to\g^D_S$, which sends
$(a_\Prin;a_+)$ to $a_\Prin -a_+$.  But then \lemref{lem:sheaf-vanishing} means
that both the kernel and the cokernel vanish and hence 
$\gDPr\oplus\g^D_{S,+}=\g^D_S$.  We have proved the lemma.
\end{pf}

Choose any open neighborhood $U'$ of $\widetilde Q_i:=\tilde q_i(S)$ which
does not intersect $(m,n)\cdot U'$ for any $(m,n)\in\Z^2\setminus\{(0,0)\}$.
Then, applying the trivialization \eqref{triv:gtwX} to $U'$, we obtain a
natural isomorphism
\begin{equation}
  \g_S^{Q_i} \simeqq \g\tensor \pr_\ast(\widehat\O_{\X|Q_i}(\ast Q_i)),
  \label{triv:gSQi-1}
\end{equation}
which does not depend on the choice of $U'$ and is defined globally on $S$.
Furthermore, using the coordinate $(\tau;z;\xi_i)$ with $\xi_i=t-z_i$, we
have the following isomorphism defined over S:
\begin{equation}
  \g_S^{Q_i} \simeqq \g\tensor \O_S((\xi_i)).
  \label{triv:gSQi-2}
\end{equation}
Under this trivialization, $\g_{S,+}^{Q_i}$ is identified with
$\g\tensor\O_S[[\xi_i]]$ and the connections on $\g_S^{Q_i}$ induced by
$\nabla$ and $\nabla_{\X/S}$ are written in the following forms:
\begin{align}
  & \nabla A
  = \pd{A}{\tau}\,d\tau
  + \sum_{i=1}^L\pd{A}{z_i}\,dz_i
  + \pd{A}{\xi_i}\,d\xi_i
  \in
  \g\tensor\Omega_S^1((\xi_i)) \oplus \g\tensor\O_S((\xi_i))\,d\xi_i,
  \label{eq:triv-conn-gtw}
  \\
  & \nabla_{\X/S} A
  = \pd{A}{\xi_i}\,d\xi_i
  \in
  \g\tensor\O_S((\xi_i))\,d\xi_i,
\end{align}
where $A\in \g\tensor\O_S((\xi_i))\simeqq \g_S^{Q_i}$.  We also obtain the
induced global trivialization of the sheaf of affine Lie algebras on $S$:
\begin{equation}
  \hat\g_S^D
  \simeqq \bigoplus_{i-1}^L \g\tensor\O_S((\xi_i))\oplus \O_S \khat.
  \label{triv:hatgSD}
\end{equation}
Under this trivialization, the bracket of $\hat\g_S^D$ is represented in
the following form:
\begin{equation}
  \bigl[(A_i\tensor f_i)_{i=1}^L,
        (B_i\tensor g_i)_{i=1}^L \bigr]
  = \bigl([A_i,B_i]\tensor f_i g_i\bigr)_{i=1}^L
  + \khat \sum_{i=1}^L(A_i|B_i)
      \Res\limits_{\xi_i=0} (df_i\cdot g),
  \label{triv:bracket-hatgSD}
\end{equation}
where $A_i,B_i\in\g$ and $f_i,g_i\in\O_S((\xi_i))$.

\subsection{Definition of the sheaves of conformal coinvariants and 
  conformal blocks}
\label{definition-of-sheaf-conf-block}

For any $\O_S$-Lie algebra ${\frak a}=\g_{S,+}^D,\hat\g_S^D,\text{etc.}$, we
denote by $U_S({\frak a})$ the universal $\O_S$-enveloping algebra of
$\frak a$ and define the category of $\frak a$-modules to be that of
$U_S({\frak a})$-modules.

\begin{defn}
\label{def:sCC-sCB}
For any $\hat\g_S^D$-module $\M$, we define 
the {\em sheaf $\sCC(\M)$ of conformal coinvariants} and 
the {\em sheaf $\sCB(\M)$ of conformal blocks} by
\begin{align}
  \sCC(\M)
  & := \M/\gDPr\M,
  \label{def:sheaf-conf-coinvariant}
  \\
  \sCB(\M)
  & := \HOM_{\O_S}(\sCC(\M), \O_S)
  \label{def:sheaf-conf-block}.
\end{align}
Namely, the $\O_S$-module $\sCC(\M)$ is generated by $\M$ with relations
\begin{equation}
    A_\Prin v \equiv 0
\label{sheaf-ward-cc}
\end{equation}
for all $A_\Prin \in \gDPr$, $v \in \M$,
and $\Phi\in\sCB(\M)$ means that $\Phi$ belongs to $\HOM_{\O_S}(\M,\O_S)$
and satisfies
\begin{equation}
    \Phi(A_\Prin v) = 0
\label{sheaf-ward-cb}
\end{equation}
for all $A_\Prin \in \gDPr$, $v \in \M$.
These equations, \eqref{sheaf-ward-cc} and \eqref{sheaf-ward-cb}, are
also called the {\em Ward identities}.
We can regard $\sCC(\cdot)$ as a covariant right exact functor from the
category of $\hat\g_S^D$-modules to that of $\O_S$-modules and similarly
$\sCB(\cdot)$ as a contravariant left exact functor.
\end{defn}

The $\hat\g_S^D$-modules of our concern are the sheaf version $\sWeyl_k(V)$ of
$\Weyl_k(V)$ in \secref{def-conf-block}. We give two
equivalent definitions of $\sWeyl_k(V)$.

\subsubsection*{First definition of $\sWeyl_k(V)$}
Fix an arbitrary complex number $k$. For each $i=1,\ldots,L$, let $M_i$ be a
representation with level $k$ of the affine Lie algebra 
$\ghat_i := \g\tensor\C((\xi_i))\oplus\C \khat$.  (Here $\ghat_i$-modules are
said to be of level $k$ if the canonical central element $\khat$ acts on them
as $k\cdot\id$.)  Assume the {\em smoothness} of $M_i$, namely, for any
$v_i\in M_i$, there exists $m\geqq 0$ such that, for
$A_{i_1},\dots,A_{i_\nu}\in\g$, $m_1,\dots,m_\nu\geqq0$, and
$\nu=0,1,2,\ldots$, 
\begin{equation}
  (A_{i_1}\tensor\xi_i^{m_1}\O_S[[\xi_i]])
  \cdots
  (A_{i_\nu}\tensor\xi_i^{m_\nu}\O_S[[\xi_i]])
  v_i=0
  \quad \text{if $m_1+\dots+m_\nu\geqq m$}.
\label{smoothness}
\end{equation}
Put $M:=\bigotimes_{i=1}^L M_i$ and $\M:=M\tensor\O_S$. 
Then $M$ is a representation with level $k$ of the affine Lie algebra
$(\g^{\oplus L})^{\wedge}
:=\bigoplus_{i=1}^L\bigl(\g\tensor\C((\xi_i))\bigr)\oplus\C\khat$ associated
to $\g^{\oplus L}$.  We can define the $\hat\g_S^D$-module structure on $\M$
by
\begin{equation}
  (A_i\tensor f_i(\xi_i))_{i=1}^L (v\tensor a)
  := \sum_{i=1}^L \sum_{m\in\Integer}
  (\rho_i(A_i\tensor \xi_i^m)v) \tensor (f_{i,m}a),
  \quad
  \khat v := kv,
  \label{def:act-hatgSD-M}
\end{equation}
where $A_i\in\g$, $f_i(\xi_i)=\sum_m f_{i,m}\xi_i^m\in\O_S((\xi_i))$,
$f_{i,m},a\in\O_S$, $v\in M$, and $\rho_i(A_i\tensor\xi_i^m)$ denotes the
action of $A_i\tensor\xi_i^m$ on the $i$-th factors in $v$.  If each $M_i$ is
the Weyl module $\Weyl_k(V_i)$ induced up from a finite-dimensional
irreducible representation $V_i$ of $\g$, then we put 
$V:=\bigotimes_{i=1}^L V_i$, $M:=\Weyl_k(V):=\bigotimes_{i=1}^L\Weyl_k(V_i)$,
and $\M:=\sWeyl_k(V):=\Weyl_k(V)\tensor\O_S$ and denote $\sCC(\M)$ and
$\sCB(\M)$ by $\sCC_k(V)$ and $\sCB_k(V)$ respectively.

\subsubsection*{Second definition of $\sWeyl_k(V)$}
Let $V_i$ be a finite-dimensional
irreducible representation of $\g$ and put $V:=\bigotimes_{i=1}^L V_i$.
Denote the constant sheaf associated to $V$ by the same symbol $V$.  Using
the trivialization \eqref{triv:hatgSD}, we can define the action of
$\hat\g_{S,+}=\g_{S,+}^{D}\oplus\O_S\khat$ on $V\tensor\O_S$ by
\begin{equation}
  (A_i\tensor f_i(\xi_i))_{i=1}^L (v\tensor a)
  := \sum_{i}(\rho_i(A_i) v) \tensor (f_i(0)a),
  \quad
  \khat v := kv
  \label{def:act-gS+Qi-Vi}
\end{equation}
where $A_i\in\g$, $f_i(\xi_i)\in\O_S[[\xi_i]]$, $a\in\O_S$, $v\in V$, and
$\rho_i(A_i)$ is the action of $A_i$ on the $i$-th factors in $v$.  The
$\hat\g_S^D$-module $\Weyl_{S,k}(V)$ induced from $V\tensor\O_S$ is defined by
\begin{equation}
  \Weyl_{S,k}(V)
  := \Ind_{\hat\g_{S,+}^D}^{\hat\g_S^D} (V\tensor\O_S)
  = U_S(\hat\g_S^D) \tensor_{U_S(\hat\g_{S,+}^D)} (V\tensor\O_S)
  \label{def:sheaf-weyl-mod}
\end{equation}
Using the decomposition
\begin{equation*}
  \hat\g_S^D
  = \left( \bigoplus_{i=1}^L \g\tensor\xi_i^{-1}\O_S[\xi_i^{-1}] \right)
  + \hat\g_{S,+}^D,
\end{equation*}
we can show that $\Weyl_{S,k}(V)$ has the following $\O_S$-free basis:
\begin{equation}
  \rho_{i_1}(A_{s_1}[m_1]) \cdots \rho_{i_\nu}(A_{s_\nu}[m_\nu]) v_j,
  \label{basis-bbbVSkV}
\end{equation}
where $\nu=0,1,2,\ldots$, $i_n=1,\dots,L$, $\{A_s\}$ is a basis of $\g$,
$\{v_j\}$ is a basis of $V$, and $m_1 \leqq \cdots \leqq m_\nu < 0$. This is
also an $\O_S$-free basis of $\sWeyl_k(V)$ and hence $\Weyl_{S,k}(V)$ is
isomorphic to $\sWeyl_k(V)$ as a $\hat\g^D_S$-module.  In the following we
identify $\sWeyl_k(V)$ with $\Weyl_{S,k}(V)$. 

This identification of the two definitions and \lemref{lem:decomp-hatgSD}
prove the sheaf version of \propref{restriction}.

\begin{prop}
\label{sheaf-restriction}
Let $V_i$ be a finite-dimensional irreducible representation of $\g$ for each
$i$ and put $V:=\bigotimes_{i=1}^L V_i$.  Then the natural inclusion
$V\tensor\O_S\injto\sWeyl_k(V)$ induces the following isomorphisms:
\begin{equation*}
  \sCC_k(V) \isoot V\tensor\O_S
  \quad\text{and}\quad
  \sCB_k(V) \isoto V^\ast\tensor\O_S.
\end{equation*}
\end{prop}
\begin{pf}
 From the second definitions of $\sWeyl_k(V)$ and
\lemref{lem:decomp-hatgSD}, it follows that
\begin{equation*}
  \sWeyl_k(V)
  = U_S(\hat\g^D_S) \tensor_{U_S(\hat\g^D_{S,+})} (V\otimes\O_S)
  \isoot U_S(\gDPr) \tensor_{\O_S} (V\otimes\O_S).
\end{equation*}
Namely, $\sWeyl_k(V)$ is freely generated by $V\tensor\O_S$ over
$U_S(\gDPr)$. 
Hence we obtain the formulae
\begin{align*}
  \sCC_k(V)
  & = \sWeyl_k(V)/\gDPr\sWeyl_k(V)
  \isoot \bigl( U_S(\gDPr)/\gDPr U_S(\gDPr) \bigr)
  \tensor_{\O_S} (V\otimes\O_S)
  \\
  & \isoot \O_S \otimes_{\O_S} (V\otimes\O_S) = V\otimes\O_S,
  \\
  \sCB_k(V)
  & = \HOM_{\O_S}(\sCC_k(V),\O_S) 
  \isoto \HOM_{\O_S}(V\otimes\O_S, \O_S) 
  = V^\ast\otimes\O_S.
\end{align*}
We have completed the proof of the proposition.
\end{pf}

\begin{cor}
\label{sheaf-coherent}
For each $i=1,\ldots,L$, let $V_i$ be a finite-dimensional irreducible
representation of $\g$ and $M_i$ a quotient module of the generalized Verma
module $\Weyl_k(V_i)$ of the affine Lie algebra $\hat\g_i$.  Put
$M:=\bigotimes_{i=1}^L M_i$ and $\M:=M\tensor\O_S$.  Then the sheaf $\sCC(\M)$
of conformal coinvariants and the sheaf $\sCB(\M)$ of conformal blocks are
$\O_S$-coherent.
\end{cor}
\begin{pf}
Since $\sCB(\M)=\HOM_{\O_S}(\sCC(\M),\O_S)$, it suffices for the proof to see
that $\sCC(\M)$ is coherent. The right exactness of the functor $\sCC(\cdot)$
and the fact that $\M$ is a quotient $\hat\g_S^D$-module of
$\sWeyl_k(V)$ imply that $\sCC(\M)$ is a quotient $\O_S$-module of
$\sCC_k(V)=\sCC(\sWeyl_k(V))$, which is $\O_S$-coherent due to
\propref{sheaf-restriction}.  Hence $\sCC(\M)$ is also $\O_S$-coherent.
\end{pf}

\section{Sheaf of the Virasoro algebras}
\label{sheaf-of-virasoro}

This section provides the sheaf version of the Virasoro algebras and
its actions on representations of the sheaf of affine Lie algebra, which
will be used in \secref{flat-connections} to endow the sheaf of conformal
coinvariants and the sheaf of conformal blocks with $\D_S$-module
structures, when the level is not critical (i.e., $\kappa=k+\hvee\ne0$).

\subsection{Definition of the sheaf of the Virasoro algebras}
\label{def-sheaf-virasoro}

We define the {\em sheaf of the Virasoro algebras} by
\begin{equation}
  \Vir^D_S := 
  \T_S\oplus\bigoplus_{i=1}^L\O_S((\xi_i))\partial_{\xi_i}\oplus\O_S\chat,
  \label{def:sheaf-vir}
\end{equation}
where $\T_S$ is the tangent sheaf of $S$. The Lie algebra structure which
we shall give to this $\O_S$-sheaf below reduces to the Virasoro
algebra structure on $\Vir^D$ \eqref{def:Vir-D}, when $S$ is replaced with
a point.

In order to define a Lie algebra structure on $\Vir^D_S$, we introduce the
following notation:
\begin{itemize}
\item For $\mu,\nu\in\T_S$, the symbol
  $[\mu,\nu]$ denotes the natural Lie bracket in $\T_S$;

\item For $\theta=(\theta_i)_{i=1}^L, \eta=(\eta_i)_{i=1}^L \in
  \bigoplus_{i=1}^L\O_S((\xi_i))\partial_{\xi_i}$, the symbol
  $[\theta,\eta]^0 = ([\theta_i,\eta_i]^0)_{i=1}^L$ denotes the natural
  Lie bracket in $\bigoplus_{i=1}^L\O_S((\xi_i))\partial_{\xi_i}$ given by
  \begin{equation*}
    [\theta_i(\xi_i)\partial_{\xi_i},\eta_i(\xi_i)\partial_{\xi_i}]^0
    =
    \bigl( \theta_i(\xi_i)\eta_i'(\xi_i)
         - \eta_i(\xi_i)\theta_i'(\xi_i) \bigr) \partial_{\xi_i}.
  \end{equation*}

\item $\cvir(\theta,\eta) := \sum_{i=1}^L
  \Res_{\xi_i=0}\bigl(\theta_i'''(\xi_i)\eta_i(\xi_i)\,d\xi_i\bigr)$.
  (The symbol $\cvir$ stands for ``Cocycle defining the Virasoro
  algebra''.);
 
\item For $\theta\in\bigoplus_{i=1}^L\O_S((\xi_i))\partial_{\xi_i}$ and
  $f\in\O_S$, the symbols
  $\mu(\theta)$ and $\mu(f)$ denote the natural actions of a vector field
  $\mu\in\T_S$ on $\bigoplus_{i=1}^L\O_S((\xi_i))\partial_{\xi_i}$
  and $\O_S$ respectively.
\end{itemize}
We define the Lie algebra structure on $\Vir^D_S$ by
\begin{equation}
\begin{aligned}
  & [(\mu;\theta;f\chat), (\nu;\eta;g\chat)]
  \\
  & :=
  ( [\mu, \nu];
    \mu(\eta) - \nu(\theta) + [\theta,\eta]^0;
    (\mu(g) - \nu(f) + \cvir(\theta,\eta))\chat ),
\end{aligned}
\label{def:Lie-alg-str-sheaf-vir}
\end{equation}
where $(\mu;\theta;f\chat), (\nu;\eta;g\chat)\in\Vir^D_S$.
Note that
$\Vir^D_S$ is not an $\O_S$-Lie algebra but a $\C_S$-Lie algebra.  We call
$\Vir^D_S$ the {\em sheaf of Virasoro algebras on $S$}.

\begin{rem}
Later representations of $\Vir^D_S$ shall be interpreted as
representations of an extension $\VirDPr$ of $\T_S$ defined below and thus 
shall be given a $\D_S$-module structure.
\end{rem}

Let $\T_\X$ denote the tangent sheaf of the total space $\X$ and
$\T_{\X/S}$ the relative tangent sheaf of the family $\pr:\X\to S$ (i.e., the
sheaf of vector fields along the fibers of $\pr$ on $\X$).  
Since $\pr:\X\to S$ is smooth, we have the following short exact sequence:
\begin{equation*}
  0 \to \T_{\X/S}(\ast D) \to \T_\X(\ast D) \to (\pr^\ast\T_S)(\ast D) \to 0.
\end{equation*}
Note that $\pr^*\T_S = \O_\X\tensor_{\pr^{-1}\O_S}\pr^{-1}\T_S$ does
{\em not} possess a natural Lie algebra structure, but
$\pr^{-1}\T_S\subset\pr^*\T_S$ does.  Defining $\T_{\X,\pr}(\ast D)$ to be
the inverse image of $\pr^{-1}\T_S$ in $\T_\X(\ast D)$, we obtain the
following Lie algebra extension: 
\begin{equation*}
  0 \to \T_{\X/S}(\ast D) \to \T_{\X,\pr}(\ast D) \to \pr^{-1}\T_S \to 0.
\end{equation*}
The direct image of this sequence by $\pr$ is also exact and gives the
following Lie algebra extension:
\begin{equation}
  0 \to \TDPr \to \VirDPr \to \T_S \to 0,
\label{eq:LieAlgExtOfTSByTDPr}
\end{equation}
where we put
\begin{equation*}
  \VirDPr := \pr_\ast\T_{\X,\pr}(\ast D),
  \quad
  \TDPr := \pr_\ast\T_{\X/S}(\ast D).
\end{equation*}
\begin{rem}
The exact sequence \eqref{eq:LieAlgExtOfTSByTDPr} is essential in the
constructions of connections on the sheaf $\sCC(\M)$ of conformal coinvariants
and the sheaf $\sCB(\M)$ of conformal blocks.  Generally, a connection is
defined to be an action of the tangent sheaf satisfying the certain axioms.
Using the exact sequence \eqref{eq:LieAlgExtOfTSByTDPr}, we can obtain a
connection if we have an actions of $\VirDPr$ whose restriction on $\TDPr$
is trivial.  (cf.\ \lemref{lem:Vir-ActsOnM}, \lemref{lem:Vir-ActsOnsCC},
\lemref{lem:Vir-ActsOnsCB}, \lemref{lem:TDPrKillsCB},
\lemref{lem:TDPrSubsetgDPr}, and \thmref{thm:D-module-str}.)
\end{rem}

\begin{lem}
\label{local-section-of-VirDPr}
For a local section $a_\Prin$ of $\VirDPr=\pr_\ast\T_{\X,\pr}(\ast D)$, its
pull-back $\tilde a_\Prin$ to $\Xtilde$ is of the form:
\begin{equation*}
    \tilde a_\Prin = \mu_0(\tau;z) \partial_\tau 
    + \sum_{i=1}^L \mu_i(\tau;z) \partial_{z_i}
    + \theta^t(\tau;z;t) \partial_t,
\end{equation*}
where $\mu_i=\mu_i(\tau;z)\in\O_S$ and $\theta^t(\tau;z;t)$ is a meromorphic
function globally defined along the fibers of $\tilde\pr$ with the following
properties:
\begin{enumerate}
\item The poles of $\theta^t(\tau;z;t)$ are contained in
  $\pr_{\Xtilde/\X}^{-1}(D)$;
\item The quasi-periodicity:
  \begin{equation}
    \theta^t(\tau;z;t+m\tau+n) = \theta^t(\tau;z;t) + m\mu_0(\tau;z)
    \quad\text{for $(m,n)\in\Z^2$}.
    \label{property-theta-t}
  \end{equation}
\end{enumerate}
\end{lem}

\begin{pf}
Let $(m,n)$ be in $\Z^2$ and $f_{m,n}$ denote the action of $(m,n)$ on
$\Xtilde$ given by $(\tau;z;t)\mapsto(\tau;z;t+m\tau+n)$. Then its
derivative $df_{m,n}$ sends $\partial_\tau$, $\partial_{z_i}$, and
$\partial_t$ to $\partial_\tau+m\partial_t$, $\partial_{z_i}$, and
$\partial_t$ respectively.  Since $\tilde a_\Prin$ induces the vector
field $a_\Prin$ in $\pr_\ast\T_{\X,\pr}(\ast D)$, we obtain a formula
\begin{equation*}
  \mu + \theta^t(t+m\tau+n)\partial_t = df_{m,n}(\tilde a_\Prin),
\end{equation*}
which is equivalent to
\begin{equation*}
  \theta^t(t+m\tau+n) = \theta^t(t) + m\mu_0,
\end{equation*}
which proves the lemma.
\end{pf}

The local section $a_\Prin$ is mapped to
$\mu=\mu_0\partial_t+\sum_{i=1}^L\mu_i\partial_{z_i}\in\T_S$ by the projection
along $\pr$ in \eqref{eq:LieAlgExtOfTSByTDPr} and belongs to
$\TDPr=\pr_\ast\T_{\X/S}(\ast D)$ if and only if $\mu = 0$.  Under the local
coordinate $\xi_i=t-z_i$, the local section $a_\Prin$ is uniquely represented
in the following form:
\begin{equation*}
  a_\Prin
  = \mu + \theta_i(\xi_i)\partial_{\xi_i}
  \in \T_S\oplus \O_S((\xi_i))\partial_{\xi_i},
\end{equation*}
where $\theta_i(\xi_i)$ is the Laurent expansion of $\theta^t(\tau;z;t)$ in
$\xi_i=t-z_i$.  Thus we obtain the embedding of $\VirDPr$ into
$\T_S\oplus\bigoplus_{i=1}^L\O_S((\xi_i))\partial_{\xi_i}\subset\Vir^D_S$
given by
\begin{equation}
  \VirDPr \injto 
  \T_S\oplus\bigoplus_{i=1}^L\O_S((\xi_i))\partial_{\xi_i}\subset\Vir^D_S,
  \quad
  a_\Prin 
  \mapsto (\mu;\theta)=(\mu;(\theta_i(\xi_i)\partial_{\xi_i})_{i=1}^L).
\label{emb:VirDPrD->VirSD}
\end{equation}
We identify $\VirDPr$ with its image in $\Vir^D_S$.  For instance,
$(\mu;\theta;0)\in\VirDPr$ means that $\mu\in\T_S$, 
$\theta\in\bigoplus_{i=1}^L\O_S((\xi_i))\partial_{\xi_i}$, and
$(\mu;\theta;0)$ belongs to the image of $\VirDPr$ in $\Vir^D_S$.  We also
identify the subsheaf $\TDPr\subset\VirDPr$ with its image in $\Vir^D_S$.  

\begin{rem}
These formulations are essentially an application of the Beilinson-Schechtman
theory in \cite{BS} to our situation.  The theory contains a natural
construction of the Kodaira-Spencer map of a family of compact Riemann
surfaces and its generalization to Virasoro algebras.  For a brief sketch, see
\appref{kodaira-spencer}.
\end{rem}

A natural question is whether or not the embeddings of $\VirDPr$
into $\Vir^D_S$ is a Lie algebra homomorphism. 

\begin{lem}
\label{Vir-toVirIsLieAlgHom}
The embedding $\VirDPr\injto\Vir^D_S$ is a Lie
algebra homomorphism.
\end{lem}
\begin{pf}
Let $a_\Prin$ and $b_\Prin$ be in
$\VirDPr=\pr_\ast\T_{\X,\pr}(\ast D)$. Denote their images in $\Vir^D_S$ by
$(\mu;\theta;0)$ and $(\nu;\eta;0)$ and their pull-backs to
$\Xtilde$ by $\tilde a_\Prin$ and $\tilde b_\Prin$ respectively.
It suffices for the proof to show that $\cvir(\theta,\eta)=0$.

Under the coordinate $(\tau;z;t)$ of $\Xtilde$, the vector fields 
$\tilde a_\Prin$ and $\tilde b_\Prin$ are represented as
\begin{equation}
  \tilde a_\Prin 
  = \mu + \theta^t(t)\partial_t,
  \quad
  \tilde b_\Prin 
  = \nu + \eta^t(t)\partial_t,
\label{tilde-a/b-pr}
\end{equation}
where we write $\mu,\nu\in\T_S$ in the following forms:
\begin{equation}
  \mu
  = \mu_0 \partial_\tau
  + \sum_{i=1}^L \mu_i\partial_{z_i},
  \quad
  \nu
  = \nu_0\partial_\tau
  + \sum_{i=1}^L \nu_i\partial_{z_i}.
\label{mu/nu}
\end{equation}
Here we omit the arguments $(\tau;z)$ for simplicity: $\mu_i=\mu_i(\tau;z)$,
$\theta^t(t)=\theta^t(\tau;z;t)$, etc. 
Because of \lemref{local-section-of-VirDPr}, we have
\begin{equation*}
  \theta^t(t+m\tau+n) = \theta^t(t) + m\mu_0,
  \quad
  \eta^t(t+m\tau+n) = \eta^t(t) + m\nu_0.
\end{equation*}
Hence we can define the relative meromorphic 1-form 
$\omega\in\pr_\ast\Omega_{\X/S}^1(\ast D)$ by 
\begin{equation}
  \omega
  = \omega(t)\,dt
  := \pd{^2\theta^t(t)}{t^2} \pd{\eta^t(t)}{t} \,dt.
\label{def:omega-relative}
\end{equation}
Here, the well-definedness of $\omega$ as a 1-form in
$\pr_\ast\Omega_{\X/S}^1(\ast D)$ follows from the fact that the definition
of $\omega(t)$ implies $\omega(t+m\tau+n)=\omega(t)$ for $m,n\in\Z$.

On the other hand, under the local coordinate $(\tau;z;\xi_i)$ of $\X$ around
$Q_i$ given by $\xi_i=t-z_i$, the vector fields $a_\Prin$ and $b_\Prin$ are
represented in the following forms:
\begin{equation*}
  a_\Prin 
  = \mu  + \theta_i(\xi_i)\partial_{\xi_i},
  \quad
  b_\Prin 
  = \nu + \eta_i(\xi_i)\partial_{\xi_i},
\end{equation*}
where $\mu$ and $\nu$ are the same as those in \eqref{tilde-a/b-pr} and
$\theta_i$ and $\eta_i$ are given in terms of $\theta^t$, $\eta^t$, $\mu_i$
and $\nu_i$ in \eqref{tilde-a/b-pr} and \eqref{mu/nu} by
\begin{equation*}
  \theta_i(\xi_i)=\theta^t(z_i+\xi_i) - \mu_i,
  \quad
  \eta_i(\xi_i)=\eta^t(z_i+\xi_i) - \nu_i.
\end{equation*}
Thus by \eqref{def:omega-relative}, we have
\begin{equation*}
  \omega =
  \pd{^2\theta_i(\xi_i)}{\xi_i^2}
  \pd{\eta_i(\xi_i)}{\xi_i} \,d\xi_i
  \quad\text{around $Q_i$.}
\end{equation*}
Hence the residue theorem leads to
\begin{equation*}
  \cvir(\theta,\eta)
  =-\sum_{i=1}^L\Res\limits_{\xi_i=0}
   \left(
     \pd{^2\theta_i(\tau;z;\xi_i)}{\xi_i^2}
     \pd{\eta_i(\tau;z;\xi_i)}{\xi_i}
     \,d\xi_i
   \right)
  =-\sum_{i=1}^L\Res_{Q_i}\omega=0.
\end{equation*}
This proves the lemma.
\end{pf}

\begin{rem}
The same question about the embedding $\VirDPr\injto\Vir^D_S$ can be
answered under a more general formulation for higher genus compact Riemann
surfaces with a projective structure.  However, in higher genus case, the
embedding $\VirDPr\injto\Vir^D_S$ is not always a Lie algebra homomorphism.
The case of genus 1 is very special.  See \appref{higher-genus} for a short
sketch of a formulation.
\end{rem}

In order to define the action of $\VirDPr$ on $\hat\g^D_S$, let us introduce
the following notation:
\begin{itemize}
\item $A=(A_i)_{i=1}^L, B=(B_i)_{i=1}^L \in
  \bigoplus_{i=1}^L\g\otimes\O_S((\xi_i))$;

\item $[A,B]^0=([A_i,B_i]^0)_{i=1}^L$ denotes the natural Lie bracket in
  $\g^D_S\simeqq\bigoplus_{i=1}^L\g\otimes\O_S((\xi_i))$ given by the base
  extension of the Lie algebra $\g$;

\item $(A;f\khat),(B;g\khat)\in\hat\g^D_S=\g^D_S\oplus\O_S\khat$;

\item For $\mu \in \T_S$, the symbol $\mu(A)$ denotes the natural actions of
  $\T_S$ on $\g^D_S$;

\item The natural action of $\theta=(\theta_i)_{i=1}^L
  \in\bigoplus_{i=1}^L\O_S((\xi_i))\partial_{\xi_i}$ on $A\in\g^D_S$ is
  defined by 
  \begin{equation*}
    \theta_i(A_i\otimes f_i(\xi_i))
    := A_i\otimes(\theta_i(\xi_i)f_i'(\xi_i))
    \quad
    \theta(A):=(\theta_i(A_i\tensor f_i))_{i=1}^L,
  \end{equation*}
  where $\theta_i = \theta_i(\xi_i)\partial_{\xi_i}$ and 
  $A=A_i\otimes f_i(\xi_i)\in\g\otimes\O_S((\xi_i))$.
\end{itemize}
Then the action of $\VirDPr$ on $\hat\g^D_S$ is defined by
\begin{equation}
  (\mu;\theta;0)\cdot(A;g\khat) := [(\mu;\theta;0), (A;g\khat)]
  \quad\text{for $(\mu;\theta;0)\in\VirDPr$ and $(A;g\khat)\in\hat\g^D_S$,}
  \label{def:Vir-ActsOnHatg}
\end{equation}
where $\VirDPr$ is identified with its image in $\Vir^D_S$ by
\eqref{emb:VirDPrD->VirSD} and the
bracket of the right-hand side is a Lie bracket in the semi-direct product
Lie algebra $\Vir^D_S\ltimes\hat\g^D_S$ defined by
\begin{equation}
  [(\mu;\theta;f\chat), (A;g\khat)]
  := (\mu(A)+\theta(A); \mu(g)\khat)
  \quad
  \text{for $(\mu;\theta;f\chat)\in\Vir^D_S$ and $(A;g\khat)\in\hat\g^D_S$.}
\label{def:semi-direct-prod}
\end{equation}

\begin{lem}
\label{lem:Vir-ActsOngDS-}
The action of $\VirDPr$ on $\hat\g^D_S$ preserves $\gDPr$. 
\end{lem}
\begin{pf}
Because of \eqref{eq:triv-conn-gtw}, under the identifications above, the
restriction of the action of $\VirDPr=\pr_\ast\T_{\X,\pi}(\ast D)$ on
$\gDPr$ comes from the action of $\T_{\X,\pi}(\ast D)$ on $\gtw_\X(\ast
D)$ via the connection $\nabla$ on $\gtw_\X$. Namely, if
$a_\Prin\in\VirDPr$, $A_\Prin\in\gDPr$, and their images in $\Vir^D_S$
and $\hat\g^D_S$ are denoted by $(\mu;\theta;0)$ and $(A;0)$ respectively,
then
\begin{equation*}
  [(\mu;\theta;0), (A;0)] 
  = 
  \left(
    \text{the image of $\nabla_{a_\Prin}A_\Prin\in\pr_\ast\gtw(\ast D)$}
  \right).
\end{equation*}
Thus we obtain $[\VirDPr,\gDPr]\subset\gDPr$.
\end{pf}

\subsection{Action of the sheaf of Virasoro algebras}
\label{action-of-vir-sheaf}

In this subsection we define an action of the Lie algebra
$\Vir^D_S\ltimes\hat\g^D_S$ on $\hat\g^D_S$ modules. 

Fix an arbitrary complex number $k$. For each $i=1,\ldots,L$, let $M_i$ be a
representation with level $k$ of the affine Lie algebra $\ghat_i$ satisfying
the smoothness condition \eqref{smoothness}. Put
$M:=\bigotimes_{i=1}^L M_i$ and $\M:=M\tensor\O_S$.  Then $M$ is a
representation with level $k$ of the affine Lie algebra
$(\g^{\oplus L})^\wedge=
\bigoplus_{i=1}^L\bigl(\g\tensor\C((\xi_i))\bigr)\oplus\C\khat$
and $\M$ is a $\hat\g^D_S$-module.

The Sugawara operators $S[m]$ acting on $M_i$'s are given by the formula
\eqref{def:sugawara} and its action on the $i$-th factor $M_i$ in $M$ is
denoted by $\rho_i(S[m])$. 
Define the Sugawara tensor field by
\begin{equation*}
  S(\xi) (d\xi)^2 := \sum_{m\in\Z} \xi^{-m-2} S[m] \,(d\xi)^2,
\end{equation*}
and its action on $M_i$ is denoted by $\rho_i(S(\xi_i)) (d\xi)^2$.  Then,
by the same way as \lemref{global-correl}, we can prove the following lemma.

\begin{lem}
\label{lem:sheaf-sug-correl}
For any $s\in S$, $\Phi\in\sCB(\M)_s$ and $v\in\M_s$, there exists a unique
$\omega\in(\pi_\ast\Omega_{\X/S}^2(\ast D))_s$ such that the expression of
$\omega$ under the coordinate $\xi_i$ coincides with
$\Phi(\rho_i(S(\xi_i))v)\,(d\xi_i)^2$ for each $i=1,\dots,L$.
\end{lem}

\begin{def}
\label{def:sheaf-sug-correl}
We denote $\omega$ in \lemref{lem:sheaf-sug-correl} by
$\Phi(S(\xi)v)\,(d\xi)^2$ or $\Phi(S(P)v)\,(dP)^2$, which is called a
{\em correlation function of the Sugawara tensor} $S(\xi)$ and $v$ under
$\Phi$, or a {\em Sugawara correlation function} for short.
\end{def}

Assume that $\kappa=k+\hvee\ne0$ and put $c_k:=k\dim\g/\kappa$.  Define the
Virasoro operators $T[m]$ and the energy-momentum tensor $T(\xi)$ by
\begin{equation}
  T[m] := \kappa^{-1} S[m],
  \quad
  T(\xi)(d\xi)^2 := \kappa^{-1} S(\xi) (d\xi)^2
  \label{def:Tm}
\end{equation}
as in \eqref{def:virasoro} and the {\em energy-momentum correlation function}
$\Phi(T(\xi)v)\,(d\xi)^2$ to be $\kappa^{-1}\Phi(S(\xi)v)\,(d\xi)^2$ as in
\eqref{def:em-correl}. 

The action $\rho_i(T[m])$ of the Virasoro operators on $M_i$ defines a
representation of the Virasoro algebra with central charge
$c_k=k\dim\g/\kappa$ (\lemref{prop:vec-field-action}). For
$v_i\tensor g\in M_i\tensor\O_S$ and
$\theta_i=\sum_{m\in\Z}\theta_{i,m}\xi_i^{m+1}\partial_{\xi_i}
\in\O_S((\xi_i))\partial_{\xi_i}$, put 
\begin{equation}
  \rho_i(T\{\theta_i\})(v_i\tensor g)
  = \sum_{m\in\Z} (\rho_i(-T[m]) v_i) \tensor (\theta_{i,m}g).
  \label{def:Ttheta1}
\end{equation}
For example, $\rho_i(T\{\xi_i^{m+1}\partial_{\xi_i}\})=\rho_i(-T[m])$.
For $\theta=(\theta_i)_{i=1}^L
\in\bigoplus_{i=1}^L\O_S((\xi_i))\partial_{\xi_i}$, the operator
$T\{\theta\}$ acting on $\M$ is defined by
\begin{equation}
  T\{\theta\} := \sum_{i=1}^L \rho_i(T\{\theta_i\}),
  \label{def:Ttheta2}
\end{equation}
where we consider $\rho_i(T\{\theta_i\})$ as an operator acting on the $i$-th
factor in $\M$.

Define the action of $(\mu;\theta;f\chat)\in\Vir^D_S$ on $\M$ by
\begin{equation}
  (\mu;\theta;f\chat)\cdot (v\tensor g)
  := v\tensor\mu(g) + T\{\theta\}(v\tensor g) + c_k v\tensor(fg)
\label{def:VirActionOnM}
\end{equation}
for $v\tensor g\in\M=M\tensor\O_S$.
The dual actions on $\M^\ast:=\HOM_{\O_S}(\M,\O_S)$ are defined by
\begin{align}
  & (\mu\Phi)(v) := \mu(\Phi(v)) - \Phi(\mu(v)),
  \label{def:vec-field-OnM*}
  \\
  & (\rho^\ast_i(T^\ast\{\theta_i\}\Phi))(v) := - \Phi(\rho_i(T\{\theta_i\})v),
  \label{def:Ti*(vec-field)OnM*}
  \\
  & T^\ast\{\theta\} := \sum_{i=1}^L \rho^\ast_i(T\{\theta_i\}),
  \label{def:T*(vec-field)OnM*}
  \\
  & ((\mu;\theta;f\chat)\cdot\Phi)(v)
  := \mu(\Phi(v)) - \Phi((\mu;\theta;f\chat)\cdot v),
  \label{def:VirActionOnM*}
\end{align}
where $\Phi\in\M^\ast$, and $v\in\M$. Since we have
\begin{equation}
  (\mu;\theta;f\chat)\cdot\Phi
  = \mu\Phi + T^\ast\{\theta\}\Phi - c_k f\Phi,
\label{VirActionOnM*}
\end{equation}
the action of $\Vir^D_S$ on $\M^\ast$ defines a representation of $\Vir^D_S$
with central charge $-c_k$.

The Virasoro operators $T[m]$ satisfy the commutation relations
\eqref{comm-rel:virasoro}. Therefore a straightforward calculation
proves the following lemma.

\begin{lem}
\label{lem:VirHatgActsOnM}
The action of $\hat\g^D_S$ \eqref{def:act-hatgSD-M} and that of $\Vir^D_S$
\eqref{def:VirActionOnM} on $\M$ induce a representation of the
Lie algebra $\Vir^D_S\ltimes\hat\g^D_S$, whose semi-direct product Lie
algebra structure is given by \eqref{def:semi-direct-prod}.
\end{lem}

Define the actions of $\VirDPr$ on $\M$ and $\M^\ast$ through the embedding
$\VirDPr\injto\Vir^D_S$ and the actions of $\Vir^D_S$. Then
\lemref{Vir-toVirIsLieAlgHom} and \lemref{lem:VirHatgActsOnM} immediately lead
to the following lemma.

\begin{lem}
\label{lem:Vir-ActsOnM}
These actions of $\VirDPr$ on $\M$ and $\M^\ast$ are
representations of the Lie algebra $\VirDPr$.
\end{lem}

\begin{lem}
\label{lem:Vir-ActsOnsCC}
The action of $\VirDPr$ on $\M$ preserves $\gDPr\M$ and hence defines a
representation on $\sCC(\M)$ of the Lie algebra $\VirDPr$.
\end{lem}
\begin{pf}
Assume that $\alpha_\Prin\in\VirDPr$, $A_\Prin\in\gDPr$, and $v\in\M$.
\lemref{lem:VirHatgActsOnM} implies that
\begin{equation*}
  \alpha_\Prin A_\Prin v
  = [\alpha_\Prin, A_\Prin]v + A_\Prin \alpha_\Prin v,
\end{equation*}
and \lemref{lem:Vir-ActsOngDS-} means that
$[\alpha_\Prin,A_\Prin]\in\gDPr$.  Hence $\alpha_\Prin A_\Prin v$ belongs
to $\gDPr\M$ and $\VirDPr\gDPr\M$ is included in $\gDPr\M$.  From
\lemref{lem:Vir-ActsOnM} it follows that the induced action of $\VirDPr$
on $\sCC(\M)=\M/\gDPr\M$ defines a representation of the Lie algebra
$\VirDPr$.
\end{pf}

As a result of \lemref{lem:Vir-ActsOnM}, \lemref{lem:Vir-ActsOnsCC} and
\eqref{def:VirActionOnM*}, we obtain the following lemma.

\begin{lem}
\label{lem:Vir-ActsOnsCB}
The action of $\VirDPr$ on $\M^\ast$ preserves the subsheaf $\sCB(\M)$ of
$\M^\ast$ and defines a representation on $\sCB(\M)$ of the Lie algebra
$\VirDPr$. 
\end{lem}

The actions of $\TDPr$ on $\M$ and $\M^*$ are also defined through the
embedding $\TDPr\injto\Vir^D_S$.  Then, as in the proof of
\lemref{global-vec-field-action}, we can show the following lemma from
\lemref{lem:sheaf-sug-correl} thanks to the existence of the
energy-momentum correlation function \eqref{def:Tm}.

\begin{lem}
\label{lem:TDPrKillsCB}
The action of $\TDPr$ on $\M^\ast$ satisfies $\TDPr\cdot\sCB(\M)=0$.
\end{lem}

Using the exact sequence \eqref{eq:LieAlgExtOfTSByTDPr},
\lemref{lem:Vir-ActsOnsCB}, and \lemref{lem:TDPrKillsCB}, we can construct a
flat connection on the sheaf $\sCB(\M)$ of conformal blocks in
\secref{flat-connections} (\remref{rem:desc-Dast}).  However, for the
construction of a flat connection on the sheaf $\sCC(\M)$ of conformal
coinvariants, we shall need the following lemma, as well as the exact sequence
\eqref{eq:LieAlgExtOfTSByTDPr} and \lemref{lem:Vir-ActsOnsCC}.

\begin{lem}
\label{lem:TDPrSubsetgDPr}
The action of $\TDPr$ on $\M$ satisfies $\TDPr\M\subset\gDPr\M$.
\end{lem}

We remark that \lemref{lem:TDPrSubsetgDPr} implies \lemref{lem:TDPrKillsCB},
but the converse does not hold.  The key point in the proof of
\lemref{lem:TDPrKillsCB} is the notion of the energy-momentum correlation
function, which is not useful for the proof of \lemref{lem:TDPrSubsetgDPr}.
Hence we must find a direct proof of \lemref{lem:TDPrSubsetgDPr} without using
the energy-momentum correlation functions.  The rest of this subsection is
devoted to the proof of this lemma along the course similar to that of
\cite{tsuchimoto:93}.

We define the $\O_S$-inner product $(\,.\,,.\,)$ on
$\g^D_S\simeqq\bigoplus_{i=1}^L\g\otimes\O_S((\xi_i))$ by
\begin{equation}
  \bigl(
    (A_i\tensor f_i)_{i=1}^L , (B_i\tensor g_i)_{i=1}^L
  \bigr)
  =
  \sum_{i=1}^L (A_i|B_i) \Res\limits_{\xi_i=0}(f_ig_i\,d\xi_i)
  \label{def:inn-prod-gDS}
\end{equation}
for $A_i,B_i\in\g$ and $f_i,g_i\in\O_S((\xi_i))$.  This inner product is
non-degenerate and allows us to regard $\g^D_S$ as the topological dual of
itself under the $\xi_i$-adic topologies.  Putting 
\begin{align}
  & R := \{\, (a,b,m,i) \mid
    (a,b)\in(\Z/N\Z)^2\setminus\{(0,0)\},\;
    m\in\Z,\;
    i=1,\ldots,L
  \,\},
  \label{def:R}
  \\
  & e_{a,b,i}^m
  = (e_{a,b,i,j}^m(\xi_j))_{j=1}^L
  := (\delta_{i,j}J_{a,b}\tensor\xi_j^m)_{j=1}^L,
  \\
  & e^{a,b,i}_{m}
  = (e^{a,b,i}_{m,j}(\xi_j))_{j=1}^L
  := (\delta_{i,j}J^{a,b}\tensor\xi_j^{-m-1})_{j=1}^L,
\end{align}
we obtain the following topological dual $\O_S$-bases of $\g^D_S$ with
respect to the inner product:
\begin{equation*}
  F_0 := \{\, e_{a,b,i}^m \mid (a,b,m,i)\in R \,\},
  \quad
  F^0 := \{\, e^{a,b,i}_m \mid (a,b,m,i)\in R \,\}.
\end{equation*}
For $A=(A_i)_{i=1}^L,B=(B_i)_{i=1}^L\in\g^D_S$, we introduce the
following notation:
\begin{align}
  & \rho(A) := \sum_{i=1}^L \rho_i(A_i),
  \\
  & \NP\rho(A)\rho(B)\NP
  := \sum_{i=1}^L \rho_i(\NP A_i B_i \NP)
   + \sum_{i\ne j} \rho_i(A_i)\rho_j(B_j).
\end{align}
Recall that, for $\theta=(\theta_i(\xi_i)\partial_{\xi_i})_{i=1}^L\in\T^D_S$,
the Virasoro operator $T\{\theta\}$ acting on $\M$ is defined by
\eqref{def:Ttheta2}, \eqref{def:Ttheta1}, \eqref{def:Tm} and
\eqref{def:sugawara}.  Using the dual bases above, we can represent the
Virasoro operator $T\{\theta\}$ in the following form:
\begin{equation}
  T\{\theta\}
  =
  - \frac{1}{2\kappa}
  \sum_{(a,b,m,i)\in R} \NP\rho(e^{a,b,i}_m\circ\theta)\rho(e_{a,b,i}^m)\NP,
  \label{Ttheta-e}
\end{equation}
where we put
\begin{equation*}
  e^{a,b,i}_m\circ\theta
  := (e^{a,b,i}_{m,j}(\xi_j)\theta_j(\xi_j))_{j=1}^L
  \in \g^D_S.
\end{equation*}
The formula \eqref{Ttheta-e} follows from the special cases with
$\theta=(\delta_{i,j}\xi_j^{n+1}\partial_{\xi_j})_{j=1}^L$ for $n\in\Z$, which
are obtained by straightforward calculations.

The bases which we really need later are however not these naively defined
bases, $F_0$ and $F^0$, but ``good'' dual frames in the sense of
\cite{tsuchimoto:93}. See \lemref{lem:JJ} and \remref{good-frames} below.
In order to construct such dual bases, we define the meromorphic functions
$w_{a,b}^n$ on $\UHP\times\C$ for $n=0,1,2,\ldots$ by derivatives of
$w_{a,b}$ \eqref{def:w-a,b}:
\begin{equation}
  w_{a,b}^n(t)
  = w_{a,b}^n(\tau;t)
  := \frac{(-1)^n}{n!} \pd{^n}{t^n} w_{a,b}(\tau;t).
  \label{def:wabn}
\end{equation}
Then the meromorphic function $w_{a,b}^n(\tau;t-z_i)$ on $\Xtilde$ can
be regarded as a global section of the line bundle $L_{a,b}(\ast Q_i)$
on $\X$ and its Laurent expansion in $\xi_i$ is written in the following form:
\begin{equation}
  w_{a,b}^n(\xi_i) = \xi_i^{-n-1} + (-1)^n w_{a,b,n} + O(\xi_i),
  \label{wabn-expansion}
\end{equation}
where $w_{a,b,n}$ is the coefficient of $t^n$ in the Laurent expansion
\eqref{w-a,b:expand} of $w_{a,b}(t)$.
For $(a,b)\in(\Z/N\Z)^2\setminus\{(0,0)\}$, $i,j=1,\dots,L$, and $m\in\Z$, we
define $f_{a,b,i,j}^m\in\O_S((\xi_j))$ by
\begin{equation}
  f_{a,b,i,j}^m(\xi_j) :=
  \begin{cases}
    \delta_{i,j} \xi_j^m          & \text{if $m\geqq0$}, \\
    w_{a,b}^{-m-1}(z_j-z_i+\xi_j) & \text{if $m<0$},
  \end{cases}
  \label{def:fabijm}
\end{equation}
and put
\begin{equation}
\begin{aligned}
  & J_{a,b,i}^m
  = (J_{a,b,i,j}^m)_{j=1}^L
  := (J_{a,b}\tensor f_{a,b,i,j}^m(\xi_j))_{j=1}^L,
  \\
  & J^{a,b,i}_m
  = (J^{a,b,i}_{m,j})_{j=1}^m
  := (J^{a,b}\tensor f_{-a,-b,i,j}^{-m-1}(\xi_j))_{j=1}^L.
\end{aligned}
\label{def:Jabim}
\end{equation}
Then the following topological $\O_S$-bases of $\g^D_S$
\begin{equation}
  F_1 := \{\, J_{a,b,i}^m \mid (a,b,m,i)\in R \,\},
  \quad
  F^1 := \{\, J^{a,b,i}_m \mid (a,b,m,i)\in R \,\},
  \label{dual-bases-J}
\end{equation}
are dual to each other by virtue of the residue theorem. Changing the
bases of expansion from $F_0$ and $F^0$ to $F_1$ and $F^1$, we obtain the
another expression of $T\{\theta\}$ from
\eqref{Ttheta-e}:
\begin{equation}
  T\{\theta\}
  =
  - \frac{1}{2\kappa}
  \sum_{(a,b,m,i)\in R} \NP\rho(J^{a,b,i}_m\circ\theta)\rho(J_{a,b,i}^m)\NP,
  \label{Ttheta-J}
\end{equation}
where we put, as above,
\begin{equation*}
  J^{a,b,i}_m\circ\theta
  := (J^{a,b,i}_{m,j}(\xi_j)\theta_j(\xi_j))_{j=1}^L
  \in \g^D_S.
\end{equation*}

The following lemma immediately follows from the definitions above and
the decomposition \eqref{decomp:gtwX} of $\gtw_\X$ to the direct sum of the
line bundles $L_{a,b}$.

\begin{lem}
\label{lem:JJ}
Under the notation above, we have the following:
\begin{enumerate}
\item If $m$ is negative, then $J_{a,b,i}^m$ is equal to the image of
  $J_{a,b}w_{a,b}^m(\tau;t-z_i)\in\gDPr$.
\item If $m$ is not negative and $\theta$ is a local section of 
  $\TDPr=\pr_\ast\T_{\X/S}(\ast D)$ with
  $\theta=\theta^t(\tau;z;t)\partial_t$, then $J^{a,b,i}_m\circ\theta$ is
  equal to the image of
  $J^{a,b}w_{-a,-b}^{-m-1}(\tau;t-z_i)\theta^t(\tau;z;t)\in\gDPr$.
\end{enumerate}
\end{lem}

\begin{rem}
\label{good-frames}
This lemma means that the topological dual bases $F_1$ and $F^1$ given by
\eqref{dual-bases-J} are {\em good dual frames} of $\g^D_S$ in the sense of
Tsuchimoto \cite{tsuchimoto:93}. 
\end{rem}

For $f_i\in\O_S((\xi_i))$, the regular part $f_{i,+}$ and the singular part
$f_{i,-}$ are uniquely defined by the conditions
\begin{equation*}
  f_i = f_{i,+} + f_{i,-},
  \quad
  f_{i,+}\in\O_S[[\xi_i]],
  \quad
  f_{i,-}\in\xi_i^{-1}\O_S[\xi_i^{-1}].
\end{equation*}
Note that the differentiation $\partial_{\xi_i}$ commutes the
operations $f_i\mapsto f_{i,\pm}$:
\begin{equation*}
  \partial_{\xi_i}(f_{i,\pm}(\xi_i)) = (f_i'(\xi_i))_\pm,
\end{equation*}
which shall be denoted by $f'_{i,\pm}(\xi_i)$.

\begin{lem}
\label{lem:NP}
For $f_i,g_i\in\O_S((\xi_i))$, the normal product 
$\NP(J^{a,b}\tensor f_i)(J_{a,b}\tensor g_i)\NP$ can be represented in the
following forms:
\begin{align*}
  \NP(J^{a,b}\tensor f_i)(J_{a,b}\tensor g_i)\NP
  &
  = (J^{a,b}\tensor f_i)(J_{a,b}\tensor g_i)
  + \khat \Res\limits_{\xi_i=0}(f_{i,+}(\xi_i)g'_{i,-}(\xi_i)\,d\xi_i)
  \\ &
  = (J_{a,b}\tensor g_i)(J^{a,b}\tensor f_i)
  - \khat \Res\limits_{\xi_i=0}(f_{i,-}(\xi_i)g'_{i,+}(\xi_i)\,d\xi_i).
\end{align*}
\end{lem}
\begin{pf}
 From the commutativity of $J_{a,b}$ and $J^{a,b}$ and the
definition of the normal product \eqref{def:normal-order}, we can find that
\begin{equation*}
  \NP(J^{a,b}\tensor f_i)(J_{a,b}\tensor g_i)\NP
  = (J^{a,b}\tensor f_i)(J_{a,b}\tensor g_{i,+})
  + (J_{a,b}\tensor g_{i,-})(J^{a,b}\tensor f_i).
\end{equation*}
Using the definition \eqref{triv:bracket-hatgSD} of the Lie algebra
structure on $\hat\g^D_S$, we obtain the formulae
\begin{align}
  & [J^{a,b}\tensor g_{i,-}, J_{a,b}\tensor f_i]
  = \khat \Res\limits_{\xi_i=0}
    \bigl( g'_{i,-}(\xi_i)f_i(\xi_i)\,d\xi_i \bigr)
  = \khat \Res\limits_{\xi_i=0}
    \bigl( f_{i,+}(\xi_i)g'_{i,-}(\xi_i)\,d\xi_i \bigr),
  \\
  & [J^{a,b}\tensor f_i, J_{a,b}\tensor g_{i,+}]
  = \khat \Res\limits_{\xi_i=0}
    \bigl( f'_i(\xi_i)g_{i,+}(\xi_i)\,d\xi_i \bigr)
  = - \khat \Res\limits_{\xi_i=0}
    \bigl( f_{i,-}(\xi_i)g'_{i,+}(\xi_i)\,d\xi_i \bigr).
\end{align}
These formulae prove the lemma.
\end{pf}

For the brevity of notation, we introduce the sets $R_\pm$ by
\begin{equation*}
  R_+:=\{\,(a,b,m,i)\in R\mid m \geqq 0 \,\},
  \quad
  R_-:=\{\,(a,b,m,i)\in R\mid m < 0 \,\},
\end{equation*}
where $R$ is defined by \eqref{def:R}.

\begin{lem}
\label{lem:Ttheta-J}
For $\theta=(\theta_i(\xi_i)\partial_{\xi_i})_{i=1}^L\in
\bigoplus_{i=1}^L\O_S((\xi_i))\partial_{\xi_i}$, we have the following
expressions of the Virasoro operator $T\{\theta\}$:
\begin{equation}
  T\{\theta\}
  = - \frac{1}{2\kappa}
  \Bigl(
  \sum_{R_+}
  \rho(J^{a,b,i}_m\circ\theta)\rho(J_{a,b,i}^m)
  +
  \sum_{R_-}
  \rho(J_{a,b,i}^m)\rho(J^{a,b,i}_m\circ\theta)
  \Bigr),
  \label{Ttheta-J1}
\end{equation}
where the symbols $\sum\limits_{R_\pm}$ denote the summations over
$(a,b,m,i)\in R_\pm$.  Moreover, substituting the definition \eqref{def:Jabim}
to this formula, we obtain
\begin{equation}
\begin{split}
  T\{\theta\}
  & = - \frac{1}{2\kappa}
  \Bigl(
  \sum_{R_+}
  \sum_{j=1}^L
  \rho_j(J^{a,b}\tensor w_{-a,-b}^m(z_j-z_i+\xi_j)\theta_j(\xi_j))
  \,
  \rho_i(J_{a,b}\tensor\xi_i^m)
  \\
  & \hphantom{= - \frac{1}{2\kappa}}\!\!
  +
  \sum_{R_-}
  \rho(J_{a,b,i}^m)
  \rho(J^{a,b,i}_m\circ\theta)
  \Bigr).
\end{split}
\label{Ttheta-J2}
\end{equation}
\end{lem}
\begin{pf}
 From the definitions \eqref{def:wabn}, \eqref{def:fabijm}, and
\eqref{def:Jabim}, we can derive the following formulae:
\begin{alignat*}{2}
  & \rho(J_{a,b,i}^m) 
  = \rho_i(J_{a,b}\tensor\xi_i^m)
  \quad & &
  \text{for $(a,b,m,i)\in R_+$},
  \\
  & \rho(J^{a,b,i}_m\circ\theta) 
  = \rho_i(J^{a,b}\tensor\xi_i^{-m-1}\theta_i(\xi_i))
  \quad & &
  \text{for $(a,b,m,i)\in R_-$},
  \\
  & J_{a,b,i,i}^m = J_{a,b}\tensor w_{a,b}^{-m-1}(\xi_i)
  \quad & & 
  \text{for $(a,b,m,i)\in R_-$},
  \\
  & \partial_{\xi_i}(w_{a,b,+}^{-m-1}(\xi)) = m w_{a,b,+}^{-m}(\xi_i)
  \quad & &
  \text{for $(a,b,m,i)\in R_-$}.
\end{alignat*}
Therefore, applying \lemref{lem:NP} to the formula \eqref{Ttheta-J}, we obtain
\begin{equation*}
  T\{\theta\}
  = (\text{the right-hand side of \eqref{Ttheta-J1}}) + \khat R(\theta),
\end{equation*}
where $R(\theta)$ is defined by
\begin{equation*}
  R(\theta)
  :=
  \sum_{R_-}
  \Res\limits_{\xi_i=0}
  \bigl(
    (\xi_i^{-m-1}\theta_i(\xi_i))_-
    \,
    (-m)w_{a,b,+}^{-m}(\xi_i)
    \,d\xi_i
  \bigr),
\end{equation*}
which is a finite sum and hence is a local section of $\O_S$.  Hence, for the
proof of \eqref{Ttheta-J1}, it is enough to show that $R(\theta)=0$ in the
case of $\theta=(\delta_{i,j}\xi_j^{-l}\partial_{\xi_j})_{j=1}^L$ for
$l=1,2,3,\ldots$ and $i=1,\dots,L$.  Using the formula \eqref{def:wabn} and
the Laurent expansion \eqref{w-a,b:expand} of $w_{a,b}(t)$, we can find that
\begin{align*}
  & \sum_{m<0}
  \Res\limits_{\xi_i=0}
  \bigl(
    (\xi_i^{-m-1}\xi_i^{-l})_-
    \,
    (-m)w_{a,b,+}^{-m}(\xi_i)
    \,d\xi_i
  \bigr)
  \\
  &=
  \sum_{n=1}^l \Res\limits_{\xi_i=0}
  \left(
    \xi_i^{n-l-1}
    \,
    n (-1)^n \sum_{\nu\geqq n}
    {\nu \choose n} w_{a,b,\nu}\xi_i^{\nu-n}
    \,d\xi_i
  \right)
  \\
  &
  =
  \left(
    \sum_{n=1}^l (-1)^n n {l \choose n}
  \right)
  w_{a,b,l}
  =
  \begin{cases}
    - w_{a,b,1} & \text{if $l=1$}, \\
    0           & \text{if $l>2$}.
  \end{cases}
\end{align*}
Hence we obtain $R(\theta)=0$ by \eqref{w-a,b,1:sum}.
We have proved the first expression \eqref{Ttheta-J1} of $T\{\theta\}$.
\end{pf}

\begin{rem}
The second expression \eqref{Ttheta-J2} of $T\{\theta\}$ is very useful for
the explicit calculations of the elliptic Knizhnik-Zamolodchikov connections
(\secref{elliptic-kz}).
\end{rem}

We are ready to prove \lemref{lem:TDPrSubsetgDPr}.  Let $v$ be in $\M$ and
$\theta$ in $\TDPr$.  Then, applying \lemref{lem:JJ} to the first expression
\eqref{Ttheta-J1} of $T\{\theta\}$ in \lemref{lem:Ttheta-J} shows that
$T\{\theta\}v \in \gDPr\M$.  Hence we have proved
\lemref{lem:TDPrSubsetgDPr}.


\section{Flat connections}
\label{flat-connections}

We keep all the notation in the previous section.  In this section, we assume
$\kappa=k+\hvee\ne0$ and shall construct $\D_S$-module structures on the sheaf
$\sCC(\M)$ of conformal coinvariants and on the sheaf $\sCB(\M)$ of conformal
blocks. We can show as a direct consequence that, under the assumption of
\lemref{sheaf-coherent}, the sheaves $\sCC(\M)$ and $\sCB(\M)$ are locally
free coherent $\O_S$-modules (i.e., vector bundles on $S$) with flat
connections and their fibers at $s\in S$ are canonically isomorphic to the
space of conformal coinvariants and that of conformal blocks respectively.  In
\secref{elliptic-kz}, we shall show that the connections on $\sCC_k(V)$ and
$\sCB_k(V)$ coincide with the elliptic Knizhnik-Zamolodchikov equations
introduced by Etingof \cite{eti:94}.  In \secref{automorphic}, we shall obtain
a proof of the modular property of the connections without referring the
explicit expressions of them.

\subsection{Construction of flat connections}
\label{const-flat-conn}

Recall that we have the Lie algebra extension $\VirDPr$ of the tangent sheaf
$\T_S$ by $\TDPr$ in \eqref{eq:LieAlgExtOfTSByTDPr}.  Since the action of
$\TDPr$ maps $\M$ into $\gDPr\M$ due to \lemref{lem:TDPrSubsetgDPr},
the representation of the Lie algebra $\VirDPr$ on $\sCC(\M)$ given by
\lemref{lem:Vir-ActsOnsCC} induces the Lie algebra action of $\T_S$ on
$\sCC(\M)$, which shall be denoted by
\begin{equation}
  \T_S\times\sCC(\M)\to\sCC(\M),
  \quad
  (\mu,\phi) \mapsto D_\mu\phi.
\label{def:conn-cc}
\end{equation}
Moreover it immediately follows from \eqref{def:VirActionOnM} that
\begin{equation*}
  D_{f\mu}v = f(D_\mu v),
  \quad
  D_\mu(fv) = \mu(f)v + f(D_\mu v)
\end{equation*}
for $f\in\O_S$, $\mu\in\T_S$, and $v\in\sCC(\M)$.  Thus we obtain the flat
connection $D$ on the sheaf $\sCC(\M)$ of conformal coinvariants.  Because of
$\sCB(\M)=\HOM_{\O_S}(\sCC(\M),\O_S)$, we obtain the dual connection $D^\ast$
on $\sCB(\M)$:
\begin{equation}
    (D_\mu^\ast \Phi)(v) := \mu (\Phi(v)) - \Phi(D_\mu v).
\label{def:conn-cb}
\end{equation}

We can summarize the results as follows.

\begin{thm}
\label{thm:D-module-str}
For each $i=1,\ldots,L$, let $M_i$ be a representation with level $k$ of the
affine Lie algebra $\ghat_i$ satisfying the smoothness condition
\eqref{smoothness}.  Put $M:=\bigotimes_{i=1}^L M_i$ and $\M:=M\tensor\O_S$.
Then $M$ is a representation with level $k$ of $(\g^{\oplus L})^\wedge$ and
$\M$ is a $\hat\g^D_S$-module.  Assume that $\kappa=k+\hvee\ne0$. Then the
action of $\VirDPr$ on $\M$ induces the $\D_S$-module structures on the sheaf
$\sCC(\M)$ of conformal coinvariants and on the sheaf $\sCB(\M)$ of conformal
blocks.
\end{thm}

\begin{rem}
\label{rem:desc-Dast}
We have another description of the dual connection $D^\ast$ on $\sCB(\M)$.
\lemref{lem:TDPrKillsCB} shows that the representation of the Lie algebra
$\VirDPr$ on $\sCB(\M)$ given by \lemref{lem:Vir-ActsOnsCB} induces the Lie
algebra action of $\T_S$ on $\sCB(\M)$, which coincides with the dual
connection $D^\ast$.  
\end{rem}

\begin{cor}
\label{locally-freeness}
For each $i=1,\ldots,L$, let $V_i$ be a finite-dimensional irreducible
representation of $\g$ and $M_i$ a quotient module of the Weyl module
$\Weyl_k(V_i)$ over the affine Lie algebra $\hat\g_i$.  Put
$M:=\bigotimes_{i=1}^L M_i$ and $\M:=M\otimes\O_S$.  Assume that
$\kappa=k+\hvee\ne0$.  Then the sheaf $\sCC(\M)$ of conformal coinvariants and
the sheaf $\sCB(\M)$ of conformal blocks are locally free coherent
$\O_S$-modules with flat connections on $S$ and dual to each other.
\end{cor}
\begin{pf}
We have already shown the $\O_S$-coherencies of $\sCC(\M)$ and $\sCB(\M)$ in
\corref{sheaf-coherent} and the existences of $\D_S$-module structures on
$\sCC(\M)$ and on $\sCB(\M)$ in \thmref{thm:D-module-str}.  It is
well-known that any $\O_S$-coherent $\D_S$-module 
is $\O_S$-locally free.
(See Theorem 6.1 in Chapter I of \cite{hotta},
Proposition 1.7 in Chapter VI of \cite{borel-etal:87},
or Theorem 1.1.25 of \cite{bjork:93}.)
\end{pf}

\begin{cor}
\label{cor:ident-fiber}
For each $i=1,\ldots,L$, let $M_i$ be a representation with level $k$
of the affine Lie algebra $\ghat_i$ satisfying the smoothness
condition \eqref{smoothness}.  Put $M:=\bigotimes_{i=1}^L M_i$ and
$\M:=M\tensor\O_S$.  Then the fiber of $\sCC(\M)$ at $s=(\tau;z)\in S$
is canonically isomorphic to the space of conformal coinvariants for
$(X;\{Q_i\}_{i=1}^L)=(X_\tau;\{q_i(s)\}_{i=1}^L)$:
\begin{equation}
  \sCC(\M)|_s \simeqq \CC(M).
  \label{id-fib-CC}
\end{equation}
Moreover, under the assumptions of \corref{locally-freeness}, the
fiber of $\sCB(\M)$ at $s\in S$ is canonically isomorphic to the space
of conformal blocks:
\begin{equation}
  \sCB(\M)|_s \simeqq \CB(M).
  \label{id-fib-CB}
\end{equation}
\end{cor}
\begin{pf}
We can identify the restriction of $\gtw_\X$ on $X=\pr^{-1}(s)$ with $\gtw$
in \secref{conf-block}. Put $Q:=D\cap X$.  Then we can find the following
canonical isomorphisms without using \corref{locally-freeness}:
\begin{align}
  &
  \M|_s 
  = (M\tensor\O_S)|_s
  \simeqq M\tensor\C
  = M,
  \label{id-fib-M}
  \\ &
  \left. (\pr_\ast\gtw_\X(m D)) \right|_s
  \simeqq H^0(X,\gtw(mQ))
  \quad\text{for $m\geqq0$},
  \label{id-fib-H0gtw}
  \\ &
  \left. \gDPr \right|_s
  = \left. (\pr_\ast\gtw_\X(\ast D)) \right|_s
  \simeqq H^0(X,\gtw(\ast Q))
  = \g^Q_\prin,
  \label{id-fib-gDpr}
  \\ &
  \left. (\gDPr\M) \right|_s
  \simeqq \g^Q_\prin M.
  \label{id-fib-gDprM}
\end{align}
The Riemann-Roch theorem shows that 
$\dim_\Comp H^0(\pr^{-1}(s),\gtw(m D)|_{\pr^{-1}(s)})$ is a constant
function of
$s\in S$ if $m\geqq0$. Therefore the existence of the isomorphism
\eqref{id-fib-H0gtw} follows from the Grauert theorem.
(For the Grauert theorem, see, for example, Corollary 12.9 in Chapter III of
\cite{hartshorne:77} and Theorem in Chapter 10 \S5.3 (p.~209) of
\cite{gra-rem:84}).
The isomorphism \eqref{id-fib-gDpr} is obtained by the inductive limit of
\eqref{id-fib-H0gtw}.  The isomorphism \eqref{id-fib-gDprM} is obtained by
using \eqref{id-fib-M}, \eqref{id-fib-gDpr}, and applying the right exact
functor $(\cdot)|_s$ to the exact sequence
\(
    \gDPr \tensor_{\O_S} \M \to \gDPr \M \to 0
\). %
Similarly the isomorphism \eqref{id-fib-CC} is obtained from
\eqref{id-fib-M} and \eqref{id-fib-gDprM}.

Under the assumption of \corref{locally-freeness}, the sheaves $\sCC(\M)$ and
$\sCB(\M)$ are locally free $\O_S$-modules of finite rank and dual to each
other.  Then we have $\sCB(\M)|_s = \Hom_\C(\sCC(\M)|_s,\C)$.
This formula together with the isomorphism \eqref{id-fib-CC} gives the
isomorphism \eqref{id-fib-CB}.
\end{pf}

Let us describe the connections $D$ and $D^\ast$ more explicitly.  For any
vector field $\mu\in\T_S$, the action of $D_\mu$ on $\sCC(\M)$ and that of
$D^\ast_\mu$ on $\sCB(\M)$ are described as follows.  Using the short exact
sequence \eqref{eq:LieAlgExtOfTSByTDPr}, we can lift $\mu$ to an element
$(\mu;\theta;0)\in\VirDPr$ at least locally on $S$, and the ambiguity in the
choice of $(\mu;\theta;0)$ is equal to $\TDPr$.  But, since the actions of
$\TDPr$ on $\sCC(\M)$ and $\sCB(\M)$ are trivial, the actions of
$(\mu;\theta;0)$ on $\sCC(\M)$ and $\sCB(\M)$ do not depend on the choice of
the lift and give $D_\mu$ and $D^\ast_\mu$:
\begin{align}
  & D_\mu v 
  = (\mu;\theta;0)\cdot v = \mu(v) + T\{\theta\}v,
  \label{def:conn-cc:explicit}
  \\
  & D^\ast_\mu \Phi 
  = (\mu;\theta;0)\cdot\Phi = \mu(\Phi) + T^\ast\{\theta\}\Phi
\label{def:conn-cb:explicit}
\end{align}
for $(\mu;\theta;0)\in\VirDPr$ and $v\in\sCC(\M)$, and $\Phi\in\sCB(\M)$.
The following lemma provides us with the explicit formulae of
$(\mu;\theta;0)\in\VirDPr$.

\begin{lem}
\label{explicit-formula-of-conn}
The following expressions in the coordinate $(\tau;z;t)$ define vector
fields in $\VirDPr=\pr_\ast\T_{\X,\pr}(\ast D)$ and are lifts of
$\partial_{z_i}$ and $\partial_\tau$ respectively:
\begin{equation}
  a_\Prin(\partial_{z_i}) = \partial_{z_i}, \quad
  a_\Prin(\partial_\tau) = \partial_\tau + Z(\tau;z;t)\partial_t,
  \label{aprin-1}
\end{equation}
where the function $Z(\tau;z;t)$ is a global meromorphic function on $\Xtilde$
satisfying the following properties:
\begin{enumerate}
\item The poles of $Z(\tau;z;t)$ are contained in $\pr_{\Xtilde/\X}^{-1}(D)$;
\item The quasi-periodicity:
  \begin{equation}
    Z(\tau;z;t+m\tau+n) = Z(\tau;z;t) + m
    \quad\text{for $(m,n)\in\Z^2$}.
    \label{property-Z}
  \end{equation}
\end{enumerate}
\end{lem}
\begin{pf}
Since the action of $(m,n)\in\Z^2$ on $\Xtilde$ sends $\partial_\tau$,
$\partial_{z_i}$, and $\partial_t$ to $\partial_\tau+n\partial_t$,
$\partial_{z_i}$, and $\partial_t$ respectively, the expressions
\eqref{aprin-1} are vector fields in $\pr_\ast\T_{\X,\pr}(\ast D)$.
\end{pf}

\begin{example}
\label{ex:Z}
We can use the following function as $Z(\tau;z;t)$ for any
$i_0\in\{1,\ldots,L\}$. (cf.\ \S3.1 of \cite{fel-wie:96}):
\begin{equation}
  Z(\tau;z;t) = Z_{1,1}(\tau;t-z_{i_0}),
  \quad
  Z_{1,1}(\tau;t)
  = 
  -\frac{1}{2\pi i}
  \frac{\theta'_{[0,0]}(t)}{\theta_{[0,0]}(t)}.
  \label{def:Z}
\end{equation}
(See \appref{theta-func} for the notation of the theta functions.)
\end{example}

\begin{lem}
\label{lem:explicit-formula-D-and-Dast}
The connections $D$ on $\sCC(\M)$ and $D^\ast$ on $\sCB(\M)$ possess the
following expressions:
\begin{align}
  & D_{\PD{}{z_i}}
  = \partial_{z_i} - \rho_i(T\{\partial_{\xi_i}\}),
  \label{explicit-Dzi}
  \\
  & D_{\PD{}{\tau}} 
  = \partial_\tau + T\{Z(\tau;z;t)\partial_t\}
  = \partial_\tau 
  + \sum_{i=1}^L\rho_i(T\{Z(\tau;z;z_i+\xi_i)\partial_{\xi_i}\}),
  \label{explicit-Dtau}
  \\
  & D^\ast_{\PD{}{z_i}}
  = \partial_{z_i} - \rho^\ast_i(T\{\partial_{\xi_i}\}),
  \label{explicit-Dastzi}
  \\
  & D^\ast_{\PD{}{\tau}} 
  = \partial_\tau + T^\ast\{Z(\tau;z;t)\partial_t\}
  = \partial_\tau 
  + \sum_{i=1}^L\rho^\ast_i(T\{Z(\tau;z;z_i+\xi_i)\partial_{\xi_i}\}).
  \label{explicit-Dasttau}
\end{align}
\end{lem}
\begin{pf}
In the local coordinate $(\tau;z;\xi_j)$ with $\xi_j=t-z_j$, the vector
fields given by \eqref{aprin-1} can be represented in the following forms
around $Q_i$:
\begin{equation}
  a_\Prin(\partial_{z_i})
  = \partial_{z_i} - \delta_{i,j}\partial_{\xi_j},
  \quad
  a_\Prin(\partial_\tau)
  = \partial_\tau + Z(\tau;z;z_j+\xi_j)\partial_{\xi_j}.
  \label{aprin-2}
\end{equation}
Namely, their images in $\Vir^D_S$ are of the following forms:
\begin{equation}
  a_\Prin(\partial_{z_i})
  = (\partial_{z_i} ; (-\delta_{i,j}\partial_{\xi_j})_{j=1}^L ; 0),
  \quad
  a_\Prin(\partial_\tau)
  = (\partial_\tau; (Z(\tau;z;z_j+\xi_j)\partial_{\xi_j})_{j=1}^L ; 0),
  \label{aprin-3}
\end{equation}
whose actions on $\sCC(\M)$ and $\sCB(\M)$ can be written in the forms
\eqref{explicit-Dzi}, \eqref{explicit-Dtau}, \eqref{explicit-Dastzi}, and
\eqref{explicit-Dasttau}.
\end{pf}

\begin{rem}
The explicit formulae \eqref{explicit-Dastzi} and \eqref{explicit-Dasttau} of
the flat connection coincide with the expressions of $\nabla_\tau$ and
$\nabla_{z_i}$ in \S3.1 of \cite{fel-wie:96}. In order to prove the flatness
of their connection, Felder and Wieczerkowski use the explicit expression of
$\nabla_{z_i}$ which corresponds to \eqref{kz-operator} in our case.  However
in our construction the flatness of the connection is a priori obvious.
\lemref{Vir-toVirIsLieAlgHom} is the key lemma for the proof the flatness.
Our proof of the flatness can be also applied to Proposition 3.4 of
\cite{fel-wie:96}.
\end{rem}

\subsection{Elliptic Knizhnik-Zamolodchikov equations}
\label{elliptic-kz}

In this subsection, we show that the flat connection on the sheaf of conformal
blocks defined in \secref{const-flat-conn} is nothing but the elliptic
Knizhnik-Zamolodchikov equation introduced by Etingof \cite{eti:94}, when
$\M = \sWeyl_k(V)$. 

Let $V_i$ be a finite-dimensional irreducible representation of $\g=sl_N(\C)$
for each $i$ and put $V:=\bigotimes_{i=1}^L V_i$.  In view of
\propref{sheaf-restriction}, the connections $D$ on $\sCC_k(V)$ and $D^\ast$ on
$\sCB_k(V)$ are regarded as connections on $V\tensor\O_S$ and on
$V^\ast\tensor\O_S$ respectively and are dual to each other.  In order to give
the explicit expressions of these connections, we define the function
$Z_{a,b}(t)=Z_{a,b}(\tau;t)$ by
\begin{equation}
  Z_{a,b}(t) :=
  \frac{w_{a,b}(t)}{4\pi i}
  \left(
    \frac{\theta'_{[a,b]}(t)}{\theta_{[a,b]}(t)}
  - \frac{\theta'_{[a,b]}   }{\theta_{[a,b]}   }
  \right).
\label{def:Z-a,b}
\end{equation}
The point $t=0$ is an apparent singularity of this function and we
can analytically continue $Z_{a,b}(t)$ to $t=0$ by
\begin{equation}
  Z_{a,b}(0)
  := 
  \frac{1}{4\pi i}
  \left(
    \frac{\theta''_{[a,b]}}{\theta_{[a,b]}}
    -
    \left( \frac{\theta'_{[a,b]}}{\theta_{[a,b]}} \right)^2
  \right).
  \label{Zab(0)}
\end{equation}

\begin{thm}
\label{kz}
As operators acting on $V\tensor\O_S$ and $V^\ast\tensor \O_S$, the operators
$D_\mu$ and $D^\ast_\mu$ for $\mu=\PD{}{z_i},\PD{}{\tau}$ possess the
following expressions:
\begin{align}
  D_{\PD{}{z_i}}
  & = \frac{\partial}{\partial z_i}
  - \frac{1}{\kappa}
  \sum_{j\neq i}\sum_{(a,b)\neq(0,0)}
    w_{a,b}(z_j - z_i) \rho_j(J_{a,b}) \rho_i(J^{a,b}),
  \label{KZ-operator}
  \\
  D_{\PD{}{\tau}}
  & = \frac{\partial}{\partial \tau}
  + \frac{1}{\kappa}
  \sum_{i,j = 1}^L \sum_{(a,b)\neq(0,0)}
    Z_{a,b}(z_j - z_i) \rho_j(J_{a,b}) \rho_i(J^{a,b}),
  \label{KZB-operator}
  \\
  D^\ast_{\PD{}{z_i}}
  & = \frac{\partial}{\partial z_i}
  + \frac{1}{\kappa}
  \sum_{j\neq i}\sum_{(a,b)\neq(0,0)}
    w_{a,b}(z_j - z_i) \rho^\ast_j(J_{a,b}) \rho^\ast_i(J^{a,b}),
  \label{kz-operator}
  \\
  D^\ast_{\PD{}{\tau}}
  & = \frac{\partial}{\partial \tau}
  - \frac{1}{\kappa}
  \sum_{i,j = 1}^L \sum_{(a,b)\neq(0,0)}
    Z_{a,b}(z_j - z_i) \rho^\ast_j(J_{a,b}) \rho^\ast_i(J^{a,b}),
  \label{kzb-operator}
\end{align}
where $\rho_i(A_i)$ and $\rho^\ast_i(A_i)$ for $A_i\in\g$ act on the $i$-th
factor $V_i$ of $V$ and on the $i$-th factor $V^\ast_i$ of $V^\ast$
respectively. Here, for each $(a,b)\neq(0,0)$, the function
$w_{a,b}(t)=w_{a,b}(\tau;t)$ is defined by \eqref{def:w-a,b} and the function
$Z_{a,b}(t)=Z_{a,b}(\tau;t)$ by \eqref{def:Z-a,b} and \eqref{Zab(0)}.
\end{thm}
\begin{pf}
Since the connections $D$ and $D^\ast$ are dual to each other, for the proof
of the proposition it suffices to obtain either the formulae for $D_\mu$
(i.e., \eqref{KZ-operator} and \eqref{KZB-operator}) or those for $D^\ast_\mu$
(i.e., \eqref{kz-operator} and \eqref{kzb-operator}).  The formulae
\eqref{kz-operator} and \eqref{kzb-operator} can be proved in the same way as
\lemref{gaudin-as-sug-correl} or the statements in \S6 of \cite{f-f-r:94}.
But we shall give the proof of \eqref{KZ-operator} and \eqref{KZB-operator}
using the expression \eqref{Ttheta-J2} of $\T\{\theta\}$ in
\lemref{lem:Ttheta-J}.

Let us fix $v\in V^\ast\tensor\O_S\subset\M:=\sWeyl_k(V)$. We rewrite
$\rho_i(T\{\partial_{\xi_i}\})v$ and $T\{Z(\tau;z;t)\partial_t\}v$ in
\lemref{lem:explicit-formula-D-and-Dast} modulo $\gDPr\M$ in terms of
operators acting on $V$, 
applying the Ward identity \eqref{sheaf-ward-cc} in the following form:
\begin{align}
  & \rho_i(J^{a,b}\tensor w_{-a,-b}^n(\xi_i)) v'
  \equiv
  - \sum_{j\ne i} w_{-a,-b}^n(z_j-z_i) \rho_j(J^{a,b})v',
  \label{ward-id-1}
  \\
  & \rho_i(J^{a,b}\tensor \xi^{-n-1})v'
  \equiv
  (-1)^{n+1} w_{-a,-b,n}\rho_i(J^{a,b})v'
  - \sum_{j\ne i} w_{-a,-b}^n(z_j-z_i) \rho_j(J^{a,b})v',
  \label{ward-id-2}
\end{align}
for $v'\in V\tensor\O_S$ and $n\geqq0$. Here we used the Laurent
expansion \eqref{wabn-expansion} and \lemref{lem:JJ}.

First we prove \eqref{KZ-operator} from \eqref{explicit-Dzi}.
The formula \eqref{Ttheta-J2} for 
$\theta=(\delta_{i,j}\partial_{\xi_j})_{j=1}^L$ together with \lemref{lem:JJ}
shows
\begin{equation}
  \rho_i(T\{\partial_{\xi_i}\})v
  \equiv
  - \frac{1}{2\kappa}
  \sum_{\begin{Sb} (a,b)\ne(0,0) \\ i'=1,\dots,L \end{Sb}}
    \rho_i(J^{a,b}\tensor w_{-a,-b}(z_i-z_{i'}+\xi_i))
    \rho_{i'}(J_{a,b})
    v.
  \label{Tpartialxi-1}
\end{equation}
Applying the Ward identity \eqref{ward-id-1} ($n=0$) to the terms with
$i'=i$ in the right-hand side of \eqref{Tpartialxi-1}, we find that
\begin{equation}
\begin{aligned}
  \rho_i(T\{\partial_{\xi_i}\})v
  & \equiv
  \frac{1}{2\kappa}
  \sum_{(a,b)\ne(0,0)}
  \sum_{j'\ne i}
  w_{-a,-b}(z_{j'}-z_i) \rho_{j'}(J^{a,b}) \rho_i(J_{a,b})v
  \\
  &
  -
  \frac{1}{2\kappa}
  \sum_{(a,b)\ne(0,0)}
  \sum_{i'\ne i}
  w_{-a,-b}(z_i-z_{i'}) \rho_i(J^{a,b}) \rho_{i'}(J_{a,b})v.
\end{aligned}
\label{Tpartialxi-2}
\end{equation}
Renumbering the indices of the first sum by $(a,b)\mapsto(-a,-b)$ and applying
\eqref{w-a,b:property} to the second, we conclude that
\begin{equation*}
  \rho_i(T\{\partial_{\xi_i}\})v
  \equiv
  \frac{1}{\kappa}
  \sum_{(a,b)\ne(0,0)}
  \sum_{j\ne i}
  w_{a,b}(z_j-z_i) \rho_j(J_{a,b}) \rho_i(J^{a,b})v.
\end{equation*}
This means that the operator $D_{\PD{}{z_i}}$ acting on $V\tensor\O_S$ is
represented as \eqref{KZ-operator}.

The expression \eqref{KZB-operator} of $D_{\PD{}{\tau}}$ can be deduced
from \eqref{explicit-Dtau} in the similar manner, which is however more
involved.  Let us use here the function $Z(t)=Z(\tau;z;t)$ given by
\exref{ex:Z}.  Then $Z(t)$ is regular at $t\ne z_{i_0}$ and the
Laurent expansion of $Z(t)$ in $\xi_{i_0}=t-z_{i_0}$ is represented as
\begin{equation}
  Z(z_{i_0}+\xi_{i_0})
  = - \frac{1}{2\pi i}(\xi_{i_0}^{-1} + Z_1 \xi_{i_0} + O(\xi_{i_0}^3)),
  \quad
  Z_1 := \frac{\theta'''_{[0,0]}}{3\theta'_{[0,0]}}.
  \label{Z-expansion}
\end{equation}
The formula \eqref{Ttheta-J2} and \lemref{lem:JJ} imply
\begin{equation}
  T\{Z(\tau;z;t)\partial_t\}v
  \equiv
  - \frac{1}{2\kappa}
  \sum_{(a,b)\ne(0,0)}
  \sum_{i,j=1}^L
  v_{a,b,i,j},
  \label{TZ-1st}
\end{equation}
where the vectors $v_{a,b,i,j}$ are defined by
\begin{equation*}
  v_{a,b,i,j}
  :=
  \rho_j(J^{a,b}\tensor w_{-a,-b}(z_j-z_i+\xi_j)Z(z_j+\xi_j))
  \rho_i(J_{a,b}) v.
\end{equation*}
Using the Laurent expansions \eqref{w-a,b:expand}, \eqref{Z-expansion}, and
the Ward identities \eqref{ward-id-1} ($n=0$) and \eqref{ward-id-2}
($n=0,1$), we can find the following expressions for the vectors
$v_{a,b,i,j}$:

\noindent $\bullet$ \enspace %
If $j\ne i$ and $j\ne i_0$, then
\begin{equation*}
  v_{a,b,i,j}
  =
  w_{-a,-b}(z_j-z_i)Z(z_j)
  \rho_j(J^{a,b})\rho_i(J_{a,b})v.
\end{equation*}

\noindent $\bullet$ \enspace %
If $j=i$ and $j\ne i_0$, then
\begin{equation*}
  v_{a,b,i,i}
  \equiv
  Z'(z_i)
  \rho_i(J^{a,b})\rho_i(J_{a,b})v
  - \sum_{j'\ne i}
  w_{-a,-b}(z_{j'}-z_i)Z(z_i)
  \rho_{j'}(J^{a,b})\rho_i(J_{a,b})v.
\end{equation*}

\noindent $\bullet$ \enspace %
If $j\ne i$ and $j=i_0$, then
\begin{align*}
  v_{a,b,i,i_0}
  & \equiv
  - \frac{1}{2\pi i}
  (- w_{-a,-b,0}w_{-a,-b}(z_{i_0}-z_i) + w'_{-a,-b}(z_{i_0}-z_i))
  \rho_{i_0}(J^{a,b})\rho_i(J_{a,b})v
  \\
  &
  + \frac{1}{2\pi i}
  \sum_{j'\ne i_0}
  w_{-a,-b}(z_{i_0}-z_i)w_{-a,-b}(z_{j'}-z_{i_0})
  \rho_{j'}(J^{a,b})\rho_i(J_{a,b})v.
\end{align*}

\noindent $\bullet$ \enspace %
If $j=i$ and $j=i_0$, then
\begin{align*}
  v_{a,b,i_0,i_0}
  & \equiv 
  - \frac{1}{2\pi i}
  (2w_{-a,-b,1} - (w_{-a,-b,0})^2  + Z_1)
  \rho_{i_0}(J^{a,b})\rho_{i_0}(J_{a,b})v
  \\
  &
  + \frac{1}{2\pi i}
  \sum_{j'\ne i_0}
  (- w'_{-a,-b}(z_{j'}-z_{i_0}) + w_{-a,-b,0}w_{-a,-b}(z_{j'}-z_{i_0}))
  \rho_{j'}(J^{a,b})\rho_{i_0}(J_{a,b})v.
\end{align*}
Note that $w_{-a,-b}^1(t)= - w'_{-a,-b}(t)$.  Substituting these expressions
into \eqref{TZ-1st}, we obtain an expression for $T\{Z(\tau;z;t)\partial_t\}v$
like
\begin{equation}
  T\{Z(\tau;z;t)\partial_t\}v
  = - \frac{1}{2\kappa}
  \sum_{\begin{Sb} (a,b)\ne(0,0) \\ i,j=1,\dots,L \end{Sb}}
  Z_{a,b,i,j} \rho_j(J_{a,b})\rho_i(J^{a,b})v
  \quad
  \text{with $Z_{a,b,j,i}=Z_{-a,-b,i,j}$}.
  \label{TZ-2nd}
\end{equation}
Note that
$\rho_j(J^{a,b})\rho_i(J_{a,b})=\rho_j(J_{-a,-b})\rho_i(J^{-a,-b})$.  It
is possible to compute all $Z_{a,b,i,j}$ directly, but we take a short
cut. Since the coefficients $Z_{a,b,i,j}$ do not depend on the
representations $V_i$ and the choice of $Z(t)=Z(\tau;z;t)$, for the
determination of all $Z_{a,b,i,j}$, we have only to calculate
$Z_{a,b,i_0,j}$ for $j=1,\dots,L$.  Picking up the terms which should be
contained in $Z_{a,b,i_0,j}$ from the expressions for $v_{a,b,i,j}$, we
obtain
\begin{align}
  & Z_{a,b,i_0,i_0}
  = - \frac{1}{2\pi i}(2w_{a,b,1} - (w_{a,b,0})^2  + Z_1),
  \label{Zab00-1st}
  \\
  & Z_{a,b,i_0,j}
  = w_{a,b}(z_j-z_{i_0})Z(z_j) 
  - \frac{1}{2\pi i}(- w_{a,b,0}w_{a,b}(z_j-z_{i_0}) + w'_{a,b}(z_j-z_{i_0})),
  \label{Zab0j-1st}
\end{align}
where we assume $j\ne i_0$.  (Here we have used the formulae
\eqref{w-a,b:property} and $w_{-a,-b,\nu}=(-1)^{\nu+1}w_{a,b,\nu}$.)  Note
that we have 
\begin{equation}
  w'_{a,b}(t)
  = w_{a,b}(t) (\log w_{a,b}(t))'
  = w_{a,b}(t) 
    \left(
      \frac{\theta'_{[a,b]}(t)}{\theta_{[a,b]}(t)}
    - \frac{\theta'_{[0,0]}(t)}{\theta_{[0,0]}(t)}
    \right).
  \label{w-a,b:deriv}
\end{equation}
Substituting \eqref{def:w-a,b}, \eqref{w-a,b:coeff}, \eqref{w-a,b:deriv}, and
the expression of $Z_1$ in \eqref{Z-expansion} to \eqref{Zab00-1st} and
\eqref{Zab0j-1st}, we can find the following results:
\begin{align}
  & Z_{a,b,i_0,i_0} = -2 Z_{a,b}(0) = -2 Z_{a,b}(z_{i_0}-z_{i_0}),
  \\
  & Z_{a,b,i_0,j} = -2 Z_{a,b}(z_j-z_{i_0}).
\end{align}
These formulae complete the proof of \eqref{KZB-operator}.
\end{pf}

\begin{rem}
Etingof found in \cite{eti:94} that certain twisted traces $F$ of vertex
operators of $\widehat{sl}_N$ satisfy $D^\ast_{\PD{}{z_i}}F = 0$ and 
$D^\ast_{\PD{}{\tau}}F = 0$, namely, that they are flat sections of the
connection $D^\ast$.
\end{rem}

\subsection{Modular invariance of the flat connections}
\label{automorphic}

In \cite{eti:94}, Etingof proved the modular invariance of the elliptic
Knizhnik-Zamolodchikov equations by explicit computation.  In this
subsection, we give a more geometric proof of this fact without use of the
explicit formulae of the connections.  A similar property was proved for the
non-twisted WZW model on elliptic curves in \cite{fel-wie:96}.

For the brevity of notation, we introduce the following symbols:
\begin{itemize}
\item Put $\Gamma:=SL_2(\Z)$;
\item In the following, we assume that the symbol $\gamma$ always denote an
  arbitrary matrix $\begin{pmatrix}a&b\\c&d\end{pmatrix}\in\Gamma$;
\item For $s=(\tau;z)\in S$, $t\in\C$, and $i=1,\dots,L$, put
  \begin{align*}
    &
    \check\tau := \frac{a\tau+b}{c\tau+d},
    \quad
    \check z = (\check z_i)_{i=1}^L 
    := \frac{z}{c\tau+d} = \left(\frac{z_i}{c\tau+d}\right)_{i=1}^L,
    \\ &
    \check t := \frac{t}{c\tau+d},
    \quad
    \check\xi_i := \frac{\xi_i}{c\tau+d},
    \quad
    \check s := (\check\tau; \check z).
  \end{align*}
\end{itemize}
As is well-known, $\gamma\in\Gamma$ acts on $S$ by
\begin{equation}
  \gamma\cdot s := \check s = (\check\tau;\check z)
  \quad\text{for $s=(\tau;z)\in S$.}
  \label{def:GammaOnS}
\end{equation}

In order to define actions of $\Gamma$ on $\X$, $\Gtw_\X$ and $\gtw_\X$,
we first extend the actions of $\Z^2$ on $\Xtilde$, $G$, and $\g$ defined
by \eqref{def:action-Z2-on-Xtilde}, \eqref{def:Gal-action-on-G}, and
\eqref{def:Gal-action-on-g} respectively to those of
the semi-direct product group $\Gammatilde:=\Z^2\rtimes\Gamma$ defined by
\begin{equation*}
  (m,n;\gamma)(m',n';\gamma')
  :=(m+m'd-n'c, n-m'b+n'a;\gamma\gamma')
\end{equation*}
for $(m,n;\gamma),(m',n';\gamma')\in\Z^2\times\Gamma$.  Then we have
\begin{equation}
  (m,n;1)(0,0;\gamma) = (0,0;\gamma)(ma+nc,mb+nd;1)
  \label{eq:(m,n)gamma}
\end{equation}
for $(m,n)\in\Z^2$ and $\gamma\in\Gamma$.

The action of $\Gammatilde$ on $\Xtilde$ is a standard one. An element
$\gamma$ of $\Gamma$ acts on $\Xtilde$ as
\begin{equation*}
  \gamma\cdot(\tau;z;t):=(\check\tau;\check z;\check t)
  \quad
  \text{for $(\tau;z;t)\in\Xtilde$}.
\end{equation*}
This action together with
\eqref{def:action-Z2-on-Xtilde} induces an action of $\Gammatilde$ on
$\Xtilde$.

The actions of $\Gammatilde$ on $G$ and $\g$ are defined via the
action of $\Gamma$ on the Heisenberg group $\Heis_N$, the central extension
of $(\Integer/N\Integer)^2$ defined by
\begin{equation*}
  \Heis_N := \C^{\times} \times (\Z/N\Z)^2.
\end{equation*}
Here the group structure of $\Heis_N$ is given by
\begin{equation*}
  (r;m,n)(r';m',n')
  := (rr'\eps^{nm'};m+m',n+n')
\end{equation*}
for $(r;m,n),(r';m',n')\in\Heis_N$.
The Heisenberg group $\Heis_N$ is isomorphic to the group generated by the
symbols $\hat\alpha$, $\hat\beta$, and $\hat r$ for $r\in\C^\times$ with
defining relations 
\begin{equation}
\begin{aligned}
  &
  \hat\alpha^N = \hat\beta^N = 1, \quad
  \hat\alpha \hat\beta = \eps \hat\beta \hat\alpha, \quad
  \\ &
  \hat r \hat\alpha = \hat\alpha \hat r, \quad
  \hat r \hat\beta = \hat\beta \hat r, \quad
  \widehat{r_1 r_2} = {\hat r}_1 {\hat r}_2
  \quad\text{for $r,r_i\in\C^\times$},
\end{aligned}
\label{Heisenberg:rel}
\end{equation}
where the identification of the group with $\Heis_N$ is given by
\begin{equation}
  (r;m,n) = \hat r {\hat\beta}^m {\hat\alpha}^n.
\label{Heisenberg:presentation}
\end{equation}
Thus the matrices $\alpha$ and $\beta$ given by \eqref{def:alpha-beta} define
a representation of $\Heis_N$ on $\C^N$ by
\begin{equation*}
  \hat\alpha v = \alpha v, \quad
  \hat\beta v = \beta v, \quad
  \hat r v = r v \quad
  \text{for $v\in\C^N$ and $r\in\C^\times$}.
\end{equation*}
Note that this representation is irreducible.  For $\gamma\in\Gamma=SL_2(\Z)$,
we can define the group automorphism $(\cdot)^\gamma$ of $\Heis_N$ as follows:
\begin{itemize}
\item If $N$ is odd, then we put
  \begin{equation}
  \begin{aligned}
    &
    (r;0,0)^\gamma := (r;0,0), \quad
    \\ &
    (1;1,0)^\gamma := (1;a,b), \quad
    (1;0,1)^\gamma := (1;c,d)
  \end{aligned}
  \label{def:gammaOnHeisenberg:odd}
  \end{equation}
  for $r\in\C^\times$.
\item If $N$ is even, then we put
  \begin{equation}
  \begin{aligned}
    &
    (r;0,0)^\gamma := (r;0,0), \quad
    \\ &
    (1;1,0)^\gamma := ((\sqrt\eps)^{ab};a,b), \quad
    (1;0,1)^\gamma := ((\sqrt\eps)^{cd};c,d)
  \end{aligned}
  \label{def:gammaOnHeisenberg:even}
  \end{equation}
  for $r\in\C^\times$, where $\sqrt\eps=\exp(\pi i/N)$, a primitive
  $(2N)$-th root of unity.
\end{itemize}
The fact that this defines an automorphism of $\Heis_N$ follows from the
presentation \eqref{Heisenberg:rel} and \eqref{Heisenberg:presentation}.
Note that this action of $\Gamma$ on $\Heis_N$ induces that of $\Gamma$ on
$(\Z/N\Z)^2$ given by $(m,n)^\gamma=(ma+nc,mb+nd)$ for $(m,n)\in(\Z/N\Z)^2$
and $\gamma\in\Gamma$ (cf.\ \eqref{eq:(m,n)gamma}). 

Twisting the representation of $\Heis_N$ on $\C^N$ by $\gamma\in\Gamma$,
we obtain another irreducible representation of $\Heis_N$ on $\C^N$:
\begin{equation*}
  \Heis_N\times\C^N\to\C^N,
  \quad
  (h,v) \mapsto h^\gamma v.
\end{equation*}
Since the Heisenberg group $\Heis_N$ has a unique irreducible
representation, up to isomorphism, in which $\hat r\in\Heis_N$ for 
$r\in\C^\times$ acts as multiplication by $r$ (the theorem of von Neumann
and Stone), using the Schur lemma, we can find that there is 
$x_\gamma\in GL_N(\C)$, uniquely determined up to scalar multiplications,
such that
\begin{equation}
  h x_\gamma v = x_\gamma h^\gamma v
  \quad
  \text{for $h\in\Heis_N$ and $\gamma\in\Gamma$.}
  \label{def:x-gamma}
\end{equation}
The mapping $\gamma\mapsto x_\gamma$ induces a group homomorphism from
$\Gamma=SL_2(\Z)$ into $PGL_N(\C)$, which does not depend on the choices
of $x_\gamma$'s. In the following, we take $x_\gamma$ from $G=SL_N(\C)$,
which is uniquely determined up to factor $\pm1$ by $\gamma$.

\begin{example}
For $\gamma = \begin{pmatrix} 0 & 1 \\ -1 & 0 \end{pmatrix}$,
$x_\gamma \propto (\eps^{-(a-1)(b-1)})_{a,b=1}^N$.
\end{example}

\begin{example}
When $N$ is odd, we can choose $x_\gamma = 1$ for 
$\gamma \in \Gamma(N) = 
\{\, \gamma\in\Gamma\,|\,\gamma\equiv 1 \mod N\,\}$.
\end{example}

The desired actions of $\Gammatilde$ on $G$ and $\g$ are now defined by
\begin{alignat*}{2}
  & (m,n;\gamma)\cdot g
  := (\beta^m\alpha^nx_\gamma)g(\beta^m\alpha^nx_\gamma)^{-1}
  & \quad &
  \text{for $(m,n;\gamma)\in\Gammatilde$ and $g\in G$},
  \\
  & (m,n;\gamma)\cdot A
  := (\beta^m\alpha^nx_\gamma)A(\beta^m\alpha^nx_\gamma)^{-1}
  & \quad &
  \text{for $(m,n;\gamma)\in\Gammatilde$ and $A\in\g$}.
\end{alignat*}
These are extensions of the actions of $\Z^2$ given by
\eqref{def:Gal-action-on-G} and \eqref{def:Gal-action-on-g}.  

The diagonal actions of $\Gammatilde$ on $\Xtilde\times G$ and
$\Xtilde\times\g$ are defined by
\begin{alignat*}{2}
  & \tilde\gamma\cdot(\tilde x;g)
  := (\tilde\gamma\cdot\tilde x; \tilde\gamma\cdot g)
  & \quad &
  \text{for $(\tilde x;g)\in\Xtilde\times G$ and
    $\tilde\gamma\in\Gammatilde$},
  \\
  & \tilde\gamma\cdot(\tilde x;A)
  := (\tilde\gamma\cdot\tilde x; \tilde\gamma\cdot A)
  & \quad &
  \text{for $(\tilde x;A)\in\Xtilde\times\g$ and
    $\tilde\gamma\in\Gammatilde$}.
\end{alignat*}
Note that these actions also do not depend on the choice of $x_\gamma$.

{}From the actions of $\Gammatilde$ defined above, we obtain the induced
actions of $\Gamma$ on $\X$, $\Gtw_\X$, and $\gtw_\X$ defined by
\begin{alignat*}{2}
  & \gamma\cdot x := \pr_{\Xtilde/\X}(\gamma\cdot\tilde x)
  & \quad &
  \text{for $x=\pr_{\Xtilde/\X}(\tilde x)\in\X$},
  \\
  & \gamma\cdot g^\twist := [\gamma\cdot(\tilde x;g)]
  & \quad &
  \text{for $g^\twist=[(\tilde x;g)]\in\Gtw_\X$},
  \\
  & \gamma\cdot A^\twist := [\gamma\cdot(\tilde x;A)]
  & \quad &
  \text{for $A^\twist=[(\tilde x;A)]\in\gtw_\X$},
\end{alignat*}
where $\tilde x\in\Xtilde$, $(\tilde x;g)\in\X\times G$, and 
$(\tilde x;A)\in\Xtilde\times \g$ are representatives of $x\in\X$,
$g^\twist\in\Gtw_\X$, and $(\tilde x;A)\in\X\times\g$ respectively
and $\gamma\in\Gamma$ is identified with $(0,0;\gamma)\in\Gammatilde$.
(For the definitions of $\X$, $\Gtw_\X$, and $\gtw_\X$, see
\eqref{def:X} and \eqref{def:GtwX-gtwX}.)
The projections $\Xtilde \onto S$, $\gtw_\X \onto \X$, etc.\ are equivariant
with respect to the actions of $\Gamma$.

Moreover we obtain the following induced equivariant actions of
$\gamma\in\Gamma$ on $\T_S$, $\g^D_S$, and $\Vir^D_S$:
\begin{itemize}
\item The biholomorphic map $\gamma:S \to S$ induces the Lie algebra
  isomorphism $(\cdot)^\gamma:\gamma^{-1}\T_S\isoto\T_S$ of vector fields:
  \begin{equation}
    \mu^\gamma
    :=
    \mu_0(\check s)
    \left(
      (c\tau+d)^2 \partial_\tau
      + \sum_{i=1}^L c(c\tau+d)z_i\partial_{z_i}
    \right)
    + \sum_{i=1}^L \mu_i(\check s)(c\tau+d)\partial_{z_i}
    \label{def:mu-gamma}
  \end{equation}
  for $\mu=\mu_0(s)\partial_\tau+\sum_{i=1}^L\mu_i(s)\partial_{z_i}
  \in\gamma^{-1}\T_S$.  (Formally $\mu^\gamma$ is obtained by the substitution
  of $\check s=(\check\tau;\check z)$ in $s=(\tau;z)$.)
\item The Lie algebra isomorphism
  $(\cdot)^\gamma:\gamma^{-1}\hat\g^D_S\isoto\hat\g^D_S$ is defined by
  \begin{equation}
    \bigl(
      (A_i(s;\xi_i))_{i=1}^L;
      f(s)\khat
    \bigr)^\gamma
    :=
    \bigl(
      (x_\gamma^{-1} A_i(\check s;\check\xi_i) x_\gamma)_{i=1}^L;
      f(\check s)\khat
    \bigr)
  \label{def:gSD-gamma}
  \end{equation}
  for $A_i(s;\xi_i)\in\gamma^{-1}(\g\tensor\O_S((\xi_i)))$ and
  $f(s)\in\gamma^{-1}\O_S$. 
\item The Lie algebra isomorphism
  $(\cdot)^\gamma:\gamma^{-1}\Vir^D_S\isoto\Vir^D_S$ is defined by
  \begin{equation}
    \bigl(
      \mu;
      (\theta_i)_{i=1}^L;
      f(s)\chat
    \bigr)^\gamma
    :=
    \bigl(
      \mu^\gamma;
      (
        \mu_0(\check s)c(c\tau+d)\xi_i\partial_{\xi_i}
        + \theta_i(\check s;\check\xi_i)(c\tau+d)\partial_{\xi_i}
      )_{i=1}^L;
      f(\check s)\chat
    \bigr)
    \label{def:Vir-gamma}
  \end{equation}
  for $\mu=\mu_0(s)\partial_\tau+\sum_{i=1}^L\mu_i(s)\partial_{z_i}
  \in\gamma^{-1}\T_S$, 
  $\theta_i=\theta_i(s;\xi_i)\partial_{\xi_i}
  \in\gamma^{-1}(\O_S((\xi_i))\partial_{\xi_i})$,
  and $f(s)\in\gamma^{-1}\O_S$.
\end{itemize}
The formula \eqref{def:Vir-gamma} reflects the fact that the vector field
\begin{equation*}
  \mu_0(\check s)\partial_{\check\tau}
  + \sum_{i=1}^L \mu_i(\check s)\partial_{\check z_i}
  + \theta_i(\check s;\check\xi_i)\partial_{\check\xi_i}
\end{equation*}
represented in the local coordinate $(\check\tau;\check z;\check\xi_i)$ is
equal to
\begin{align*}
  & \mu_0(\check s)
  \left(
    (c\tau+d)^2\partial_\tau
    + \sum_{i=1}^L c(c\tau+d)z_i\partial_{z_i}
    + c(c\tau+d)\xi_i\partial_{\xi_i}
  \right)
  \\
  & + \sum_{i=1}^L \mu_i(\check s)
  (c\tau+d)\partial_{z_i}
  + \theta_i(\check s;\check\xi_i)
  (c\tau+d)\partial_{\xi_i}
\end{align*}
represented in the local coordinate $(\tau;z;\xi_i)$.  Therefore the
following lemma is a direct consequence of the definitions above.

\begin{lem}
\label{lem:gamma-1}
For $\gamma\in\Gamma$, the isomorphisms above satisfy the following:
\begin{enumerate}
\item The isomorphisms induce the following Lie algebra isomorphism:
  \[
     (\cdot)^\gamma:\gamma^{-1}\left(\Vir^D_S\ltimes\hat\g^D_S\right)
     \isoto
     \Vir^D_S\ltimes\hat\g^D_S.
  \]
\item $(\gamma^{-1}\gDPr)^\gamma = \gDPr$.
\item $(\gamma^{-1}\VirDPr)^\gamma=\VirDPr$ and
  $(\gamma^{-1}\TDPr)^\gamma=\TDPr$.
\end{enumerate}
\end{lem}

Under these preparations, we show the modular property of the sheaf of
conformal coinvariants $\sCC(\M)$ and the sheaf of conformal blocks
$\sCB(\M)$ coming from quotients of the Weyl modules.

For each $i=1,\dots,L$, let $V_i$ be a finite-dimensional irreducible
representation of $\g$ and $M_i$ a quotient of the Weyl module
$\Weyl_k(V_i)$.  Put $M:=\bigotimes_{i=1}^LM_i$ and $\M:=M\tensor\O_S$.
We denote by $c_i$ the eigenvalue of the Casimir operator $C_i$ acting on
$V_i$ given by \eqref{def:casimir} and put $\Delta_i:=\kappa^{-1}c_i$.
The Virasoro operator $\rho_i(T[0])$ acting on $M$ is diagonalizable and
each of its eigenvalues is of the form $\Delta_i+m$, where $m$ is a
non-negative integer.  Thus, fixing branches of the holomorphic functions
$(c\tau+d)^{C_i/\kappa}$ on the upper half plane $\UHP$, we obtain an
operator $(c\tau+d)^{\rho(T[0])}=\prod_{i=1}^L(c\tau+d)^{\rho_i(T[0])}$
acting on $\M$, where we put $\rho(T[0]):=\sum_{i=1}^L\rho_i(T[0])$.

For $\gamma\in\Gamma$, define the isomorphism
$(\cdot)^\gamma:\gamma^{-1}\M\isoto\M$ by
\begin{equation}
  v(s)^\gamma := x_\gamma^{-1} (c\tau+d)^{\rho(T[0])} v(\check s)
  \quad
  \text{for $v(s)\in\gamma^{-1}\M$},
  \label{def:v-gamma}
\end{equation}
regarding $x_\gamma^{-1}$ as an automorphism of $\M$ through
the natural diagonal action of $G=SL_N(\C)$ on $M$.
This isomorphism \eqref{def:v-gamma} induces the isomorphism
$(\cdot)^\gamma:\gamma^{-1}(\M^\ast)\isoto\M^\ast$ given by
\begin{equation}
  \Phi(s)^\gamma := x_\gamma^{-1}(c\tau+d)^{\rho^\ast(T[0])}\Phi(\check s)
  \quad\text{for $\Phi(s)\in\gamma^{-1}(\M^\ast)$},
  \label{def:Phi-gamma}
\end{equation}
where we put $\rho^\ast(T[0]):=\sum_{i=1}^L\rho_i^\ast(T[0])$.

\begin{lem}
\label{lem:gamma-2}
For $P\in\Vir^D_S\ltimes\hat\g^D_S$, $v\in\M$, and $\gamma\in\Gamma$, we
have
\begin{equation*}
  (P\cdot v)^\gamma = P^\gamma\cdot v^\gamma
\end{equation*}
\end{lem}

\begin{pf}
Since $\rho_i(T\{\xi_i\partial_{\xi_i}\})=-\rho_i(T[0])$, we have the
following identity of operators acting on $\M$ (cf.\ \eqref{def:Vir-gamma}):
\begin{equation}
\begin{aligned}
  & \mu_0(\check s)(c\tau+d)^2\partial_\tau
  + \sum_{i=1}^L
  \rho_i(T\{ \mu_0(\check s)c(c\tau+d)\xi_i\partial_{\xi_i} \})
  \\
  &\quad =
  (c\tau+d)^{\rho(T[0])}\cdot
  (\mu_0(\check s)(c\tau+d)^2\partial_\tau)\cdot
  (c\tau+d)^{-\rho(T[0])}.
\end{aligned}
\label{Ad(rhoT[0])mu0}
\end{equation}
Similarly, the definition of $\Vir^D_S\ltimes\hat\g^D_S$
\eqref{def:semi-direct-prod} leads to the following identities:
\begin{align}
  & e^{-r\rho_i(T[0])}
  \rho_i(A_i\tensor f_i(\xi_i))
  e^{r\rho_i(T[0])}
  = \rho_i(A_i\tensor f_i(e^r\xi_i)),
  \label{Ad(rhoT[0])Af}
  \\
  & e^{-r\rho_i(T[0])} 
  \rho_i(\theta_i(\xi_i)\partial_{\xi_i})
  e^{r\rho_i(T[0])}
  = \rho_i(\theta_i(e^r\xi_i)e^{-r}\partial_{\xi_i}),
  \label{Ad(rhoT[0])theta}
  \\
  & g\, \rho_i(\theta_i(\xi_i)\partial_{\xi_i})\,
  = \rho_i(\theta_i(\xi_i)\partial_{\xi_i})\, g
  \label{g*theta=theta*g}
\end{align}
for $A_i\tensor f_i(\xi_i)\in\g\tensor\O_S((\xi_i))$, $r\in\C$,
$\theta_i(\xi_i)\partial_{\xi_i}\in\O((\xi_i))\partial_{\xi_i}$, and
$g\in G$.  Hence we have
\begin{equation}
\begin{split}
  & (c\tau+d)^{-\rho(T[0])}x_\gamma
  \cdot
  \bigl(
    (A_i(s;\xi_i))_{i=1}^L;
    f(s)\khat
  \bigr)^\gamma
  \cdot
  x_\gamma^{-1}(c\tau+d)^{\rho(T[0])}
  \\
  & \quad = 
  \bigl(
    (A_i(\check s;\xi_i))_{i=1}^L;
    f(\check s)\khat
  \bigr),
\end{split}
\end{equation}
for $A_i(s;\xi_i)\in\gamma^{-1}(\g\tensor\O_S((\xi_i)))$,
$f(s)\in\gamma^{-1}\O_S$ because of \eqref{def:gSD-gamma}, and
\begin{equation}
\begin{split}
  & (c\tau+d)^{-\rho(T[0])}x_\gamma
  \cdot
  \bigl(
    \mu;
    (\theta_i(s;\xi_i))_{i=1}^L;
    g(s)\chat
  \bigr)^\gamma
  \cdot
  x_\gamma^{-1}(c\tau+d)^{\rho(T[0])}
  \\
  & \quad = 
  \bigl(
    \mu^\gamma;
    (\theta_i(\check s;\xi_i)\partial_{\xi_i})_{i=1}^L;
    g(\check s)\chat
  \bigr)
\end{split}
\end{equation}
for
$\mu=\mu_0(s)\partial_\tau+\sum_{i=1}^L\mu_i(s)\partial_{\xi_i}
\in\gamma^{-1}\T_S$, 
$\theta_i(s;\xi_i)\partial_{\xi_i}
\in\gamma^{-1}(\O_S((\xi_i))\partial_{\xi_i})$,
$g(s)\in\gamma^{-1}\O_S$ because of \eqref{def:Vir-gamma}.
These formulae prove the lemma. 
\end{pf}

{}From \lemref{lem:gamma-1} and \lemref{lem:gamma-2} follow the modular
invariance of $\sCC(\M)$ and $\sCB(\M)$, and the modular transformations
of connections $D$ \eqref{def:conn-cc} and $D^\ast$ \eqref{def:conn-cb}.
The results are summarized in the following theorem.

\begin{thm}
\label{thm:modular-inv-1}
For each $i=1,\dots,L$, let $V_i$ be a finite-dimensional irreducible
representation of $\g$ and $M_i$ a quotient of the Weyl module $\Weyl_k(V_i)$.
Put $M:=\bigotimes_{i=1}^L M_i$ and $\M:=M\otimes\O_S$. 
Let $\gamma$ be in $\Gamma=SL_2(\Z)$.  Then the isomorphisms
$(\cdot)^\gamma:\gamma^{-1}\M\isoto\M$ and
$(\cdot)^\gamma:\gamma^{-1}(\M^\ast)\isoto\M^\ast$ defined by
\eqref{def:v-gamma} and \eqref{def:Phi-gamma} induce the isomorphisms
\begin{alignat*}{2}
  & (\cdot)^\gamma : \gamma^{-1}\sCC(\M) \isoto \sCC(\M),
  & \quad &
  v(s) \mapsto x_\gamma^{-1} (c\tau+d)^{\rho(T[0])} v(\check s),
  \\
  & (\cdot)^\gamma : \gamma^{-1}\sCB(\M) \isoto \sCB(\M),
  & \quad &
  \Phi(s) \mapsto x_\gamma^{-1}(c\tau+d)^{\rho^\ast(T[0])} \Phi(\check s),
\end{alignat*}
where $x_\gamma$ is an element of $G$ satisfying \eqref{def:x-gamma}.
Furthermore these isomorphisms correspond local flat sections with
respect to the connections $D$ and $D^\ast$ to local flat sections.  Namely,
denoting the subsheaf of local flat sections of $\sCC(\M)$ and $\sCB(\M)$ by
$\sCC(\M)^D$ and $\sCB(\M)^{D^\ast}$ respectively, we obtain the following
induced isomorphisms: 
\begin{align*}
  & (\cdot)^\gamma : \gamma^{-1}(\sCC(\M)^D) \isoto \sCC(\M)^D,
  \\
  & (\cdot)^\gamma : \gamma^{-1}(\sCB(\M)^{D^\ast}) \isoto \sCB(\M)^{D^\ast}.
\end{align*}
The modular transformations of the connections are represented as
\begin{align}
  &
  \begin{aligned}
    D_{\mu^\gamma}
    &=
    \mu^\gamma
    + \sum_{i=1}^L 
    \rho_i(T\{ \theta_i(\check s;\check\xi_i)\partial_{\check\xi_i} \})
    - \mu_0(\check s)\,c(c\tau+d)\rho(T[0])
    \\
    &= \left[ D_\mu \right]_{s\mapsto\check s}
    - \mu_0(\check s)\,c(c\tau+d)\rho(T[0]),
  \end{aligned}
  \label{modular:D}
  \\ &
  \begin{aligned}
    D^\ast_{\mu^\gamma}
    &=
    \mu^\gamma
    + \sum_{i=1}^L 
    \rho^\ast_i(T\{ \theta_i(\check s;\check\xi_i)\partial_{\check\xi_i} \})
    - \mu_0(\check s)\,c(c\tau+d)\rho^\ast(T[0])
    \\
    &= \left[ D^\ast_\mu \right]_{s\mapsto\check s}
    - \mu_0(\check s)\,c(c\tau+d)\rho^\ast(T[0]),
  \end{aligned}
  \label{modular:D*}
\end{align}
where $\mu=\mu_0(s)\partial_\tau+\sum_{i=1}^L\mu_i(s)\partial_{z_i}
\in\gamma^{-1}\T_S$,
$\theta_i=\theta_i(s;\xi_i)\partial_{\xi_i}
\in\gamma^{-1}(\O_S((\xi_i))\partial_{\xi_i})$,
$(\mu;(\theta_i)_{i=1}^L;0)\in\gamma^{-1}(\VirDPr)$, and
$[\,\cdot\,]_{s\mapsto\check s}$ denotes the substitution of $\check s$ in
$s$, namely, $[\mu+A(s)]_{s\mapsto\check s}=\mu^\gamma+A(\check s)$ for 
$\mu\in\gamma^{-1}\T_S$ and $A(s)\in\gamma^{-1}\END_{\O_S}(\sCC(\M))$ or
$A(s)\in\gamma^{-1}\END_{\O_S}(\sCB(\M))$.
\end{thm}
Indeed \eqref{modular:D} and \eqref{modular:D*} follow from the explicit
expressions \eqref{def:conn-cc:explicit} and \eqref{def:conn-cb:explicit}
and the definition \eqref{def:Vir-gamma}.

Using \propref{sheaf-restriction}, we can identify $\sCC_k(V)$ and $\sCB_k(V)$
with $V\tensor\O_S$ and $V^\ast\tensor\O_S$ respectively.  Put
$\Delta:=\sum_{i=1}^L\Delta_i$.  Then the theorem above implies the following
corollary.

\begin{cor}
\label{cor:modular-inv-2}
Then, under the identifications $\sCC_k(V)=V\tensor\O_S$ and
$\sCB_k(V)=V^\ast\tensor\O_S$, the isomorphisms
$(\cdot)^\gamma:\gamma^{-1}\sCC_k(V)\isoto\sCC_k(V)$ and
$(\cdot)^\gamma:\gamma^{-1}\sCB_k(\M)\isoto\sCB_k(V)$ are of the following
forms: 
\begin{alignat*}{2}
  & (\cdot)^\gamma : \gamma^{-1}(V\tensor\O_S) \isoto V\tensor\O_S,
  & \quad &
  v(s) \mapsto x_\gamma^{-1} (c\tau+d)^\Delta v(\check s),
  \\
  & (\cdot)^\gamma :  \gamma^{-1}(V^\ast\tensor\O_S) \isoto V^\ast\tensor\O_S,
  & \quad &
  F(s) \mapsto x_\gamma^{-1} (c\tau+d)^{-\Delta} F(\check s).
\end{alignat*}
The modular transformations of the connections expressed as in \thmref{kz} are
represented as
\begin{align}
  & D_{\mu^\gamma}
  =
  \left[ D_\mu \right]_{s\mapsto\check s}
  - \mu_0(\check s)c(c\tau+d)\Delta,
\label{modular:D:fin}
  \\
  & D^\ast_{\mu^\gamma}
  =
  \left[ D^\ast_\mu \right]_{s\mapsto\check s}
  + \mu_0(\check s)c(c\tau+d)\Delta,
\label{modular:D*:fin}
\end{align}
where $\mu=\mu_0(s)\partial_\tau+\sum_{i=1}^L\mu_i(s)\partial_{\xi_i}
\in\gamma^{-1}\T_S$ and 
$[\,\cdot\,]_{s\mapsto\check s}$ denotes the substitution of $\check s$ in
$s$, namely, $[\mu+A(s)]_{s\mapsto\check s}=\mu^\gamma+A(\check s)$ for 
$\mu\in\gamma^{-1}\T_S$ and $A(s)\in\gamma^{-1}\END_{\O_S}(V\tensor\O_S)$ or
$A(s)\in\gamma^{-1}\END_{\O_S}(V^\ast\tensor\O_S)$.
\end{cor}

\begin{example}
Applying \eqref{modular:D*:fin} to $\mu = \partial_{z_i}$ and 
$\mu = \partial_\tau$, we obtain
\begin{align*}
  & D^\ast_{(\PD{}{z_i})^\gamma}
  =
  \left[ D^\ast_{\PD{}{z_i}} \right]_{s\mapsto\check s},
  \\
  & D^\ast_{(\PD{}{\tau})^\gamma}
  =
  \left[ D^\ast_{\PD{}{\tau}} \right]_{s\mapsto\check s}
  + c(c\tau+d)\Delta,
\end{align*}
which were found by Etingof \cite{eti:94}, \S4.
\end{example}

\section{Concluding remarks}
\label{conclusion}

We have examined a twisted WZW model on elliptic curves which gives the
XYZ Gaudin model at the critical level and Etingof's elliptic KZ equations
at the off-critical level. We make several comments on the related
interesting problems to be solved.

\subsection*{Factorization}
We have studied the twisted WZW model only on a family of smooth pointed
elliptic curves, but as in \cite{tu-ue-ya:89} we can also consider the
model on a family of stable pointed elliptic curves. By extending the
connections acting on the sheaves of conformal blocks to those with
regular singularities at the boundary of the family, we shall be able to
establish the equivalence between our geometric approach and Etingof's
approach by means of twisted traces of the products of twisted vertex
operators on the Riemann sphere.  Furthermore we shall be able to obtain a
dimension formulae for the spaces of conformal blocks. A detailed
investigation shall be given in a forthcoming paper.

\subsection*{Generalization to higher genus}
Bernard generalized the KZB equations to higher genus Riemann surfaces in
\cite{ber:88-2}. In \cite{fel:96}, Felder established the geometric
interpretation of the KZB equations on Riemann surfaces by the non-twisted
WZW models clarifying the notion of the dynamical $r$-matrices in higher
genus cases. Our formulation for the twisted WZW model is also valid for
arbitrary Riemann surfaces. See \appref{higher-genus} for details.

\subsection*{Discretization}
Felder calls his interpretation of the KZB equations in \cite{fel:96}
``the first step of the `St.\ Petersburg $q$-deformation recipe' in higher
genus cases''. We hope that our twisted WZW model on elliptic curves can
also be $q$-deformed. The resulting ``elliptic $q$-KZ equations'', for
example, would be related to the difference equations proposed in
\cite{takebe:97}.

\subsection*{Intertwining vectors}
Note that the Boltzmann weights of the $A^{(1)}_{N-1}$ face model can be
expressed by the elliptic quantum $R$-matrix and the intertwining vectors
(\cite{bax:73}, \cite{JMO:87}, \cite{djkmo}).  Therefore it can be
expected that there exists a quasi-classical limit of the intertwining
vectors by means of which the relation of the non-twisted and twisted WZW
models on elliptic curves will be clarified.  

In addition, the intertwining vectors play an important role in constructing
Bethe vectors of the XYZ spin chain models (\cite{bax:73},
\cite{takebe:xyz}) and the integral solution of the difference equations
in \cite{takebe:97}. They are introduced as a kind of technical tools
there, but our approach from the twisted WZW model, i.e., from the
classical limit should reveal their algebro-geometric meaning.


\section*{Acknowledgements}
TT is supported by Postdoctoral Fellowship for Research abroad of Japan
Society for the Promotion of Science.  He expresses his gratitude to
Benjamin~Enriquez, Giovanni~Felder, Edward~Frenkel, Ian~Grojnowski,
Takeshi~Ikeda, Feodor~Malikov, Hirosi~Ooguri, Nicolai~Reshetikhin for
comments and discussions.  The authors also express their gratitude to
Yoshifumi Tsuchimoto, who explained to them a detail of the proof of the
main theorem in \cite{tsuchimoto:93}.  Parts of this work were done while
TT was visiting the Department of Mathematics of the University of
California at Berkeley, the Department of Mathematics of Kyoto University,
the Erwin Schr\"odinger Institute for Mathematical Physics, the Landau
Institute for Theoretical Physics and Centre Emile Borel - Institut Henri
Poincar\'e - UMS 839, CNRS/UPMC. He thanks these institutes for
hospitality.

%
%
%
%
\appendix
\section{Theta functions with characteristics}
\label{theta-func}

Here we collect properties of theta functions of one variable used in this
paper.

Following \cite{mum}, we use the notation:
\begin{equation}
    \theta_{\kappa,\kappa'}(t;\tau)
    =
    \sum_{n\in\Z} 
    e^{\pi i (n+\kappa)^2 \tau + 2 \pi i (n+\kappa)(t+\kappa')}.
\label{def:theta}
\end{equation}
for the theta functions with characteristics. Here $t$ is a complex
number, $\tau$ belongs to the upper half plane $\UHP$ and $\kappa$,
$\kappa'$ are rational numbers. They are related with each other by shifts
of $t$:
\begin{equation}
    \theta_{\kappa_1+\kappa_2, \kappa'_1+\kappa'_2}(t;\tau)
    =
    e^{\pi i \kappa_2^2 \tau + 2\pi i \kappa_2(t+\kappa'_1+\kappa'_2)}
    \theta_{\kappa_1, \kappa'_1}(t + \kappa_2 \tau + \kappa'_2;\tau).
\label{theta:shift-of-char}
\end{equation}
Since the zero set of $\theta_{00}(t;\tau)$ is 
$\{1/2 + \tau/2 \} + \Z + \Z \tau$, the zero set of
$\theta_{\kappa,\kappa'}(t;\tau)$ is
\begin{equation}
    \left\{\frac12 -\kappa + \left(\frac12 - \kappa'\right)\tau \right\}
    + \Z + \Z \tau,
\label{theta:zero}
\end{equation}
because of \eqref{theta:shift-of-char}. 
Fundamental properties of the function $\theta_{\kappa,\kappa'}$ are the
quasi-periodicity with respect to $t$:
\begin{equation}
\begin{aligned}
    \theta_{\kappa,\kappa'}(t+1;\tau) 
    &= e^{2\pi i \kappa} \theta_{\kappa,\kappa'}(t;\tau),
\\
    \theta_{\kappa,\kappa'}(t+\tau;\tau) 
    &= e^{- \pi i \tau - 2\pi i (t + \kappa')}
       \theta_{\kappa,\kappa'}(t;\tau),
\end{aligned}
\label{theta:quasi-period}
\end{equation}
and the automorphic property:
\begin{equation}
\begin{aligned}
    \theta_{\kappa,\kappa'}(t;\tau+1)
    &=
    e^{- \pi i \kappa(\kappa+1)}
    \theta_{\kappa,\kappa+\kappa'+ \frac12}(t;\tau),
\\
    \theta_{\kappa,\kappa'}\left(\frac{t}{\tau};-\frac{1}{\tau}\right)
    &=
    (-i\tau)^{1/2}
    e^{2 \pi i \kappa \kappa'} e^{\pi i t^2/\tau}
    \theta_{\kappa',-\kappa}(t;\tau).
\end{aligned}
\label{theta:modular}
\end{equation}
The formulae
\begin{align}
    \theta_{-\kappa,-\kappa'}(t;\tau)
    &=
    \theta_{\kappa,\kappa'}(-t;\tau),
\label{theta:odd}
\\
    \theta_{\kappa+m,\kappa'+n}(t;\tau) 
    &= e^{2\pi i \kappa n}
    \theta_{\kappa,\kappa'}(t;\tau),
\label{theta:char-period}
\end{align}
are easily deduced from the definition \eqref{def:theta},
where $m$ and $n$ are integers.

We use mostly the following special characteristics. Let $N \geqq 2$ be a
natural number and $a,b$ arbitrary integers. We denote
\begin{equation}
    \theta_{[a,b]}(t;\tau)
    :=
    \theta_{\frac{a}{N} - \frac12, -\frac{b}{N} + \frac12}(t;\tau).
\label{def:theta[a,b]}
\end{equation}
The standard abbreviations,
\begin{equation*}
    \theta_{[a,b]} := \theta_{[a,b]}(0;\tau), \quad
    \theta'_{[a,b]} := \left. \od{}{t}\theta_{[a,b]}(t;\tau)\right|_{t=0},
\end{equation*}
and $\theta''_{[a,b]}$ etc.\ likewise defined are also used.

\begin{lem}
\label{jacobi}
\begin{equation}
    \frac{N^2-1}{6} \frac{\theta'''_{[0,0]}}{\theta'_{[0,0]}}
    -
    \frac{1}{2}
    \sum_{(a,b) \neq (0,0)}
    \frac{\theta_{[a,b]}''}{\theta_{[a,b]}}
    = 0,
\end{equation}
where the indices in the second term run through $a = 0, \ldots, N-1$,
$b = 0, \ldots, N-1$ and $(a,b) \neq (0,0)$.
\end{lem}

This is a generalization of the well-known formula
\begin{equation}
    \frac{\theta_{1/2,1/2}'''}{\theta_{1/2,1/2}'}
    =
    \frac{\theta_{1/2,  0}''}{\theta_{1/2,  0}} +
    \frac{\theta_{  0,  0}''}{\theta_{  0,  0}} +
    \frac{\theta_{  0,1/2}''}{\theta_{  0,1/2}},
\end{equation}
which is the case $N=2$ of \lemref{jacobi}.
\begin{pf}
It is easy to show that
\begin{equation}
    N \prod_{a=0}^{N-1} \prod_{b=0}^{N-1} \theta_{[a,b]}(t;\tau)
    =
    \left(\prod_{(a,b)\neq(0,0)} \theta_{[a,b]}\right)
    \theta_{[0,0]}(Nt;\tau).
\label{lem:prod}
\end{equation}
In fact we have only to compare the periodicity and zeros of the both sides by
using \eqref{theta:quasi-period} and \eqref{theta:zero}. The over-all
coefficient can be determined by the first term (namely the coefficient of
$t$) in the Taylor expansion around $t=0$.

The coefficients of $t^2$ of the Taylor expansion of \eqref{lem:prod}
give
\begin{equation*}
    \sum_{(a,b)\neq(0,0)} \frac{\theta'_{[a,b]}}{\theta_{[a,b]}} = 0,
\end{equation*}
and using this equality, we can rewrite the terms of order $t^3$ in
\eqref{lem:prod} as follows:
\begin{equation}
    \frac{N^2-1}{6} \frac{\theta'''_{[0,0]}}{\theta'_{[0,0]}}
    -
    \frac{1}{2}
    \sum_{(a,b) \neq (0,0)}
    \frac{\theta_{[a,b]}''}{\theta_{[a,b]}}
    = 
    - \sum_{(a,b) \neq (0,0)}
    \left(\frac{\theta'_{[a,b]}}{\theta_{[a,b]}}\right)^2.
\label{jacobi:tmp}
\end{equation}
Therefore, in order to prove \lemref{jacobi}, we have to show that the
right-hand side of \eqref{jacobi:tmp} is zero. Let us denote it by
$f(\tau)$ as a function of $\tau$. It has following properties:
\begin{itemize}
\item
$f(\tau)$ is a holomorphic function on the upper half plane $\UHP$.
(\eqref{theta:zero})

\item
$f(\tau)$ is bounded when $\Im\tau \to +\infty$. (\eqref{def:theta})

\item
$f(\tau+1) = f(\tau)$, $f(-1/\tau) = \tau^2 f(\tau)$.
(\eqref{theta:modular}, \eqref{theta:char-period})
\end{itemize}
Hence $f(\tau)$ is an integral modular form of weight 2 and level 1,
which is nothing but zero. (See, for example, Th\'eor\`em 4 (i), \S3
Chapitre VII of \cite{ser:70}, Proposition 2.26 of
\cite{shimura:71} or Theorem 14 in Chapter II of \cite{sch:74}.)  This
proves the lemma.
\end{pf}

%
%
%
%
\section{The Kodaira-Spencer map of a family of Riemann surfaces}
\label{kodaira-spencer}

Let $\pr:\X\to S$ be a family of compact Riemann surfaces over a complex
manifold $S$ and $q_i:S\to\X$ a holomorphic section of $\pr$ for each
$i=1,\dots,L$.  Put $Q_i:=q_i(S)$ and assume that $Q_i\cap Q_j=\emptyset$ if
$i\ne j$.  Then $D:=\bigcup_{i=1}^L Q_i$ is a divisor of $\X$ \'etale over
$S$.  We call $(\pr:\X\to S; q_1,\dots,q_L)$ a family of pointed compact
Riemann surfaces.  Denote by $\T_\X(-\log D)$ the sheaf of vector fields
tangent to $D$.  As in \secref{def-sheaf-virasoro}, let $\T_{\X,\pr}(-\log D)$
be the inverse image of $\pr^{-1}\T_S\subset\pr^\ast\T_S$ in $\T_\X(-\log D)$.
Then we obtain the following short exact sequence:
\begin{equation*}
  0 \to \T_{\X/S}(-D) \to \T_{\X,\pr}(-\log D) \to \pr^{-1}\T_S \to 0,
\end{equation*}
which is a Lie algebra extension.  The derived direct image of this sequence
produces the following long exact sequence:
\begin{equation*}
  \cdots 
  \to \pr_\ast\T_{\X,\pr}(-\log D) 
  \to \T_S
  \to R^1\pr_\ast\T_{\X/S}(-D)
  \to R^1\pr_\ast\T_{\X,\pr}(-\log D)
  \to \cdots.
\end{equation*}
The connecting homomorphism $\T_S\to R^1\pr_\ast\T_{\X/S}(-D)$ is called the
{\em Kodaira-Spencer map} of the family $(\pr:\X\to S;q_1,\dots,q_L)$.

For an $\O_\X$-module $\F$ and a closed analytic subspace $Z$ of $\X$, denote
by $\F^\wedge_Z$ the completion of $\F$ at $Z$.  Consider the following exact
sequences:
\begin{align*}
  & 0
  \to (\T_{\X/S}(-D))^\wedge_D
  \to (\T_\X(-\log D))^\wedge_D
  \to (\pr^\ast\T_S)^\wedge_D
  \to 0,
  \\
  & 0
  \to (\T_{\X/S}(\ast D))^\wedge_D
  \to (\T_\X(\ast D))^\wedge_D
  \to (\pr^\ast\T_S)^\wedge_D
  \to 0.
\end{align*}
As above, the inverse images of
$\pi^{-1}\T_S|_D\subset(\pr^\ast\T_S)^\wedge_D$
in $(\T_\X(-\log D))^\wedge_D$ and $(\T_\X(\ast D))^\wedge_D$
is denoted by $\T_\pr(-\log D)$ and $\T_\pr(\ast D)$ respectively.
Then we obtain the Lie algebra extensions below:
\begin{align*}
  & 0
  \to (\T_{\X/S}(-D))^\wedge_D
  \to \T_\pr(-\log D)
  \to \pr^{-1}\T_S|_{D}
  \to 0,
  \\
  & 0
  \to (\T_{\X/S}(\ast D))^\wedge_D
  \to \T_\pr(\ast D)
  \to \pr^{-1}\T_S|_{D}
  \to 0.
\end{align*}
The direct images of these exact sequences to $S$ are also exact. 

Then we have the following commutative diagram:
\begin{equation}
\minCDarrowwidth 5mm
\begin{CD}
  @. 0 @. 0 @. 0 @. (\sharp) @.
  \\
  @. @VVV @VVV @VVV @VVV @.
  \\
  0
  @>>> \pr_\ast\T_{\X/S}^{-D}
  @>>> T_\Prin^{\ast D} \oplus T^{-D}
  @>>> T^{\ast D}
  @>>> R^1\pr_\ast\T_{\X/S}^{-D}
  @>>> 0,
  \\
  @. @VVV @VVV @VVV @. @.
  \\
  0
  @>>> \pr_\ast\T_{\X,\pr}^{-\log D}
  @>>> T_{\pr,\Prin}^{\ast D} \oplus T_\pr^{-\log D}
  @>>> T_\pr^{\ast D}
  @.
  @.
  \\
  @. @VVV @VVV @VVV @. @.
  \\
  0
  @>>> \T_S
  @>>> \T_S\oplus(\T_S)^L
  @>>> (\T_S)^L       
  @.
  @.
  \\
  @. @VVV @VVV @VVV @. @.
  \\
  @. (\sharp) @. 0 @. 0 @. @.
\end{CD}
\label{cd:kod-sp}
\end{equation}
where we put
\begin{alignat*}{3}
  & \T_{\X/S}^{-D} := \T_{\X/S}(-D), \quad &
  & \T_{\X,\pr}^{-\log D} := \T_{\X,\pr}(-\log D), \quad &
  \\
  & T_\Prin^{\ast D}:= \pr_\ast\T_{\X/S}(\ast D), \quad &
  & T^{-D}  := \pr_\ast(\T_{\X/S}(-D))^\wedge_D, \quad &
  & T^{\ast D}:= \pr_\ast(\T_{\X/S}(\ast D))^\wedge_D,
  \\
  & T_{\pr,\Prin}^{\ast D} := \pr_\ast\T_{\X,\pi}(\ast D), \quad &
  & T_\pr^{-\log D} := \pr_\ast\T_\pr(-\log D), \quad &
  & T_\pr^{\ast} := \pr_\ast\T_\pr(\ast D).
\end{alignat*}
The horizontal and vertical sequences in the diagram \eqref{cd:kod-sp} are all
exact and the arrow from $\T_S$ to
$R^1\pr_\ast\T_{\X/S}^{-D}=R^1\pr_\ast\T_{\X/S}(-D)$ through $(\sharp)$ is
nothing but the Kodaira-Spencer map, which is described as follows.  For
$\mu\in\T_S$, chasing the diagram above, we can choose $(a_\Prin,a_+)$,
$\alpha$ and $[\alpha]$ in order:
\begin{enumerate}
\item $(a_\Prin,a_+)\in T_{\pr,\Prin}^{\ast D}\oplus T_\pr^{-\log D}$, whose
  image in $\T_S\oplus(\T_S)^L$ is equal to
  $(\mu;(\mu)_{i=1}^L)\in\T_S\oplus(\T_S)^L$;
\item $\alpha\in T^{\ast D}$, whose image in $T_\pr^{\ast D}$ is equal to
  $a_\Prin - a_+\in T_\pr^{\ast D}$;
\item $[\alpha]\in R^1\pr_\ast\T_{\X/S}^{-D}$, which is the image of $\alpha$
  in $R^1\pr_\ast\T_{\X/S}^{-D}$.
\end{enumerate}
Then the cohomology class $[\alpha]\in R^1\pr_\ast\T_{\X/S}(-D)$ does not
depend on the choice of $(a_\Prin,a_+)$ and $\alpha$.  The Kodaira-Spencer
map sends $\mu\in\T_S$ to $[\alpha]\in R^1\pr_\ast\T_{\X/S}(-D)$.  Recall
that $T_\Prin^{\ast D}$ and $T_{\pr,\Prin}^{\ast D}$ are denoted by
$\TDPr$ and $\VirDPr$ respectively in \secref{flat-connections}. The short
exact sequence \eqref{eq:LieAlgExtOfTSByTDPr} is included in the second
vertical exact sequence in \eqref{cd:kod-sp} and the lifting from $\T_S$
to $\VirDPr$ (cf.\ \secref{const-flat-conn}) essentially corresponds to the
operation 1 above.  Hence we can see that the constructions in
\secref{flat-connections} originate in the above description of the
Kodaira-Spencer map.

However we considered not only $T_\pr^{*D}$ but also its extension $\Vir^S_D$
by $\O_S\chat$. This is a difference between the description of the
Kodaira-Spencer map and the constructions in \secref{flat-connections}.  We
remark that Beilinson and Schechtman give the intrinsic 
(i.e., coordinate-free) description of the Virasoro algebras in \cite{BS}.

%
%
%
%
\section{On a formulation for higher genus Riemann surfaces}
\label{higher-genus}

In this appendix, we shall comment on a formulation of twisted
Wess-Zumino-Witten models on higher genus Riemann surfaces.

Let $\pr:\X\to S$, $q_i:S\to\X$, $Q_i=q_i(S)$, and $D=\sum Q_i$ be the same as
\appref{kodaira-spencer}. Suppose that, for each $i$, we can take a
holomorphic function $\xi_i$ on an open neighborhood $U_i$ of $Q_i$ with the
property that the mapping $U_i\to S\times\xi_i(U_i)$ given by
$x\mapsto(\pr(x),\xi_i(x))$ is biholomorphic and $\xi_i(Q_i)=\{0\}$.

Then, in precisely the same way as \secref{flat-connections}, we can
define $\T_{\X,\pr}(\ast D)$, $\T_{\X/S}(\ast D)$, $\Vir^D_S$, 
$\TDPr$, $\VirDPr$, etc. We define the action of $\pr_\ast\T_{\X,\pr}(\ast D)$
on $\Vir^D_S$ by
\begin{equation}
  a_\Prin\cdot\alpha := [a_\Prin, \alpha]
  \quad\text{for $a_\Prin\in\VirDPr$, $\alpha\in\Vir^D_S$},
  \label{def:Vir-ActsOnVir}
\end{equation}
where $\VirDPr$ is identified with its image in $\Vir^D_S$ and the bracket
of the right-hand side is the Lie algebra structure of $\Vir^D_S$.

We remark that the embedding $\VirDPr$ and that of $\TDPr$
into $\Vir^D_S$ are not always Lie algebra homomorphisms and the action of
$\VirDPr$ on $\Vir^D_S$ does not always preserve $\TDPr$. Thus we must add a
supplementary structure on the Riemann surface. 

Suppose that we can take an open covering $\{U_\lambda\}_{\lambda\in\Lambda}$
of $\X$ and a family $\{\xi_\lambda:U_\lambda\to \C\}_{\lambda\in\Lambda}$ of
holomorphic functions satisfying the following properties:
\begin{enumerate}
\item For each $\lambda\in\Lambda$, the mapping $U_\lambda\to S\times\C$ given
  by $x\mapsto(\pr(x),\xi_\lambda(x))$ is a biholomorphic mapping from
  $U_\lambda$ onto an open subset of $S\times\C$.
\item For any $\lambda,\lambda'\in\Lambda$ with 
  $U_\lambda\cap U_\lambda'\ne\emptyset$, there exists
  $a,b,c,d\in\O_S(S)$ with the property that
  $\xi_{\lambda'}=(a\xi_\lambda+b)/(c\xi_\lambda+d)$ on $U_\lambda\cap
  U_\lambda'$.
\end{enumerate}
We call $\{\xi_\lambda\}_{\lambda\in\Lambda}$ a {\em projective structure}
on the family $\pr:\X\to S$ of Riemann surfaces.  Moreover assume that
$\{\xi_\lambda\}_{\lambda\in\Lambda}\cup\{\xi_i\}_{i=1}^L$ is also a projective
structure on the family.

\begin{lem}
\label{lem:Vir-ActsOnT-}
Under the assumption above, the action of $\VirDPr$ on $\Vir^D_S$ preserves
$\TDPr$ and in particular the embedding $\TDPr\injto\Vir^D_S$ is a Lie
algebra homomorphism.
\end{lem}
\begin{pf}
It suffices to show that $\cvir(\theta,\eta)=0$ for
$(\mu;\theta;0)\in\VirDPr$ and $(0;\nu;0)\in\TDPr$.  For this purpose, as
in the proof of \lemref{Vir-toVirIsLieAlgHom}, it is enough to show that, for
$(\mu;\theta;0)\in\VirDPr$ and $(0;\nu;0)\in\TDPr$, we can take
$\omega\in\pr_\ast\Omega_{\X/S}^1(\ast D)$ with the property that 
$\omega = \theta'''(\xi_i)\eta(\xi_i)\,d\xi_i$ near $Q_i$.  To do this, we
calculate the transformation property of
$\theta'''(\xi_\lambda)\eta(\xi_\lambda)\,d\xi_\lambda$ under coordinate
changes.  Take any $\xi$, $\zeta$ from
$\{\xi_\lambda\}_{\lambda\in\Lambda}\cup\{\xi_i\}_{i=1}^L$.  Then there is
$a,b,c,d\in\O_S(S)$ with $\zeta=(a\xi+b)/(c\xi+d)$.  Fix a local coordinate
$s=(s_i)_{i=1}^M$ on $S$.  Then we obtain two local coordinates $(s;\xi)$ and
$(s;\zeta)$ on $\X$. Let $a_\Prin$ be in $\pr_\ast\T_{\X,\pr}(\ast D)$ and
$\alpha$ in $\pr_\ast\T_{\X/S}(\ast D)$.  Then we can represent $a_\Prin$ and
$\alpha$ in the two local coordinates $(s;\xi)$ and $(s;\zeta)$:
\begin{equation*}
  a_\Prin =
  \begin{cases}
    \mu + \theta^\xi(s;\xi)\partial_\xi & \text{in $(s;\xi)$}, \\
    \mu + \theta^\zeta(s;\zeta)\partial_\zeta & \text{in $(s;\zeta)$},
  \end{cases}
  \quad
  \alpha =
  \begin{cases}
    \eta^\xi(s;\xi)\partial_\xi & \text{in $(s;\xi)$}, \\
    \eta^\zeta(s;\zeta)\partial_\zeta & \text{in $(s;\zeta)$},
  \end{cases}
\end{equation*}
where $\mu=\sum_{i=1}^M\mu_i(s)\partial_{s_i}$.  Then a straightforward
calculation shows that
\begin{equation*}
    \pd{^3\theta^\xi(s;\xi)}{\xi^3} \eta^\xi(s;\xi) \,d\xi
  = \pd{^3\theta^\zeta(s;\zeta)}{\zeta^3} \eta^\zeta(s;\zeta) \,d\zeta.
\end{equation*}
Hence there is a unique $\omega\in\pr_\ast\Omega_{\X/S}^1(\ast D)$ such that
the
representation of $\omega$ under the coordinate $(s;\xi)$ is equal to 
$(\partial_\xi^3\theta^\xi(s;\xi))\eta^\xi(s;\xi)\,d\xi$
for any $\xi\in\{\xi_\lambda\}_{\lambda\in\Lambda}\cup\{\xi_i\}_{i=1}^L$.
Thus we have completed the proof.
\end{pf}

\begin{example}
Assume that $\pr:\X\to S$ denotes the family of elliptic curves defined in
\secref{family-of-curves}.  Then the local coordinate $t$ along the fibers
gives a projective structure on the family.  For each holomorphic section
$q:S\to\X$ of $\pr:\X\onto S$, put $\xi_q:=t-q$, which is regarded as a
holomorphic function on a sufficiently small open neighborhood of $q(S)$.
Then the family $\{\xi_q\}$ is a projective structure on the family and
contains $\{\xi_i\}_{i=1}^L=\{\xi_{q_i}\}_{i=1}^L$.
\end{example}

\begin{example}
A family of compact Riemann surfaces given by the Schottky
pa\-ram\-e\-tri\-za\-tion possess a natural projective structure, because each
compact Riemann surface in that family is represented as the quotient space of
the punctured Riemann sphere by the action of a Schottky group, which is
generated by a certain finite set consisting of fractional linear
transformations.  The Schottky parametrization is used in \cite{ber:88-2}.
\end{example}

We can generalize the statement of \lemref{lem:Vir-ActsOngDS-} for the family
of Riemann surfaces with projective structures. Let $\gtw_\X$ denote an
$\O_\X$-Lie algebra which is locally $\O_\X$-free of finite rank with
holomorphic flat connection $\nabla$ and suppose that the action of a vector
field on $\gtw_\X$ via the connection is a Lie algebra derivation:
\begin{equation*}
  \nabla[A,B] = [\nabla A, B] + [A, \nabla B]
  \quad
  \text{for $A,B\in\gtw_\X$}.
\end{equation*}
Assume that the fibers of $\gtw_\X$ are (non-canonically) isomorphic to a
simple Lie algebra $\g$ over $\C$.  The $\O_\X$-inner product $(.|.)$ is
defined by 
\begin{equation*}
  (A|B)
  := \frac{1}{2\hvee} \tr_{\gtw_\X}(\ad A \ad B)
  \quad
  \text{for $A,B\in\gtw_\X$},
\end{equation*}
where $\hvee$ denotes the dual Coxeter number of $\g$.  Then the inner product
is non-degenerate and invariant under the translations along $\nabla$.  We can
take a local trivialization $\gtw_\X\simeqq\g\tensor\O_\X$, under which the
connection is represented as the exterior derivative on $\X$ (i.e., the
trivial connection).  We assume that we can take such a trivialization of
$\gtw_\X$ on some neighborhood of the divisor $D$.  Under this situation, the
constructions in \secref{sheaf-conf-block} goes through in precisely the same
way and then \lemref{lem:Vir-ActsOngDS-} also holds.

However, \lemref{lem:Vir-ActsOnM} does not always hold.  The action of
$\VirDPr$ on $\M$ is not a representation but a projective representation in
general, because the embedding $\VirDPr\injto\Vir^D_S$ is not always a Lie
algebra homomorphism but so is the composition of the embedding and the
natural projection $\Vir^D_S\onto\Vir^D_S/\O\chat$.  Nevertheless
\lemref{lem:Vir-ActsOnsCC} and \lemref{lem:Vir-ActsOnsCB} also hold in our
case.  Namely, the Lie algebra $\VirDPr$ acts on both the sheaf $\sCC(\M)$ of
conformal coinvariants and the sheaf $\sCB(\M)$ of conformal blocks.
Furthermore \lemref{lem:TDPrKillsCB} can be proved in the same way.  Therefore
we conclude that $\sCB(\M)$ possess a projectively flat connection.  For the
non-twisted Wess-Zumino-Witten models, the coordinate-free version of
\lemref{lem:Ttheta-J} for $\theta\in\TDPr$ is used in the proof of the main
theorem 4.2 in \cite{tsuchimoto:93}.  Since the analogue of
\lemref{lem:Ttheta-J} for $\theta\in\TDPr$ can be also proved in our
situation, we can obtain the projectively flat connections on $\sCC(\M)$.

We can generalize the setting above in various ways: 
\begin{enumerate}
\item We can replace the family of pointed compact Riemann surfaces by that of
  stable pointed curves in the course of \cite{tu-ue-ya:89}, \cite{ueno:95},
  and \cite{tsuchimoto:93}.  Then we shall be able to show the factorization
  property of conformal blocks under appropriate assumptions.
\item We can consider not only deformations of pointed Riemann surfaces but
  also deformations of $\Gtw_\X$-torsors (or principal bundles).  For example,
  the Knizhnik-Zamolodchikov-Bernard equations on Riemann surfaces (cf.\ 
  \cite{ber:88-1}, \cite{ber:88-2}, \cite{fel-wie:96}, \cite{fel:96}) can be
  formulated on a family of pairs of compact Riemann surfaces and principal
  $G$-bundles on them, where $G$ is a finite-dimensional simple algebraic
  group over $\C$.
\item Furthermore we can also consider Borel subgroup bundles which are
  subbundles of the restriction of $\Gtw_\X$ on $D$ and their deformations.
  (More generally, we can consider a family of quasi parabolic structures on
  $\Gtw$-torsors.)  Then we can define the notion of highest weight
  representations of the sheaf of affine Lie algebras with respect to the
  Borel subgroup bundles.  Anyway a choice of a Borel subgroup is required
  by the definition of the category $\O$ of representations, which contains
  the Verma modules, their irreducible quotients, and especially the Wakimoto
  modules.   Note that the constructions of the Wakimoto modules (cf.\ 
  \cite{fei-fre:90-1}, \cite{fei-fre:90-2}, and \cite{kuroki}) essentially
  depend on the choice of a triangular decomposition (equivalently that of a
  Borel subalgebra) of a finite-dimensional semi-simple Lie algebra over $\C$.
\item We can replace the holomorphic flat connection on $\gtw_\X$ by a
  meromorphic flat connection with regular singularity along the divisor $D$.
  Assume that the local monodromy group of the connection around $D$ is
  finite.  Then we can construct a sheaf of twisted affine Lie algebras at $D$
  and can define the notion of conformal blocks for representations of the
  twisted affine Lie algebras.
\end{enumerate}
Detailed expositions shall be given in forthcoming papers.

%
%
%
%
%
%

\end{document}